\input amstex
\magnification=1200
\loadeufm
\loadeusm
\UseAMSsymbols
\input amssym.def
\hsize=7.00 true in
\hoffset=.10 true in
\voffset=-0.1 true in
\vsize=9.50 true in

\font\BIGtitle=cmr10 scaled\magstep3
\font\bigtitle=cmr10 scaled\magstep2
\font\boldsectionfont=cmb10 scaled\magstep1

\def\scr#1{{\fam\eusmfam\relax#1}}
\def\scrA{{\scr A}}
\def\scrB{{\scr B}}
\def\scrC{{\scr C}}

\def\scrD{{\scr D}}
\def\scrE{{\scr E}}
\def\scrF{{\scr F}}
\def\scrG{{\scr G}}
\def\scrH{{\scr H}}
\def\scrI{{\scr I}}
\def\scrK{{\scr K}}
\def\scrJ{{\scr J}}
\def\scrL{{\scr L}}
\def\scrM{{\scr M}}
\def\scrN{{\scr N}}

\def\scrP{{\scr P}}

\def\scrQ{{\scr Q}}
\def\scrS{{\scr S}}
\def\scrU{{\scr U}}
\def\scrV{{\scr V}}
\def\scrW{{\scr W}}

\def\scrX{{\scr X}}
\def\scrY{{\scr Y}}
\def\scrZ{{\scr Z}}
\def\scrR{{\scr R}}
\def\scrT{{\scr T}}
\def\gr#1{{\fam\eufmfam\relax#1}}

\def\grA{{\gr A}}	
	
	\def\grc{{\gr c}}
	\def\grd{{\gr d}}

\def\grG{{\gr G}}	\def\grg{{\gr g}}

\def\grl{{\gr l}}

	\def\grs{{\gr s}}

	\def\grw{{\gr w}}
	\def\grx{{\gr x}}
	\def\gry{{\gr y}}
	\def\grz{{\gr z}}

\def\db#1{{\fam\msbfam\relax#1}}

\def\dbA{{\db A}} 
\def\dbC{{\db C}} \def\dbD{{\db D}}
 \def\dbF{{\db F}}
\def\dbG{{\db G}} \def\dbH{{\db H}}
\def\dbI{{\db I}} 
 \def\dbL{{\db L}}
 \def\dbN{{\db N}}
\def\dbO{{\db O}} 
\def\dbQ{{\db Q}} \def\dbR{{\db R}}
\def\dbS{{\db S}} \def\dbT{{\db T}}
\def\dbU{{\db U}} 
 
 \def\dbZ{{\db Z}}

\def\eps{{\varepsilon}}
\def\Ker{\text{Ker}}

\def\der{\text{der}}

\def\sic{\text{sc}}
\def\Res{\text{Res}}
\def\ab{\text{ab}}
\def\ad{\text{ad}}

\def\Gal{\text{Gal}}
\def\Hom{\text{Hom}}
\def\End{\text{End}}
\def\Aut{\text{Aut}}
\def\Spec{\text{Spec}}
\def\Tr{\text{Tr}}

\def\Lie{\text{Lie}}

\def\leaderfill{\leaders\hbox to 1em
     {\hss.\hss}\hfill}
\def\nspace{\lineskip=1pt\baselineskip=12pt\lineskiplimit=0pt}

\def\finishproclaim{\par\rm
     \ifdim\lastskip<\medskipamount\removelastskip
     \penalty55\medskip\fi}
\def\endproof{$\hfill \square$}

\def\references#1{\par
  \centerline{\boldsectionfont References}%
     \parindent=#1pt\nspace}
\def\Ref[#1]{\par\hang\indent\llap{\hbox to\parindent
     {[#1]\hfil\enspace}}\ignorespaces}
\def\Item#1{\par\hang\indent\llap{\hbox to\parindent
     {#1\hfill$\,\,$}}\ignorespaces}
\def\ItemItem#1{\par\indent\hangindent2\parindent
     \hbox to \parindent{#1\hfill\enspace}\ignorespaces}

\def\arrowsim{\,\smash{\mathop{\to}\limits^{\lower1.5pt
  \hbox{$\scriptstyle\sim$}}}\,}

\def\doublemaprights#1#2#3#4{\raise3pt\hbox{$\mathop{\,\,\hbox to
     #1pt{\rightarrowfill}\kern-30pt\lower3.95pt\hbox to
     #2pt{\rightarrowfill}\,\,}\limits_{#3}^{#4}$}}

\def\rightcapdownarrow{\raise9pt\hbox{$\ssize\cap$}\kern-7.75pt
     \Big\downarrow}

\def\rcapmapdown#1{\rightcapdownarrow\kern-1.0pt\vcenter{
     \hbox{$\scriptstyle#1$}}}

\def\rmapdown#1{\Big\downarrow\kern-1.0pt\vcenter{
     \hbox{$\scriptstyle#1$}}}
\def\rightsubsetarrow#1{{\ssize\subset}\kern-4.5pt\lowe r2.85pt
     \hbox to #1pt{\rightarrowfill}}
\def\longtwoheadedrightarrow#1{\raise2.2pt\hbox to #1pt{\hrulefill}
     \!\!\!\twoheadrightarrow}

\def\Gal{\operatorname{\hbox{Gal}}}
\def\Hom{\operatorname{\hbox{Hom}}}

\def\im{\hbox{Im}}

\NoBlackBoxes
\parindent=25pt
\document
\footline={\hfil}

\bigskip\bigskip
\centerline{\BIGtitle Mod p classification of Shimura F-crystals}
\vskip 0.2in
\centerline{\bigtitle Adrian Vasiu}
\medskip
\centerline{\bigtitle to appear in Math. Nachr.}
\footline={\hfill}
\medskip\noindent
{\bf ABSTRACT}. Let $k$ be an algebraically closed field of positive characteristic $p$. We first classify the $D$-truncations mod $p$ of Shimura $F$-crystals over $k$ and then we study stratifications defined by inner isomorphism classes of these $D$-truncations. This generalizes previous works of Kraft, Ekedahl, Oort, Moonen, and Wedhorn. As a main tool we introduce and study Bruhat $F$-decompositions; they generalize the combined form of Steinberg theorem and of classical Bruhat decompositions for reductive groups over $k$.

\footline={\hss\tenrm \folio\hss}
\pageno=1

\bigskip
\noindent
{\boldsectionfont \S1. Introduction}

\bigskip
A smooth, affine group scheme $E$ over an affine scheme $\Spec(R)$ is called {\it reductive}, if its fibres are reductive groups over fields and thus are connected. Let $E^{\ad}$ be the adjoint group scheme of $E$ i.e., the quotient of $E$ by its center. Let $D$ be a closed subgroup scheme of $E$; if $R$ is not a field, we assume $D$ is smooth. Let $\Lie(D)$ be the {\it Lie algebra} over $R$ of $D$. As $R$-modules, we identify $\Lie(D)=\Ker(D(R[x]/x^2)\to D(R))$, where the $R$-epimorphism $R[x]/(x^2)\twoheadrightarrow R$ maps $x$ to $0$. The Lie bracket on $\Lie(D)$ is obtained by taking the total differential of the commutator morphism $D\times_{\Spec(R)} D\to D$ at identity sections. If $\Lie(D)$ is a free $R$-module, we view it as a $D$-module via the adjoint representation of $D$. For a free $R$-module $K$ of finite rank, let $\pmb{GL}_K$ be the reductive group scheme over $\Spec(R)$ of linear automorphisms of $K$. If $f_1,f_2\in\End_{\dbZ}(K)$, let $f_1f_2:=f_1\circ f_2\in\End_{\dbZ}(K)$.

Let $p$ be a prime. Let $k$ be an algebraically closed field of characteristic $p$. Let $W(k)$ be the ring of Witt vectors with coefficients in $k$. Each reductive group scheme over $\Spec(W(k))$ is split i.e., is a Chevalley group scheme. Let $B(k)$ be the field of fractions of $W(k)$. Let $\sigma$ be the Frobenius automorphism of $k$, $W(k)$, and $B(k)$.

\medskip
\noindent
{\bf 1.1. Shimura $F$-crystals over $k$}

Let $(M,\phi)$ be a {\it Dieudonn\'e module} over $k$. Thus $M$ is a free $W(k)$-module of finite rank and $\phi:M\to M$ is a $\sigma$-linear endomorphism such that we have $pM\subseteq\phi(M)$. Let $\vartheta:M\to M$ be the Verschiebung map of $(M,\phi)$; we have $\phi \vartheta=\vartheta\phi=p1_M$. We denote also by $\phi$ the $\sigma$-linear automorphism of $\End_{B(k)}(M[{1\over p}])$ that maps $e\in\End_{B(k)}(M[{1\over p}])$ to $\phi e\phi^{-1}=\vartheta^{-1} e\vartheta\in\End_{B(k)}(M[{1\over p}])$. Let $G$ be a reductive, closed subgroup scheme of $\pmb{GL}_M$ such that $\phi(\Lie(G_{B(k)}))=\Lie(G_{B(k)})$ and there exists a direct sum decomposition $M=F^1\oplus F^0$ with the properties that: (i) $\phi(M+{1\over p}F^1)=M$ and (ii) the cocharacter of $\pmb{GL}_M$ that fixes $F^0$ and that acts on $F^1$ as the inverse of the identical character of $\dbG_m$, factors through $G$. We denote this factorization as
$$\mu:\dbG_m\to G.$$ 
Triples as $(M,\phi,G)$ show up in the study of special fibres of good integral models of Shimura varieties of Hodge type in mixed characteristic $(0,p)$ (see [19], [16], [35], [39], [40], etc.). We call them {\it Shimura $F$-crystals} over $k$. 

See Subsection 4.4 for the definition of the simple factors of $(\Lie(G^{\ad}),\phi)$ and of their types. Each simple factor of $(\Lie(G^{\ad}),\phi)$ is either trivial or non-trivial. Moreover, each non-trivial simple factor of $(\Lie(G^{\ad}),\phi)$ is of one of the following five types introduced in the context of Shimura varieties of Hodge type in [5]: $A_n$ with $n\ge 1$, $B_n$ with $n\ge 1$, $C_n$ with $n\ge 1$, $D_n^{\dbH}$ with $n\ge 4$, and $D_n^{\dbR}$ with $n\ge 4$.

For $g\in G(W(k))$ let 
$$\scrC_g:=(M,g\phi,G)$$
and let $\bar g\in G(k)$ be $g$ mod $p$. As $\phi^{-1}(pM)=(g\phi)^{-1}(pM)=F^1+pM$, the subgroup 
$$\dbL:=\{h\in G(W(k))|h(F^1+pM)=F^1+pM\}\leqslant G(W(k))$$
is intrinsically associated to the {\it family} $\{\scrC_g|g\in G(W(k))\}$ of Shimura $F$-crystals. If $g_1$, $g_2\in G(W(k))$, then by an {\it inner isomorphism} between $\scrC_{g_1}$ and $\scrC_{g_2}$ we mean an element $h\in G(W(k))$ such that $hg_1\phi=g_2\phi h$; it is easy to see that $h\in\dbL$. The normalizer $P$ of $F^1$ in $G$ is a {\it parabolic} subgroup scheme of $G$, cf. Subsubsection 3.3.4. We have $P(W(k))=\{h\in G(W(k))|h(F^1)=F^1\}\leqslant\dbL=\{h\in G(W(k))|\bar h\in P(k)\}$. 

Let $<\scrC_g>:=\{\scrC_{\dag}|\dag\in G(W(k)),\,\scrC_{\dag}\,\text{is}\,\text{inner}\,\text{isomorphic}\,\text{to}\,\scrC_g\}$ be the inner isomorphism class of $\scrC_g$. Let 
$$\dbI_{\infty}:=\{<\scrC_g>|g\in G(W(k))\}$$
be the set of inner isomorphism classes associated to the family $\{\scrC_g|g\in G(W(k))\}$. 

In [42] we studied the rational classification of $\scrC_g$'s. In this paper we study the mod $p$ classification of $\scrC_g$'s. 

\smallskip
\noindent
{\bf 1.1.1. On $D$-truncations mod $p$}

Let $\bar M:=M/pM$. Let $\bar{\phi}:\bar M\to\bar M$ and $\bar \vartheta:\bar M\to\bar M$ be the reductions mod $p$ of $\phi$ and $\vartheta$ (respectively). By the {\it $D$-truncation mod $p$} of $\scrC_g$ we mean the quadruple 
$$\scrC_{\bar g}:=(\bar M,\bar g\bar{\phi},\bar \vartheta\bar g^{-1},G_k).$$ 
Here $D$ stands for Dieudonn\'e and it is inserted in order not to create confusion with the triple $(\bar M,\bar g\bar{\phi},G_k)$ that is the usual truncation mod $p$ of the triple $(M,g\phi,G)$ which is an $F$-crystal over $k$ with a group. By an inner isomorphism between $\scrC_{\bar g_1}$ and $\scrC_{\bar g_2}$ we mean an element $\bar h\in P(k)$ such that $\bar h\bar g_1\bar{\phi}=\bar g_2\bar{\phi}\bar h$ and $\bar h\bar \vartheta\bar g_1^{-1}=\bar \vartheta\bar g_2^{-1}\bar h$. Let $<\scrC_{\bar g}>:=\{\scrC_{\bar \dag}|\bar \dag\in G(k),\,\scrC_{\bar \dag}\,\text{is}\,\text{inner}\,\text{isomorphic}\,\text{to}\,\scrC_{\bar g}\}$ be the inner isomorphism class of $\scrC_{\bar g}$. Let 
$$\dbI:=\{<\scrC_{\bar g}>|\bar g\in G(k)\}$$
be the set of inner isomorphism classes associated to the family $\{\scrC_{\bar g}|\bar g\in G(k)\}$ of $D$-truncations mod $p$. Let $\text{Aut}^{\text{red}}_{\bar g}$ be the reduced group of the group scheme $\text{Aut}_{\bar g}$ over $k$ of inner automorphisms of $\scrC_{\bar g}$, cf. Subsubsection 5.1.1. We have
$\text{Aut}^{\text{red}}_{\bar g}(k)=\text{Aut}_{\bar g}(k)=\{\bar h\in P(k)|\bar h\bar g\bar{\phi}=\bar g\bar{\phi}\bar h\;\text{and}\;\bar h\bar \vartheta\bar g^{-1}=\bar \vartheta\bar g^{-1}\bar h\}.$

\smallskip
\noindent
{\bf 1.1.2. Weyl elements}

We fix a maximal {\it torus} $T$ of $P$ through which the cocharacter $\mu$ factors (see Subsubsection 3.3.1). It is easy to see that there exists an element $g_{\text{corr}}\in G(W(k))$ such that $g_{\text{corr}}\phi$ normalizes $\Lie(T)$ (see Subsection 4.2). As $\{\scrC_g|g\in G(W(k))\}=\{\scrC_{gg_{\text{corr}}}|g\in G(W(k))\}$, for the study of $\scrC_g$'s and of their $D$-truncations mod $p$ we can replace $\phi$ by $g_{\text{corr}}\phi$. In other words, to ease notations we will assume that $g_{\text{corr}}=1_M$ i.e., we will assume that 
$$\phi(\Lie(T))=\Lie(T).\leqno (1)$$ 
For the sake of functoriality, flexibility, and simplicity of computations, we will {\it not} assume the existence of a Borel subgroup scheme $B$ of $G$ such that $T\leqslant B\leqslant P$ and $\phi(\Lie(B))\subseteq \Lie(B)$.

Let $N_T$ be the normalizer of $T$ in $G$. Let $W_G:=(N_T/T)(W(k))$ be the {\it Weyl group} of $G$. The Weyl group of $P$ is $W_P:=((N_T\cap P)/T)(W(k))\leqslant W_G$. For each $w\in W_G$, we fix a representative $g_w\in N_T(W(k))$ of it.
We note that $<\scrC_{g_w}>$ depends only on $w$ and not on the choice of $g_w$, cf. Lemma 4.1. 

\smallskip
\noindent
{\bf 1.1.3. Basic combinatorial data}

Let $\scrS_{W_P}$ be the set of subgroups of $W_P$. Our basic combinatorial data consists in three maps 
$$\sigma:W_G\arrowsim W_G,\;\;\;\;\;\;\scrS:W_G\to\scrS_{W_P},\;\;\;\text{and}\;\;\;\dbS:W_G\to\dbN\cup\{0\}$$
described as follows. To the quintuple $(M,\phi,G,\mu,T)$ one associates a $\dbZ_p$-structure $(G_{\dbZ_p},T_{\dbZ_p})$ of $(G,T)$, cf. Subsection 4.2. Let $\sigma:W_G\arrowsim W_G$ be the action of $\sigma$ on the Weyl group $W_G$ with respect to this $\dbZ_p$-structure. 

Let $l_0\in\End_{W(k)}(M)$ be the unique semisimple element such that we have $F^0=\{x\in M|l_0(x)=0\}$ and $F^1=\{x\in M|l_0(x)=-x\}$. Let $\bar l_0\in\End_k(\bar M)$ be the reduction mod $p$ of $l_0$. For $w\in W_G$, let 
$$\scrS(w):=\{w_1\in W_G|g_{w_1}\;\;\text{commutes}\;\; \text{with}\;\;(g_w\phi)^i(l_0)\;\;\text{for}\;\;\text{all}\;\;i\in\dbZ\}\leqslant W_P.$$ 
The root decomposition of $\Lie(G)$ relative to $T$ encodes the map $\dbS:W_G\to\dbN\cup\{0\}$, cf. Subsubsection 4.1.1. 

Our main goal is to prove the following three Basic Theorems and to get several applications of them (including the Basic Theorem D of Subsection 12.2). 

\medskip
\noindent
{\bf 1.2. Basic Theorem A (short form)} 

{\it For each element $\bar g\in G(k)$, the identity component of the group $\text{Aut}^{\text{red}}_{\bar g}$ is a unipotent group over $k$ (i.e., has no torus of positive dimension). Moreover for each element $w\in W_G$, we have an equality $\dim(\text{Aut}^{\text{red}}_{\bar g_w})=\dbS(w)$.}   

\medskip
\noindent
{\bf 1.3. Basic Theorem B}

{\it Let $\scrR=\scrR_{\phi,\mu}\subseteq W_G^2$ be the relation such that for $w_1,w_2\in W_G$ we have $(w_1,w_2)\in \scrR$  if and only if there exist elements $w_3\in W_P$ and $w_4\in\scrS(w_2)\leqslant W_P$ with $w_1=w_3w_4w_2\sigma(w_3^{-1})$. Then the following three things hold: 

\medskip
{\bf (a)} Let $w_1, w_2\in W_G$. We have $<\scrC_{g_{w_1}}>=<\scrC_{g_{w_2}}>$ if and only if $<\scrC_{\bar g_{w_1}}>=<\scrC_{\bar g_{w_2}}>$. 

\smallskip
{\bf (b)} 
We have $(w_1,w_2)\in \scrR$ if and only if $<\scrC_{\bar g_{w_1}}>=<\scrC_{\bar g_{w_2}}>$.

\smallskip
{\bf (c)} The relation $\scrR$ on $W_G$ is an equivalence relation.}

\medskip
\noindent
{\bf 1.4. Basic Theorem C}

{\it {\bf (a)} There exists a subset $R_G$ of $W_G$ such that for each element $g\in G(W(k))$, there exist elements $w\in R_G$ and $g_1\in \Ker(G(W(k))\to G(k))$ with the property that $<\scrC_g>=<\scrC_{g_1g_w}>$. 

\smallskip
{\bf (b)} The smallest number of elements a set $R_G$ as in (a) can have is $[W_G:W_P]$.} 

\medskip
Lemma 4.5 (c) checks that we have $<\scrC_{\bar g_1}>=<\scrC_{\bar g_2}>$ if and only if there exists an element $h_{12}\in \Ker(G(W(k))\to G(k))$ such that $<\scrC_{g_1}>=<\scrC_{h_{12}g_2}>$. Thus the statement 1.4 (a) is equivalent to the fact that for each element $\bar g\in G(k)$, there exists $w\in W_G$ such that we have $<\scrC_{\bar g}>=<\scrC_{\bar g_w}>$. Therefore from properties 1.3 (b) and 1.4 (b), we get the following listing of inner isomorphism classes in $\dbI$: 

\smallskip
\noindent
{\bf Corollary 1.1.}
{\it The map $a_{\text{can}}:\scrR\backslash W_G\to\dbI$ that takes $[w]\in \scrR\backslash W_G$ to $<\scrC_{\bar g_w}>\in\dbI$ is a well defined bijection.}$\vfootnote{1}{Strictly speaking, the map $a_{\text{can}}$ is indeed canonical only if one assumes that there exists a Borel subgroup scheme $B$ of $G$ such that $T\leqslant B\leqslant P$ and $\phi(\Lie(B))\subseteq \Lie(B)$. But we repeat that it is very desirable and practical not to make such an assumption.}$
 {\it The sets $\scrR\backslash W_G$ and $\dbI$ have $[W_G:W_P]$ elements.}

\medskip
From Corollary 1.1 and the property 1.3 (a) we get:

\smallskip
\noindent
{\bf Corollary 1.2.} {\it Let $\text{Mod}_p:\dbI_{\infty}\twoheadrightarrow\dbI$ be the surjective map that takes $<\scrC_g>$ to $<\scrC_{\bar g}>$. Then the injective map $\text{Lift}:\dbI\hookrightarrow \dbI_{\infty}$ that takes $<\scrC_{\bar g_w}>$ to $<\scrC_{g_w}>$ is a well defined section of $\text{Mod}_p$.}

\medskip
Corollary 11.1 (d) checks that the section $\text{Lift}:\dbI\hookrightarrow \dbI_{\infty}$ is canonical (i.e., its image does not depend on the maximal torus $T$ of $G$); we call this Teichm\"uller type of section as the {\it toric section} of $\text{Mod}_p$.

\medskip
\noindent
{\bf 1.5. On literature} 

The existence of the set $R_{\pmb{GL}_M}$ is only an adequate translation of the classification of truncated Barsotti--Tate groups of level 1 over $k$ obtained in [17] (see Lemma 4.4 and Corollary 4.6). 
Recent works of Ekedahl--Oort and Moonen extended [17] to the classical context of (principally quasi-polarized) truncated Barsotti--Tate groups of level 1 over $k$ endowed with certain semisimple $\dbF_p$-algebras of endomorphisms (see [26] and [22]; see also [43] for some weaker extensions). These extensions can be used to regain Basic Theorem C for the particular cases related to this classical context; all these cases are such that the simple factors of $(\Lie(G^{\ad}),\phi)$ are either of $A_n$ type or of very particular (the so called totally non-compact) $C_n$ or $D_n^{\dbH}$ type. In [26] no semisimple $\dbF_p$-algebras show up. In [22] the case $p=2$ is almost entirely excluded and the proofs are very long and do not generalize.

Formula $\dim(\text{Aut}^{\text{red}}_{\bar g_w})=\dbS(w)$ of Basic Theorem A is a new (type of) dimension formula for reduced groups of inner automorphisms. It seems to us that this formula is more practical than its analogue of [23] that pertains to the mentioned classical context. For instance, we use this formula to generalize  to the context of the family $\{\scrC_g|g\in G(W(k))\}$, the well known fact that the Barsotti--Tate group of a finite product of supersingular elliptic curves over $k$ is uniquely determined by its truncation of level $1$ (see Theorem 8.3 and Corollary 11.1 (c)). 

If $G=\pmb{GL}_M$, then either [27] or [41, Thm. 1.6] can be used to identify the (small) subset of $\scrR\backslash W_G$ formed by those elements $[w]$ for which the following stronger form of the property 1.3 (a) holds: if $g\in G(W(k))$, then we have $<\scrC_{g_w}>=<\scrC_{g}>$ if and only if $<\scrC_{\bar g_w}>=<\scrC_{\bar g}>$. But even if $G=\pmb{GL}_M$, the property 1.3 (a) does not follow from [27]. We do not know any other relevant literature that pertains to Basic Theorem B. 

\medskip
\noindent
{\bf 1.6. On contents and proofs} 

Section 2 includes complements on algebraic groups over $k$. In particular, it presents a conjecture of us which is a more general form of the combination of: (i) Steinberg theorem for reductive groups over $k$ that are equipped with finite endomorphisms and (ii) Bruhat decompositions with respect to Borel subgroups of reductive groups over $k$. We refer to this conjecture as {\it the Bruhat $F$-decompositions for reductive groups over $k$}; here $F$ stands for a finite (Frobenius) endomorphism. This conjecture is proved in this paper for all those cases that are related to Shimura $F$-crystals, cf. Corollary 7 (a). Well after this paper was submitted, in [12] it is claimed without any details that [12, Cor. 4.1] is an equivalent form of the conjecture (see also [31]). Section 3 recalls properties of reductive group schemes. Section 4 lists basic properties of Shimura $F$-crystals like $\dbZ_p$-structures, types, the existence of the set $R_{\pmb{GL}_M}$, etc. Bruhat $F$-decompositions are naturally associated to certain {\it group actions}. 

In Subsection 5.1 we present our first main new idea that a group action of the form
$$\dbT_{G_k,\sigma}:H_k\times_k G_k\to G_k$$ 
(with $G_k$ viewed only as a variety over $k$), {\it governs the refined structures of $\scrC_{\bar g}$'s}. In particular, if $\dbO$ is the set of orbits of $\dbT_{G_k,\sigma}$, then we have a natural bijection (see Lemma 5.1 (a))
$$b_{\text{can}}:\dbI\arrowsim\dbO.$$ 
The group $H_k$ over $k$ is connected, smooth, and affine; it is not reductive except for the trivial cases when $\mu:\dbG_m\to G$ factors through the center of $G$.
We use the action $\dbT_{G_k,\sigma}$ to prove (see Subsection 5.2) a longer form of the Basic Theorem A. What we do is to ``reinterpret'' the groups $\text{Aut}^{\text{red}}_{\bar g}$ in terms of {\it reduced stabilizer subgroups} of $\dbT_{G_k,\sigma}$ and to prove the analogue of Basic Theorem A for such stabilizer subgroups. 

Subsection 6.1 recalls the classical analogue $\dbT_{G_k}^{\text{cl}}$ of $\dbT_{G_k,\sigma}$ that produces the Bruhat decomposition for the quotient set $P(k)\backslash G(k)$. Theorem 6.1 ``compares'' $\dbT_{G_k}^{\text{cl}}$ and $\dbT_{G_k,\sigma}$ to prove (via Basic Theorem A) that if the set $\dbO$ is finite, then $\dbO$ has precisely $[W_G:W_P]$ elements and the values of the function $\dbS$ are computable in terms of values of the classical analogue $\dbS^{\text{cl}}:W_G\to\dbN\cup\{0\}$ of $\dbS$. We emphasize that Sections 5 and 6 (and thus also Basic Theorem A) are elementary in nature i.e., they rely mainly on basic properties of group actions.

In Section 7 we present our second main new idea: we introduce and study {\it the zero space} $\grz_{\bar g}$ of $\scrC_{\bar g}$. This zero space $\grz_{\bar g}$ is an $\dbF_p$-Lie subalgebra of $\Lie(P_k)$. The role of $\grz_{\bar g}$ is triple folded: 

\medskip
{\bf (a)}  It is the ``$\sigma$-Lie algebra'' of $\text{Aut}_{\bar g}$; thus $\text{Aut}^{\text{red}}_{\bar g}$ normalizes $\grz_{\bar g}$ (cf. Lemma 7.1 (c)). 

\smallskip
{\bf (b)} A suitable $\dbF_p$-Lie subalgebra $\grx_{\bar g}$ of $\grz_{\bar g}$ is the Lie algebra of an $\dbF_p$-structure of a (uniquely determined) reductive subgroup of $G_k$ of the same rank as $G_k$. See Theorem 7.6 for $\grx_{\bar g_w}=\{x\in \grz_{\bar g_w}|[\bar l_0,x]=0\}$. 

\smallskip
{\bf (c)} The $k$-span $\grw_{\bar g}$ of $\grz_{\bar g}$ records all ``good'' $P(k)$-conjugates of the cocharacter $\mu_k:\dbG_m\to G_k$ (i.e., of the element $\bar l_0\in\Lie(P_k)$) with respect to $\scrC_{\bar g}$. For instance, if $w\in W_G$, then Theorem 9.1 shows that a $P(k)$-conjugate $\bar l_1$ of $\bar l_0$ is also an $\text{Aut}^{\text{red}}_{\bar g_w}$-conjugate of $\bar l_0$ if and only if $\bar l_1\in\grw_{\bar g_w}$. Moreover, $\grw_{\bar g}$ is the Lie algebra of a (uniquely determined) connected, smooth subgroup of $G_k$ of the same rank as $G_k$ (see Theorem 7.8 for $\grw_{\bar g_w}$).

\medskip
Theorem 9.1 relies on Basic Theorem A and Section 7 and it is the main ingredient in proving Basic Theorem B in Subsection 9.1.  

Section 10 combines Sections 4 and 7 and Theorem 9.1 to show that under some conditions, the set $R_G$ exists if its analogue $R_{G_1}$ exists, where $G_1$ is a reductive, closed subgroup scheme of $\pmb{GL}_M$ that contains $G$ and such that the triple $(M,\phi,G_1)$ is a Shimura $F$-crystal over $k$ (see Theorem 10.1). This standard idea is also used in [22] and [43] which worked with $G_1=\pmb{GL}_M$. As a main difference from [22] and [43], in the proof of Theorem 10.1 we use direct sum decompositions $\Lie(G_1)=\Lie(G)\oplus\Lie(G)^\perp$ of $G$-modules or adjoint variants of them and we do not appeal to either $\dbZ_p$- or $\dbF_p$-algebras or the canonical stratifications of $(\bar M,\bar g\bar\phi,\bar \vartheta\bar g^{-1})$'s introduced in [17] and [26]. Moreover, we use (reductions mod $p$) of Verschiebung maps only in Subsection 4.3 and (to introduce notations) in Subsection 5.1.1 and Section 12. See Subsection 3.2 for basic properties of trace maps (most common such direct sum decompositions $\Lie(G_1)=\Lie(G)\oplus\Lie(G)^\perp$ are produced by them). 

Let $P_1$ be the normalizer of $F^1$ in $G_1$ and let $g\in G(W(k))$. The guiding idea of the proof of Theorem 10.1 is: as the set $R_{G_1}$ exists, there exists an element $\bar l_1\in \Lie(G_{1k})$ that is defined by a ``good'' $P_1(k)$-conjugate of $\mu_k:\dbG_m\to G_k$ with respect to $(\bar M,\bar g\bar\phi,\bar\vartheta\bar g^{-1},G_{1k})$; we use the mentioned direct sum decompositions to define a component of $\bar l_1$ in $\Lie(G_k)$. This component is defined by a ``good'' $P(k)$-conjugate of $\mu_k:\dbG_m\to G_k$ with respect to $\scrC_{\bar g}$ and thus we can exploit (b) to show that $<\scrC_{\bar g}>$ is $<\scrC_{\bar g_w}>$ for some $w\in W_G$. 

Section 11 uses Sections 4 to 7, 9, and 10 to prove Basic Theorem C. The main four steps are as follows.

\medskip
{\bf (i)} The existence of the set $R_G$ is of intrinsic nature i.e., is encoded in the adjoint analogue $\dbT_{G_k^{\ad},\sigma}$ of the action $\dbT_{G_k,\sigma}$. Thus to prove Basic Theorem C, we can assume that the pair $(\Lie(G^{\ad}),\phi)$ is simple and we can appeal to the {\it replacement process} which allows us to replace $(M,\phi,G)$ by any other Shimura $F$-crystal $(M_*,\phi_*,G_*)$ over $k$ for which we have an identity $(\Lie(G^{\ad}_*),\phi_*)=(\Lie(G^{\ad}),\phi)$; see properties 11.1 (i) and (ii).

\smallskip
{\bf (ii)}  The statement 1.4 (a) implies the statement 1.4 (b), cf. Corollary 6.2.

\smallskip
{\bf (iii)} Suppose we are not dealing with an exceptional case. Based on the replacement process and on the classification of cocharacters of the form $\mu:\dbG_m\to G\hookrightarrow \pmb{GL}_M$ (see [30] over $B(k)$), we can assume that there exists a reductive, closed subgroup scheme $G_1$ of $\pmb{GL}_M$ such that the triple $(M,\phi,G_1)$ is a Shimura $F$-crystal and the hypotheses of Theorem 10.1 hold (in particular, the set $R_{G_1}$ exists). Thus based on Theorem 10.1, we {\it inductively} prove that the statement 1.4 (a) holds outside the exceptional cases. 

\smallskip
{\bf (iv)} The exceptional cases are when $p=2$ and $(\Lie(G^{\ad}),\phi)$ is of $B_n$, $C_n$, $D_n^{\dbH}$, or $D_n^{\dbR}$ type. To prove the statement 1.4 (a) for these cases, we use a {\it shifting process} that is a standard combinatorial application of Basic Theorem B and Theorem 6.2. This shifting process says that it suffices to prove Basic Theorem C for $p>>0$.  

\medskip
\noindent
{\bf 1.7. Motivations and extra applications} 

Our main reason in proving Basic Theorems A to C is to get tools to study: (i) the set $\dbI_{\infty}$ and (ii) different {\it stratifications} of special fibres of integral canonical models of Shimura varieties of Hodge type in mixed characteristic $(0,p)$ that are proved to exist in [35] and [39] or are (assumed to exist and) considered in [40, Sect. 8] and [42, Sect. 5]. For instance, based on Subsections 8.1, 12.2, 12.3, and 12.4 (a), the methods of [25] can be adapted to show that each rational stratification of [42, Subsect. 5.3] has a unique closed stratum. Based on either [35, Subsect. 5.4] or [42, Sect. 5], the generalization of [26, Sect. 1] and [43, Thm. of Introd.] to the mentioned special fibres is automatic (see  Remark 12.4 (a)). Thus for the sake of generality and of not making the paper too long by recalling Shimura varieties of Hodge type, in Section 12 we introduce stratifications that generalize loci cit. directly in an axiomatized context. For the context of reductive groups, our stratifications are more general than the level $1$ stratifications formalized in [38, Subsubsect. 4.2.3]. Example 12.2 and Remark 12.4 (a) explain how the axiomatized context applies naturally to special fibres as above. Subsections 12.2 and 12.3 contain many tools one needs in applications (the number and dimensions of strata, specializations between strata, etc.). 

\medskip
\noindent
{\bf 1.8. Generalizations} 

Let $a$, $b\in\dbZ$ with $b\ge a$. The methods of Sections 5 to 7 and 9 to 12  work even if we replace the roles of the quadruple $(M=F^1\oplus F^0,\phi,\mu:\dbG_m\to G,\dbL)$ and of $D$-truncations by the ones of an analogue quadruple 
$$(\tilde M=\oplus_{i=a}^b \tilde F^i,\tilde\phi,\tilde\mu:\dbG_m\to \tilde G,\widetilde{\dbL}:=\{\tilde g\in\tilde G(W(k))|\tilde g(\oplus_{i=a}^b p^{-i}\tilde F^i)=\oplus_{i=a}^b p^{-i}\tilde F^i\})\leqno (2)$$ 
and of $F$-truncations defined in [37, Subsubsect. 3.2.9]. The parts of Sections 10 to 12 that rely on the decomposition $M=F^1\oplus F^0$, extend to the context of (2) only if $b-a\le p-1$; but a similar shifting process shows that to generalize Basic Theorem C to the context of (2) we can assume $p>>0$. Thus to generalize Basic Theorems A to C to the context of (2) one only needs to: (i) prove based on Remark 9.5 (a) the analogue of Basic Theorem C for the case when $\tilde G=\pmb{GL}_{\tilde M}$ and (ii) apply the methods of Sections 4 to 7 and 9 to 12. 

\bigskip
\noindent
{\boldsectionfont \S2. Complements on algebraic groups over $k$}

\bigskip
We include complements on a connected, smooth, affine  group $\scrG$ over $k$. We recall that two maximal tori of $\scrG$ are $\scrG(k)$-conjugate, cf. [1, Ch. IV, Cor. 11.3 (1)].

\medskip
\noindent
{\bf 2.1. Two classical results} 

Let $F:\scrG\to\scrG$ be a finite endomorphism. Let $\scrG^F(k)$ be the subgroup of $\scrG(k)$ formed by elements fixed by $F$. Steinberg theorem asserts (cf. [32, Thm. 10.1]): if the group $\scrG^F(k)$ is finite, then $\scrG(k)=\{hF(h)^{-1}|h\in\scrG(k)\}$. 

Often, $\scrG$ has a model $\scrG_{\dbF_{p^m}}$ over a finite field $\dbF_{p^m}$ and $F:\scrG\to\scrG$ is the Frobenius endomorphism with respect to it; the group $\scrG^F(k)=\scrG_{\dbF_{p^m}}(\dbF_{p^m})$ is finite and thus we have $\scrG(k)=\{hF(h)^{-1}|h\in\scrG(k)\}$. This last identity is called Lang theorem and is equivalent to the fact that the set $H^1(\dbF_{p^m},\scrG_{\dbF_{p^m}})$ has only one class.  

From now on we will take $\scrG$ to be reductive. Let $\scrT$ be a maximal torus of a Borel subgroup $\scrB$ of $\scrG$. Let $\scrU^{\scrB}$ be the unipotent radical of $\scrB$. Let $W_{\scrG}$ be the Weyl group of $\scrG$ with respect to $\scrT$. For $w\in\scrW_{\scrG}$, let $h_w\in\scrG(k)$ be an element that normalizes $\scrT$ and that represents $w$. The classical Bruhat decomposition for the pair $(\scrG,\scrB)$ asserts that we have a disjoint union $\scrG(k)=\sqcup_{w\in W_{\scrG}} \scrB(k)h_w\scrU^{\scrB}(k)$, cf. [1, Ch. IV,  Subsect. 14.12].

\medskip
\noindent
{\bf 2.2. Bruhat $F$-decompositions} 

We recall that $\scrG$ is reductive. 
Let $F:\scrG\to \scrG$ be a finite endomorphism. Let $\scrG^F(k),\scrT, W_{\scrG}$, and $(h_w)_{w\in W_{\scrG}}$ be as in Subsection 2.1. Let $\scrP$ be a parabolic subgroup of $\scrG$ that contains $\scrT$. Let $\scrU$ be the unipotent radical of $\scrP$. Let $\scrL$ be the unique Levi subgroup of $\scrP$ that contains $\scrT$, cf. [6, Vol. III, Exp. XXVI, Prop. 1.12 (ii)]. Let $\scrQ$ be the parabolic subgroup of $\scrG$ that is the opposite of $\scrP$ with respect to $\scrT$ (i.e., we have $\scrQ\cap\scrP=\scrL$). Let $\scrN$ be the unipotent radical of $\scrQ$. In Spring 2001 we conjectured:

\smallskip
\noindent
{\bf Conjecture 2.1 (Bruhat $F$-decompositions for reductive groups over $k$).} {\it Suppose the group $\scrG^F(k)$ is finite. If $\scrP\neq\scrG$, we also assume that $F(\scrT)=\scrT$. Then each element $h\in\scrG(k)$ can be written as a product $h=h_1h_2h_wF(h_2^{-1})F(h_3^{-1})$,
where $(h_1,h_2,h_3)\in \scrU(k)\times\scrL(k)\times\scrN(k)$ and where $w\in W_{\scrG}$. The minimal number of such elements $w\in W_{\scrG}$  we need is precisely $[W_{\scrG}:W_{\scrP}]$, where $W_{\scrP}:=\{w\in W_{\scrG}|h_w\scrP h_w^{-1}=\scrP\}$.}
 
\medskip
Since 2001, there have been three approaches to the Conjecture 2.1.

\medskip
$\bullet$ The inductive approach presented in this paper and which is carried on in Corollary 7 (a) for those cases that are related to Shimura $F$-crystals over $k$.

$\bullet$ The elementary approach of [31] based on group actions (that works in very particular cases). 

$\bullet$ The approach based on partial flag varieties associated to $\scrG$ presented in [12] (an equivalent form of Conjecture 2.1 is claimed in [12, Cor. 4.1]; but the passage from loc. cit. to Conjecture 2.1 requires substantial work which is not presented in [12]).

\medskip
For the sake of being self-contained, of providing more algebraic geometric 
 background to Shimura $F$-crystals, and of presenting the original approach to Conjecture 2.1, in this paper we do not rely on [12] or [31].

\smallskip
\noindent
{\bf 2.2.1. The $\dbT$-actions} 

Behind Conjecture 2.1 there exists a group action $$\dbT_{\scrG,\scrP,\scrT,F}:\scrH\times_k \scrG\to\scrG$$ 
on the $k$-variety $\scrG$ defined as follows. Let $\scrH:=\scrU\times_k \scrL\times_k \scrN$ be endowed with a new group structure via the rules: if $(r_1,r_2,r_3)$, $(s_1,s_2,s_3)\in\scrH(k)$, then
$$(r_1,r_2,r_3)\cdot (s_1,s_2,s_3):=(r_1r_2s_1r_2^{-1},r_2s_2,r_3r_2s_3r_2^{-1})\in \scrH(k)\leqno (3)$$
and $(r_1,r_2,r_3)^{-1}:=(r_2^{-1}r_1^{-1}r_2,r_2^{-1},r_2^{-1}r_3^{-1}r_2)\in \scrH(k)$.

If $\scrP\neq\scrG$, then this new group structure (on $\scrH$) is truly new i.e., it is not the natural product group structure (on $\scrH$). We call this new group structure (on $\scrH$) the {\it semidirect product group structure}. Note that $\scrH$ is isomorphic to the semidirect product $(\scrU\times_k \scrN)\rtimes_k \scrL$ defined by the inner conjugation action of $\scrL$ on $\scrU\times_k \scrN$. Thus, if $\scrH\twoheadrightarrow \scrL$ is the epimorphism that maps $(r_1,r_2,r_3)\in\scrH(k)$ to $r_2\in \scrL(k)$, then we have a short exact sequence
$$1\to \scrU\times_k \scrN\to \scrH\to \scrL\to 1.\leqno (4)$$ The new group structure on $\scrU\times_k \scrL$ (resp. on $\scrL\times_k\scrN$) is isomorphic to $\scrP$ (resp. to $\scrQ$) via the isomorphism which on $k$-valued points maps $(r_1,r_2)\in \scrU(k)\times\scrL(k)$ to $r_1r_2\in\scrP(k)$ (resp. maps $(r_2,r_3)\in \scrL(k)\times\scrN(k)$ to $r_3r_2\in\scrQ(k)$). The product $k$-morphism $\scrU\times_k\scrL\times_k \scrN\to \scrG$ is an open embedding around the identity elements, cf. [6, Vol. III, Exp. XXII, Prop. 4.1.2]. Thus $\dim(\scrH)=\dim(\scrU)+\dim(\scrL)+\dim(\scrN)=\dim(\scrG)$. 

We define the $\dbT$-action $\dbT_{\scrG,\scrP,\scrT,F}$ via the following rule on $k$-valued points $$\dbT_{\scrG,\scrP,\scrT,F}((r_1,r_2,r_3),h):=r_1r_2hF(r_2^{-1})F(r_3^{-1}).\leqno (5)$$ 
Each orbit of $\dbT_{\scrG,\scrP,\scrT,F}$ is a connected, smooth, locally closed subscheme of $\scrG$, cf. [1, Ch. I, Subsect. 1.8, Prop.] and the fact that $\scrH$ is connected. Our first motivations for Conjecture 2.1 are based on the following Proposition. 

\smallskip
\noindent
{\bf Proposition 2.2.} {\it Suppose the group $\scrG^F(k)$ is finite. If $\scrP\neq\scrG$, we also assume that $F(\scrT)=\scrT$. We have:

\medskip
{\bf (a)} If $h_0\in \scrG(k)$, then the group $\scrG^{h_0Fh_0^{-1}}(k)=\{h\in\scrG(k)|h_0F(h)h_0^{-1}=h\}$ is finite. 

\smallskip
{\bf (b)} Conjecture 2.1 holds if and only if the action $\dbT_{\scrG,\scrP,\scrT,F}$ has precisely $[W_{\scrG}:W_{\scrP}]$ orbits and each such orbit contains $h_w$ for some element $w\in W_{\scrG}$. For each $w\in W_{\scrG}$ the orbit of $h_w$ under $\dbT_{\scrG,\scrP,\scrT,F}$ depends only on $w\in W_G$ (i.e., it contains $\scrT(k)h_w$).

\smallskip
{\bf (c)} If $\scrP=\scrG$, then Conjecture 2.1 holds.

\smallskip
{\bf (d)} Let $\scrP'$ and $\scrT'$ be the images of $\scrP$ and $\scrT$ (respectively) in $\scrG^{\ad}$. Let $F':\scrG^{\ad}\to\scrG^{\ad}$ be the endomorphism induced naturally by $F$. Then Conjecture 2.1 holds for the quadruple $(\scrG,\scrP,\scrT,F)$ if and only if it holds for the quadruple $(\scrG^{\ad},\scrP',\scrT',F')$.

\smallskip
{\bf (e)} If $h_0\in \scrG(k)$ normalizes $\scrT$, then Conjecture 2.1 holds for the finite endomorphism $F:\scrG\to\scrG$ if and only if it holds for the finite endomorphism $h_0Fh_0^{-1}:\scrG\to\scrG$. 

\smallskip
{\bf (f)} If $\scrP$ is a Borel subgroup $\scrB$ of $\scrG$, then Conjecture 2.1 holds.}
\medskip

{\bf Proof:}
As we have a central isogeny from $\scrG$ to $\scrG^{\ad}$ times the abelianization of $\scrG$, it suffices to prove (a) when $\scrG$ is either a torus or an adjoint group. The case when $\scrG$ is a torus is trivial. We now assume $\scrG=\scrG^{\ad}$. As $F$ maps a simple factor of $\scrG$ to a simple factor of $\scrG$, (by replacing $F$ with positive integral powers of itself) it suffices to check (a) under the extra hypothesis that $\scrG$ is simple. As $\scrG^F(k)$ is finite, $F$ is not an automorphism (cf. [32, Cor. 10.13]). Thus $h_0Fh_0^{-1}$ is not an automorphism. Thus the group $\scrG^{h_0Fh_0^{-1}}(k)$ is finite (cf. loc. cit.) i.e., (a) holds. 

We prove (b) and (c). The first part of (b) is obvious. If $\scrP=\scrG$, then the second part of (b) and (c) are only a reformulation of the first paragraph of Subsection 2.1 (for $\scrG$ reductive). Thus (c) holds. We check the second part of (b) for the case when $\scrP\neq\scrG$. For $\tilde h_2\in \scrT(k)$ there exists $h_2\in \scrT(k)$ such that $\tilde h_2=h_2 h_wF(h_2^{-1})h_w^{-1}$, cf. (c) applied to the finite endomorphism of $\scrT$ defined by $h_wFh_w^{-1}$. Thus $\scrT(k)h_w$ belongs to the orbit of $h_w$ under $\dbT_{\scrG,\scrP,\scrT,F}$. Therefore (b) holds.

We prove (d). The action $\dbT_{\scrG^{\ad},\scrP',\scrP',F'}$ is naturally a quotient of the action $\dbT_{\scrG,\scrP,\scrP,F}$. Thus under the epimorphism $\scrG\twoheadrightarrow\scrG^{\ad}$, each orbit of $\dbT_{\scrG,\scrP,\scrP,F}$ surjects onto an orbit of $\dbT_{\scrG^{\ad},\scrP',\scrP',F'}$. Due to this and the fact that the center of $\scrG$ is contained in $\scrT$, part (d) follows from (b).

We prove (e). The orbit of $h_w$ under $\dbT_{\scrG,\scrP,\scrT,h_0Fh_0^{-1}}$ is the right translation by $h_0^{-1}$ of the orbit of $h_wh_0$ under $\dbT_{\scrG,\scrP,\scrT,F}$. From this (e) follows, cf. (b). 

We prove (f). We have $\scrL=\scrT$. By replacing $F$ with some $h_0Fh_0^{-1}$, we can assume $F(\scrN)=\scrU$ (cf. (e)). Let $w\in W_\scrG$. Due to the second part of (b), the orbit of $h_w$ under $\dbT_{\scrG,\scrP,\scrT,F}$ is $\scrU\scrT h_w\scrU=\scrP h_w\scrU$. Thus (f) is equivalent to the existence of the Bruhat decomposition of Subsection 2.1 applied with $(\scrB,\scrU^{\scrB})=(\scrP,\scrU)$.\endproof

\smallskip
\noindent
{\bf Lemma 2.3.} {\it Suppose $\scrG$ is a simple, adjoint group over $k$. 

\medskip
{\bf (a)} The center of $\Lie(\scrG)$ is trivial.

\smallskip
{\bf (b)} Let $\scrP$ be a maximal parabolic subgroup of $\scrG$ such that its unipotent radical $\scrU$ is commutative. Then the centralizer $\grc$ of $\Lie(\scrU)$ in $\Lie(\scrP)$ is $\Lie(\scrU)$.}
\medskip

{\bf Proof:}
As $\scrG$ is adjoint, its center is trivial. Thus (a) holds. We check (b). As $\scrU$ is commutative, it centralizes $\Lie(\scrU)$. Thus $\Lie(\scrU)\subseteq \grc$. To end the proof of (b) it suffices to show that the assumption that $\Lie(\scrU)\subsetneq \grc$ leads to a contradiction. This assumption implies that $\grc$ contains a non-zero element of $\Lie(\scrL)$, where $\scrL$ is a Levi subgroup of $\scrP$. But the representation $\scrL\to \pmb{GL}_{\Lie(\scrU)}$ induced by the inner conjugation action of $\scrP$ on $\Lie(\scrU)$ is a closed embedding, cf. [36, Prop. 5.1]. This implies that $\grc$ contains no element of $\Lie(\scrL)$. Contradiction.\endproof 

\smallskip
\noindent
{\bf Proposition 2.4.} {\it 
Let $\bar M$ be a finite dimensional vector space over $k$. We assume $\scrG$ is a reductive subgroup of $\pmb{GL}_{\bar M}$. Let $S_0$ be a connected, smooth subgroup of $\scrG$ that contains a maximal torus of $\scrG$. Let $S_1$ and $S_2$ be tori of $\pmb{GL}_{\bar M}$ generated by images of cocharacters of $\pmb{GL}_{\bar M}$ of weights $\{-1,0\}$. 

\medskip
{\bf (a)} The tori $S_1$ and $S_2$ commute if and only if $[\Lie(S_1),\Lie(S_2)]=0$.

\smallskip
{\bf (b)} Suppose $\Lie(S_1)\subseteq\Lie(S_0)$ and $S_2$ is a maximal torus of $\scrG$. Then $S_1$ centralizes a maximal torus of $S_0$. If moreover $p>2$, then $S_1$ normalizes $\Lie(S_0)$.}
\medskip

{\bf Proof:}
We prove (a). The only if part is trivial. To check the if part we can assume $S_1$ and $S_2$ are of rank 1. But this case is trivial. Thus (a) holds.

We prove the first part of (b). We will consider the following process. If $\Lie(S_1)$ is not contained in the center of $\Lie(S_0)$, let $x_0\in\Lie(S_1)$ be such that it does not centralize $\Lie(S_0)$. Let $C_0$ be the identity component of the centralizer of $x_0$ in $S_0$. As $x_0$ is a semisimple element of $\Lie(S_0)$, the group scheme $C_0$ is smooth and $\Lie(C_0)$ is the centralizer of $x_0$ in $\Lie(S_0)$ (cf. [1, Ch. III, Subsect. 9.1]). Thus $\Lie(S_1)\subseteq\Lie(C_0)$. As $x_0$ is contained in the Lie algebra of a torus of $S_0$ (cf. [1, Ch. IV, Prop. 11.8]), the group $C_0$ contains a maximal torus of $S_0$. As the maximal tori of $S_0$ are $S_0(k)$-conjugate, $C_0$ also contains a maximal torus of $\scrG$. If $\Lie(S_1)$ is not contained in the center of $\Lie(C_0)$, we repeat this process in the context of the pair $(S_1,C_0)$ instead of the pair $(S_1,S_0)$. 

After a finite number of steps that perform this process, we reach a context in which $\Lie(S_1)$ is contained in the center of the Lie algebra of a connected, smooth subgroup of $S_0$ that contains a maximal torus $\tilde S_2$ of $\scrG$. As $\tilde S_2$ and $S_2$ are $\scrG(k)$-conjugate, $\tilde S_2$ is also generated by images of cocharacters of $\pmb{GL}_{\bar M}$ of weights $\{-1,0\}$. As $[\Lie(S_1),\Lie(\tilde S_2)]=0$, from (a) we get that $S_1$ centralizes the maximal torus $\tilde S_2$ of $S_0$.    

To check the last part of (b) we can assume $S_1$ has rank 1. Let $\bar l_{S_1}\in\Lie(S_1)$ be the unique non-zero semisimple element whose eigenvalues (as an endomorphism of $\bar M$) are $0$ and $-1$. Let $\Lie(S_0)=\grc_{-1}\oplus\grc_0\oplus\grc_1$ be the direct sum decomposition such that for $i\in\{-1,0,1\}$ and $x\in \grc_i$ we have $[\bar l_{S_1},x]=ix$. This makes sense as $p>2$. Now it is easy to see that $S_1$ normalizes $\grc_i$ for all $i\in\{-1,0,1\}$. Thus $S_1$ normalizes $\Lie(S_0)$. \endproof

\smallskip
\noindent
{\bf Lemma 2.5.} {\it 
Suppose $k$ is an algebraic closure of a finite field $\kappa$ of characteristic $p$. Let $\scrZ$ be a connected, smooth group over $\kappa$. Let $\scrV$ be a smooth subgroup of $\scrZ$ whose identity component $\scrV^0$ is unipotent. Then the quotient $\kappa$-varieties $\scrX:=\scrZ/\scrV$ and $\scrX_0:=\scrZ/\scrV^0$ have the same number of $\kappa$-valued points.}
\medskip

{\bf Proof:} If $\scrZ$ is affine and $\scrV/\scrV^0$ is a constant \'etale group in which $\scrV(\kappa)$ surjects (this is the case needed in this paper), then one can check this Lemma via a direct count of $\kappa$-valued points. 

The proof below is due to Serre and Gabber (independently). Let $l$ be a prime different from $p$. Let $\text{H}^{*}(\ddag)$ (resp. $\text{H}_{c}^{*}(\ddag)$) be the $l$-adic cohomology (resp. the $l$-adic cohomology with proper support) with coefficients in $\dbQ_l$ of the extension to $k$ of the $\kappa$-variety $\ddag$. We claim that the $\dbQ_l$-linear map  $\text{H}_c^{*}(\scrX)\to \text{H}_c^{*}(\scrX_0)$  associated to the quotient morphism $\scrX_0\to \scrX$, is an isomorphism. If we know this, then $\scrX(\kappa)$ and  $\scrX_0(\kappa)$  have the same number of elements equal to the Lefschetz number of the Frobenius acting on the cohomology with compact support. 

As $\scrX_k$ is the quotient of $\scrX_{0k}$ under the right action of $C:=\scrV(k)/\scrV_0(k)$ on $\scrX_{0k}$, the mentioned $\dbQ_l$-linear map $\text{H}^{*}(\scrX)\to \text{H}^{*}(\scrX_0)$ is injective and its image is made up of the
elements of $\text{H}^{*}(\scrX_0)$ fixed by $C$; the same holds for $\text{H}^{*}_c$. This is a very special case of the Hochschild--Serre spectral
sequence, cf. [44, Exp. VIII, p. 406]. 

As $\scrV_0$ is unipotent, it has a normal series whose quotient factors are $\dbG_a$ groups (cf. [6, Vol. II, Exp. XVII, Cor. 4.1.3]). Thus the geometric fibres of the quotient morphism $\scrZ\to\scrX_0$ are isomorphic to $\dbA^d$, where $d:=\dim(\scrV_0)$. The $\dbQ_l$-linear map $\text{H}^{*}(\scrX_0)\to \text{H}^{*}(\scrZ)$ associated to $\scrZ\to\scrX_0$ is an isomorphism and the same holds for  $\text{H}^{*}_c$ with a Tate-twist by $(d)$.
  This is a general fact on smooth, surjective, separated, finite type morphisms whose geometric fibres are affine spaces and thus are acyclic for $l$ (see [44, Exp. XVI, Cor. 1.16 and Thm. 2.1] and [44, Exp. XVII, Cor. 1.2] for the case of $\text{H}^{*}$). 

Based on the last two paragraphs, to prove the claim it suffices to show that the natural right action of $\scrZ(k)$ on $\text{H}^{*}(\scrZ)$ is trivial. But as $\scrZ$ is connected, this is implied by [7, Cor. 6.5].\endproof

\bigskip
\noindent
{\boldsectionfont \S3. Elementary properties of reductive group schemes}

\bigskip
Let $R$, $E$, $D$, and $K$ be as before Subsection 1.1. Let $E^{\der}$, $Z(E)$, and $E^{\ad}$ be the derived group scheme of $E$, the center of $E$, and the adjoint group scheme of $E$ (respectively). Thus we have $E^{\ad}=E/Z(E)$ and $E^{\ab}:=E/E^{\der}$ is the abelianization of $E$, cf. [6, Vol. III, Exp. XXII, Def. 4.3.6 and Thm. 6.2.1]. Let $E^{\sic}$ be the simply connected semisimple group scheme cover of $E^{\der}$. The group scheme $Z(E)$ is of multiplicative type, cf. [6, Vol. III, Exp. XXII, Cor. 4.1.7]. Let $Z^0(E)$ be the maximal torus of $Z(E)$; the quotient group scheme $Z(E)/Z^0(E)$ is finite, flat and of multiplicative type. If $\tilde\mu:\dbG_m\to E$ is a cocharacter of $E$, let $d\tilde\mu:\Lie(\dbG_m)\to\Lie(E)$ be the Lie homomorphism defined by $\tilde\mu$. A cocharacter $\mu_K:\dbG_m\to\pmb{GL}_K$ is said to have weights $\{-1,0\}$ if it produces a direct sum decomposition $K=K_0\oplus K_1$ such that for $i\in\{0,1\}$ we have $K_i\neq 0$ and $\dbG_m$ acts through $\mu_K$ on $K_i$ via the $-i$-th power of the identical character of $\dbG_m$. 

A parabolic subgroup scheme of $E$ is a smooth, closed subgroup scheme $P_E$ of $E$ whose fibres are parabolic subgroups over fields. The unipotent radical $U_E$ of $P_E$ is the maximal unipotent, smooth, normal, closed subgroup scheme of $P_E$; it has connected fibres. The quotient group scheme $P_E/U_E$ exists, cf. [6, Vol. III, Exp. XXVI, Prop. 1.6]. We identify $U_E$ with the closed subgroup scheme of $E^{\ad}$ which is the unipotent radical of the parabolic subgroup scheme of $E^{\ad}$ that is the image of $P_E$ in $E^{\ad}$. A Levi subgroup scheme $L_E$ of $P_E$ is a reductive, closed subgroup scheme of $P_E$ such that the natural homomorphism $L_E\to P_E/U_E$ is an isomorphism. 

For a finite, flat monomorphism $R_0\hookrightarrow R$, let $\Res_{R/R_0} D$ be the group scheme over $\Spec(R_0)$ obtained from $D$ through the Weil restriction of scalars (see [4, Subsect. 1.5] and [3, Ch. 7, Subsect. 7.6]). If $R_1$ is a commutative $R_0$-algebra, we have a canonical and functorial group identification $\Res_{R/R_0} D(R_1)=D(R_1\otimes_{R_0} R)$. The pull back  via $\Spec(R)\to\Spec(R_0)$ of an object or a morphism $Y$ or $Y_{R_0}$ (resp. $Y_*$ with $*$ an arbitrary index different from $R_0$) of the category of $\Spec(R_0)$-schemes, is denoted by $Y_R$ (resp. by $Y_{*R}$). 

Until Subsection 3.3 we recall known properties of $D$, $E$, and trace maps. Subsection 3.3 specializes to the concrete context of the reductive, closed subgroup scheme $G$ of $\pmb{GL}_M$ and we list basic notations and properties. 

See [24, Prop. 2.1] for a proof of the following direct application of Lang theorem. 

\smallskip
\noindent
{\bf Lemma 3.1.} {\it
Let $D$ be a smooth group scheme over $\Spec(\dbZ_p)$ that has a connected,  affine special fibre $D_{\dbF_p}$. Let the automorphism $\sigma$ of $W(k)$ act on $D(W(k))$ in the natural way. Then for each element $h\in D(W(k))$, there exists an element $\tilde h\in D(W(k))$ such that we have $h=\tilde h\sigma(\tilde h^{-1})$.}
\medskip

\smallskip
\noindent
{\bf Proposition 3.2.} {\it  See Subsection 1.1 for $M$. Let $M_{\dbZ_p}$ be a $\dbZ_p$-structure of $M$ i.e., a free $\dbZ_p$-module such that $M=M_{\dbZ_p}\otimes_{\dbZ_p} W(k)$. Let $D$ be a smooth, closed subgroup scheme of $\pmb{GL}_M$ such that: (i) its generic fibre $D_{B(k)}$ is connected and (ii) there exists a direct summand $\grd$ of $\End_{\dbZ_p}(M_{\dbZ_p})$ with $\Lie(D)=\grd\otimes_{\dbZ_p} W(k)$. Then there exists a unique smooth, closed subgroup scheme $D_{\dbZ_p}$ of $\pmb{GL}_{M_{\dbZ_p}}$ such that $D=D_{\dbZ_p}\times_{\Spec(\dbZ_p)} \Spec(W(k))$.}
\medskip

{\bf Proof:}
Let $q_1, q_2:\Spec(B(k)\otimes_{\dbQ_p} B(k))\to \Spec(B(k))$ be the two projections. As $B(k)\otimes_{\dbQ_p} B(k)$ is a reduced $\dbQ_p$-algebra, the smooth closed subgroup schemes $q_1^*(D_{B(k)})$ and $q_2^*(D_{B(k)})$ of $\pmb{GL}_{M_{\dbZ_p}\otimes_{\dbZ_p} B(k)\otimes_{\dbQ_p} B(k)}$ are reduced and their Lie algebras are equal to  $\grd\otimes_{\dbZ_p} B(k)\otimes_{\dbQ_p} B(k)$. Thus the fibres of $q_1^*(D_{B(k)})$ and $q_2^*(D_{B(k)})$ are naturally identified, cf. [1, Ch. II, Subsect. 7.1]. Therefore we can identify $q_1^*(D_{B(k)})=q_2^*(D_{B(k)})$. Due to this identification, we have a canonical descent datum on the subgroup $D_{B(k)}$ of $\pmb{GL}_{M[{1\over p}]}$ with respect to the morphism $\Spec(B(k))\to\Spec(\dbQ_p)$. As $D_{B(k)}$ is affine, the descent datum is effective (cf. [3, Ch. 6, Thm. 6]); thus $D_{B(k)}$ is the pull back of a subgroup $D_{\dbQ_p}$ of $\pmb{GL}_{M_{\dbZ_p}[{1\over p}]}$. Obviously $D_{\dbQ_p}$ is connected. As $\Lie(D_{\dbQ_p})\otimes_{\dbQ_p} B(k)=\Lie(D_{B(k)})=\grd[{1\over p}]\otimes_{\dbQ_p} B(k)$, we get $\Lie(D_{\dbQ_p})=\grd[{1\over p}]$. Let $D_{\dbZ_p}$ be the Zariski closure of $D_{\dbQ_p}$ in $\pmb{GL}_{M_{\dbZ_p}}$. As $D$ is the Zariski closure of $D_{B(k)}$ in $\pmb{GL}_M$, we have $D=D_{\dbZ_p}\times_{\Spec(\dbZ_p)} \Spec(W(k))$. Thus $D_{\dbZ_p}$ is smooth. The uniqueness of $D_{\dbZ_p}$ is implied by the fact (cf. [1, Ch. II, Subsect. 7.1]) that $D_{\dbQ_p}$ is the unique connected subgroup of $\pmb{GL}_{M_{\dbZ_p}[{1\over p}]}$ whose Lie algebra is $\grd[{1\over p}]$.\endproof

\smallskip
\noindent
{\bf Proposition 3.3.} {\it  Let $R\in \{k,\dbZ_p,W(k)\}$. Let $E$ be a reductive group scheme over $\Spec(R)$. Let $X_E$ and $\tilde X_E$ be two smooth, closed subgroup schemes of $E$ that have connected fibres and that contain a fixed torus $T_E$ of $E$. 

\medskip
{\bf (a)} The centralizer $C_E$ of $T_E$ in $E$ is a reductive, closed subgroup scheme of $E$. The normalizer $N_E$ of $T_E$ in $E$ is a smooth, closed subgroup scheme of $E$  whose identity component is $C_E$.

\smallskip
{\bf (b)} If $T_E$ is a maximal torus of $E$, then the normalizer of $X_E$ in $E$ is represented by a smooth, closed, subgroup scheme of $E$ whose identity component is $X_E$. 

\smallskip
{\bf (c)} Suppose $T_E$ is a maximal torus of $E$. If $\Lie(X_E)=\Lie(\tilde X_E)$, then $X_E=\tilde X_E$. 

\smallskip
{\bf (d)} If $R=W(k)$, $X_{Ek}=\tilde X_{Ek}$, and $T_E$ is a maximal torus of $E$, then $X_E=\tilde X_E$.

\smallskip
{\bf (e)} The reductive group scheme $E$ has maximal tori.}
\medskip

{\bf Proof:} 
Part (a) follows from [6, Vol. III, Exp. XIX, Subsect. 2.8 and Prop. 6.3]. Part (b) is a particular case of [6, Vol. III, Exp. XXII, Cors. 5.3.10 and 5.3.18 (ii)]. Parts (c) and (d) are particular cases of [6, Vol. III, Exp. XXII, Cors. 5.3.5 and 5.3.7]. If $R=k$ or $R=W(k)$, then (e) is a particular case of [6, Vol. III, Exp. XIX, Prop. 6.1]. We check that (e) holds if $R=\dbZ_p$. Let $*$ be a maximal torus of $E_{\dbF_p}$, cf. [1, Ch. V, Thm. 18.2 (i)]. As the functor on the category of $\Spec(\dbZ_p)$-schemes that parameterizes maximal tori of $E$ is representable by a smooth $\Spec(\dbZ_p)$-scheme (cf. [6, Vol. II, Exp. XII, Cor. 1.10]), we get that $E$ has a maximal torus that lifts $*$.\endproof

\smallskip
\noindent
{\bf Lemma 3.4.} {\it Let $E$ be a reductive subgroup scheme of $\pmb{GL}_M$. Then $E$ is a closed subgroup scheme of $\pmb{GL}_M$; thus $\Lie(E)$ as a $W(k)$-module is a direct summand of $\End_{W(k)}(M)$. Moreover $E$ is the unique reductive subgroup scheme of $\pmb{GL}_M$ whose Lie algebra is $\Lie(E)$.}
\medskip

{\bf Proof:}
The first part is a particular case of [6, Vol. II, Exp. XVI, Cor. 1.5 a)]. Let $E_1$ be a reductive subgroup scheme of $\pmb{GL}_M$ such that $\Lie(E_1)=\Lie(E)$. We have $E_{1B(k)}=E_{B(k)}$, cf. [1, Ch. II, Subsect. 7.1]. Thus $E_1$ and $E$ are both equal to the Zariski closure in $\pmb{GL}_M$ of $E_{1B(k)}=E_{B(k)}$ . \endproof

\medskip
\noindent
{\bf 3.1. Decompositions} 

Let $*\in\{\dbF_p,k\}$. Let $E$ be an adjoint group scheme over $\Spec(W(*))$. Let $E_*=\prod_{j\in J(E)} E^j_*$ be the product decomposition into simple, adjoint groups over $*$ (cf. [33, Subsubsect. 3.1.2]). This decomposition lifts uniquely to a product decomposition $E=\prod_{j\in J(E)} E^j$ into simple, adjoint group schemes over $\Spec(W(*))$ (cf. [6, Vol. III, Exp. XXIV, Prop. 1.21]). Suppose now $*=\dbF_p$ and $E_{\dbF_p}$ is a simple, adjoint group over $\dbF_p$. There exists a unique number $q\in\dbN$ such that $E_{\dbF_p}$ is isomorphic to $\Res_{\dbF_{p^q}/\dbF_p} \tilde E_{\dbF_{p^q}}$, where $\tilde E_{\dbF_{p^q}}$ is an absolutely simple, adjoint group over $\dbF_{p^q}$ (cf. [33, Subsubsect. 3.1.2]). From [6, Vol. III, Exp. XXIV, Prop. 1.21] we get that $\tilde E_{\dbF_{p^q}}$ lifts uniquely to an adjoint group scheme $\tilde E$ over $\Spec(W(\dbF_{p^q}))$ and that $E$ is isomorphic to $\Res_{W(\dbF_{p^q})/\dbZ_p} \tilde E$.

\medskip
\noindent
{\bf 3.2. Trace maps} 

Let $\scrM$ be a free $W(k)$-module of finite rank. Let $\scrE$ be a reductive subgroup scheme of $\pmb{GL}_{\scrM}$. The Lie subalgebra $\Lie(\scrE)$ of $\End_{W(k)}(\scrM)$, when viewed as a $W(k)$-submodule, is a direct summand (cf. Lemma 3.4). Let $\rho:\scrE\hookrightarrow \pmb{GL}_{\scrM}$ be the resulting faithful representation. We refer to [9, Part III, Sect. 20] for spin representations over $\dbC$ and to [13, Ch. 27] for their versions over $\Spec(\dbZ)$. For $x$, $y\in\Lie(\scrE)$ let $\Tr(x,y)$ be the trace of the endomorphism $xy$ of $\scrM$. We get a symmetric bilinear map $\Tr:\Lie(\scrE)\times\Lie(\scrE)\to W(k)$.

\smallskip
\noindent
{\bf Lemma 3.5.} {\it The trace map $\Tr$ is perfect i.e., it induces naturally a $W(k)$-linear isomorphism from $\Lie(\scrE)$ onto $\Hom_{W(k)}(\Lie(\scrE),W(k))$, if any one of the following disjoint five conditions holds:

\medskip
{\bf (i)} $\rho$ is an isomorphism;

\smallskip
{\bf (ii)} $\rho$ is the standard faithful representation $\pmb{SL}_n\hookrightarrow \pmb{GL}_n$ and $p$ does not divide $n$;

\smallskip
{\bf (iii)} $p>2$ and $\rho$ is the spin faithful representation $\pmb{GSpin}_{2n+1}\hookrightarrow \pmb{GL}_{2^n}$ ($n\ge 1$);

\smallskip
{\bf (iv)} $p>2$ and there exists $n\in 2\dbN$ such that $\rho$ is either $\pmb{Sp}_{n}\hookrightarrow \pmb{GL}_{n}$ or $\pmb{SO}_{n}\hookrightarrow \pmb{GL}_{n}$;

\smallskip
{\bf (v)} $p>2$ and $\rho$ is a spin faithful representation $\pmb{GSpin}_{2n}\hookrightarrow \pmb{GL}_{2^{n}}$.}
\medskip

{\bf Proof:}
The Lemma is well known if (i) holds. Thus we can assume that one of the conditions (ii) to (v) holds. Let $\Tr_k:\Lie(\scrE_k)\times \Lie(\scrE_k)\to k$ be the reduction mod $p$ of $\Tr$. Let 
$$\scrJ:=\{x\in\Lie(\scrE_k)|\Tr_k(x,y)=0\;\;\forall y\in\Lie(\scrE_k)\}.$$ As for $x$, $y$, $z\in\Lie(\scrE)$ we have $\Tr(x,[y,z])=\Tr([x,y],z)$, $\scrJ$ is an ideal of $\Lie(\scrE_k)$. If (ii) or (iv) holds, then $\Lie(\scrE_k)$ has no proper ideals (cf. [14, Subsect. (0.13)]) and it is easy to see that $\scrJ\neq \Lie(\scrE_k)$; thus $\scrJ=0$.  

Suppose (iii) or (v) holds. For simplicity, we refer to the case (iii) as (v) is the same. Let $\tilde\mu:\dbG_m\to\scrE$ be a cocharacter that defines a cocharacter of $\pmb{GL}_{\scrM}$ of weights $\{-1,0\}$, cf. [30, p. 186] applied over $B(k)$ and the fact that, as $\scrE$ is split, each cocharacter of $\scrE_{B(k)}$ is $\scrE(B(k))$-conjugate to the generic fibre of a cocharacter of $\scrE$. As $p>2$, we have a direct sum decomposition $\Lie(\scrE)=\Lie(\scrE^{\der})\oplus\Lie(Z(\scrE))$ and $\Lie(\scrE^{\der})$ is the ideal of $\Lie(\scrE)$ of trace $0$ elements. Let $\tilde l_0\in\Lie(\scrE)$ be the image under $d\tilde\mu$ of the standard generator of $\Lie(\dbG_m)$. The multiplicities of the eigenvalues $-1$ and $0$ of $\tilde l_0$ are $2^{n-1}$, cf. [28, Table 4.2]. It is easy to see that $\Tr(\tilde l_0+{1\over 2}1_{\scrM},\tilde l_0+{1\over 2}1_{\scrM})=2^{n-2}$. Thus $\tilde l_0+{1\over 2}1_{\scrM}$ is an element of $\Lie(\scrE^{\der})$ (as its trace is $0$) whose reduction mod $p$ does not belong to $\scrJ$ (as $p>2$). Moreover $\Lie(Z(\scrE_k))=k1_{\scrM/p\scrM}$ is not included in $\scrJ$. But $\Lie(\scrE^{\der}_k)$ and $\Lie(Z(\scrE_k))$ are the only proper ideals of $\Lie(\scrE_k)=\Lie(\scrE^{\der}_k)\oplus\Lie(Z(\scrE_k))$, cf. [14, Subsect. (0.13)]. Thus $\scrJ=0$. 

As in all cases we have $\scrJ=0$, the trace map $\Tr$ is perfect.\endproof

\medskip
\noindent
{\bf 3.3. On $G$ and its closed subgroup schemes} 

Let $M=F^1\oplus F^0$, $\mu:\dbG_m\to G$, and $l_0$ be as in Subsections 1.1 and 1.1.3. We will recall basic properties of closed subgroup schemes of $G$ and we introduce the basic notations we will often need in what follows.  

\smallskip
\noindent
{\bf 3.3.1. Root decomposition} 

The cocharacter $\mu$ is either trivial or a closed embedding monomorphism. Thus $\im(\mu)$ is a torus of $G$. The centralizer of $\im(\mu)$ (equivalently of $l_0$) in $G$ is a reductive, closed subgroup scheme $L$ of $G$, cf. Proposition 3.3 (a). Let $T$ be a maximal torus of $L$ (cf. Proposition 3.3 (e)); it is also a maximal torus of $G$. Let
$$\Lie(G)=\Lie(T)\oplus\bigoplus_{\alpha\in\Phi} {\grg}_{\alpha}$$ 
be the root decomposition relative to $T$: $\Phi$ is a root system and $T$ acts on the rank 1 direct summand ${\grg}_{\alpha}$ of $\Lie(G)$ via its character $\alpha$. For $\alpha\in\Phi$, there exists a unique $\dbG_a$-subgroup scheme $\dbG_{a,\alpha}$ of $G$ which is closed, normalized by $T$, and with $\Lie(\dbG_{a,\alpha})=\grg_{\alpha}$ (cf. [6, Vol. III, Exp. XXII, Thm. 1.1]). Let $\bar\grg_{\alpha}:=\grg_{\alpha}/p\grg_{\alpha}\subseteq\Lie(G_k)$. 

Let $\Phi_c$ be a closed subset of $\Phi$; we have $(\Phi_c+\Phi_c)\cap\Phi\subseteq\Phi_c\subseteq\Phi$. We recall that there exists a unique connected, smooth subgroup $H_{\Phi_c}$ of $G_k$ that contains $T_k$ and whose Lie algebra is $\Lie(T_k)\oplus\bigoplus_{\alpha\in\Phi_c} \bar\grg_{\alpha}$, cf. [6, Vol. III, Exp. XXII, Thm. 5.4.7]. The group $H_{\Phi_c}$ is generated by $T_k$ and by ${\dbG_{a,\alpha}}_k$'s with $\alpha\in\Phi_c$, cf. [6, Vol. III, Exp. XXII, Lemma 5.4.4]. Let $\Phi_{d}$ be another closed subset of $\Phi$. We have  $H_{\Phi_c}\leqslant H_{\Phi_d}$ if and only if $\Phi_c\subseteq\Phi_d$, cf. either [6, Vol. III, Exp. XXII, Prop. 5.4.5 (ii)] or the mentioned generation property.   

\smallskip
\noindent
{\bf 3.3.2. On tori} 

For $m\in\dbN$ let $W_m(k):=W(k)/(p^m)$. Let $T_1$ and $T_2$ be two tori of $G$ such that $T_{1k}\leqslant T_{2k}$. There exists a unique subtorus $T_3$ of $T_2$ such that $T_{3k}=T_{1k}$. Let $h_m\in {\Ker}(G(W_m(k))\to G(k))$ be such that $h_mT_{1W_m(k)}h_m^{-1}=T_{3W_m(k)}$, cf. [6, Vol. II, Exp. IX, Thm. 3.6]. We can assume $h_{m+1}$ lifts $h_m$, cf. loc. cit. Let $h\in\Ker(G(W(k))\to G(k))$ be the unique element that lifts $h_m$ for all $m\in\dbN$. We have $hT_1h^{-1}=T_3\leqslant T_2$.

As any two maximal tori of $G_k$ are $G(k)$-conjugate, from the previous paragraph we also get that any two maximal tori of $G$ are $G(W(k))$-conjugate.

\smallskip
\noindent
{\bf 3.3.3. Defining $N$ and $U$} 

We identify $\Hom_{W(k)}(F^1,F^0)$ with the the maximal direct summand of $\End_{W(k)}(M)$ on which $\mu$ acts via the identical character of $\dbG_m$. Let $N^{\text{big}}$ be the smooth, closed subgroup scheme of $\pmb{GL}_M$ defined by: if $R$ is a commutative $W(k)$-algebra, then $N^{\text{big}}(R):=1_{M\otimes_{W(k)} R}+\Hom_{W(k)}(F^1,F^0)\otimes_{W(k)} R$. We have an identity $\Lie(N^{\text{big}})=\Hom_{W(k)}(F^1,F^0)$. Let $N:=N^{\text{big}}\cap G$. As $\mu$ factors through $G$, the intersection $\Hom_k(F^1/pF^1,F^0/pF^0)\cap\Lie(G_k)$ equals $(\Hom_{W(k)}(F^1,F^0)\cap\Lie(G))\otimes_{W(k)} k$. From this we get that:
$$N(R)=1_{M\otimes_{W(k)} R}+(\Hom_{W(k)}(F^1,F^0)\cap\Lie(G))\otimes_{W(k)} R.$$ 
The group scheme $N$ is connected, smooth, commutative, and unipotent. Moreover we have an identity $\break\Lie(N)=\Hom_{W(k)}(F^1,F^0)\cap\Lie(G)$. If $x\in\Lie(N)$, then $x^2:=x\circ x=0$. 

Similarly, the rule $U(R)=1_{M\otimes_{W(k)} R}+(\Hom_{W(k)}(F^0,F^1)\cap\Lie(G))\otimes_{W(k)} R$ defines a connected, smooth, commutative, unipotent, closed subgroup scheme $U$ of $G$. We have $\Lie(U)=\Hom_{W(k)}(F^0,F^1)\cap\Lie(G)$. If $y\in\Lie(U)$, then $y^2:=y\circ y=0$.

\smallskip
\noindent
{\bf 3.3.4. Defining $Q$ and $P$}

As the cocharacter $\mu$ factors through $G$ and produces a direct sum decomposition $M=F^1\oplus F^0$, we have a direct sum decomposition of $T$-modules
$$\Lie(G)=\Lie(N)\oplus\Lie(L)\oplus\Lie(U).$$
This implies that we have a disjoint union $\Phi=\Phi_U\cup\Phi_{L}\cup\Phi_N$, where 
$$\Phi_N:=\{\alpha\in\Phi|\grg_{\alpha}\subseteq\Lie(N)\},\;\;\Phi_{L}:=\{\alpha\in\Phi|\grg_{\alpha}\subseteq\Lie(L)\},\;\;\text{and}\;\;\Phi_U:=\{\alpha\in\Phi|\grg_{\alpha}\subseteq\Lie(U)\}.$$
Moreover $\Phi_U=-\Phi_N$. It is easy to see that we have product decompositions
$$N=\prod_{\alpha\in\Phi_N} \dbG_{a,\alpha}\;\;\;\text{and}\;\;\;U=\prod_{\alpha\in\Phi_U} \dbG_{a,\alpha}.$$
\indent
The direct sum $\Lie(N)\oplus\Lie(L)$ (resp. $\Lie(U)\oplus\Lie(L)$) is a Lie subalgebra of $\Lie(G)$. If $\alpha$, $\beta$, $\alpha+\beta\in\Phi$, then $[\grg_{\alpha}[{1\over p}],\grg_{\beta}[{1\over p}]]=\grg_{\alpha+\beta}[{1\over p}]$ (cf. Chevalley rule). The last two sentences imply that $\Phi_Q:=\Phi_N\cup\Phi_{L}$ (resp. $\Phi_P:=\Phi_U\cup\Phi_{L}$) is a closed subset of $\Phi$. Moreover $\Phi_Q\cup -\Phi_Q=\Phi_P\cup -\Phi_P=\Phi$. Thus $\Phi_Q$ (resp. $\Phi_P$) is a parabolic subset of $\Phi$. Thus $\Lie(Q)$ (resp. $\Lie(P)$) is the Lie algebra of a unique parabolic subgroup scheme $Q$ (resp. $P$) of $G$ that contains $T$, cf. [6, Vol. III, Exp. XXVI, Prop. 1.4]. As $\Phi_N=\Phi_Q\setminus -\Phi_Q$ (resp. $\Phi_U=\Phi_P\setminus -\Phi_P$), from [6, Vol. III, XXVI, Prop. 1.12] we get that the unipotent radical of $Q$ (resp. $P$) is $N$ (resp. $U$). Moreover $L$ is the unique Levi subgroup scheme of $Q$ (resp. $P$) that contains $T$, cf. loc. cit. 

Let $\tilde P$ be the normalizer of $F^1$ in $G$. Its Lie algebra $\Lie(\tilde P_{B(k)})=\{x\in\Lie(G_{B(k)})|x(F^1[{1\over p}])\subseteq F^1[{1\over p}]\}$ equals $\Lie(P_{B(k)})=\Lie(L_{B(k)})\oplus\Lie(U_{B(k)})$. Thus $P_{B(k)}$ is the identity component of $\tilde P_{B(k)}$, cf. [1, Ch. II, Subsect. 7.1] and the fact that $P_{B(k)}$ is connected. Thus $\tilde P_{B(k)}$ is a parabolic subgroup of $G_{B(k)}$ and therefore it is also connected. Thus $\tilde P_{B(k)}=P_{B(k)}$ and $P\leqslant \tilde P$. A similar argument shows that $P_k=\tilde P_k$. Thus $\tilde P=P$ i.e., $P$ is the normalizer of $F^1$ in $G$. A similar argument shows that $Q$ is the normalizer of $F^0$ in $G$. Thus $P\cap Q=L$. This implies that $P$ and $Q$ are opposite parabolic subgroups of $G$, cf. [6, Vol. III, Exp. XXVI, Thm. 4.3.2]. 

The product $W(k)$-morphism $P\times_{\Spec(W(k))} N\to G$ is an open embedding at identity sections (cf. [6, Vol. III, Exp. XXII, Prop. 4.1.2]) and $P\cap N$ is the identity section. Thus $P\times_{\Spec(W(k))} N\to G$ is an open embedding.

\bigskip
\noindent
{\boldsectionfont \S4. Basic properties of Shimura $F$-crystals over $k$}

\bigskip
Subsection 4.1 lists notations that will be often used throughout the paper. Subsection 4.2 introduces some $\dbZ_p$-structures. Subsection 4.3 recalls a basic result of [17] and includes simple properties of $D$-truncations that lead to the existence of the set $R_{\pmb{GL}_M}$. Subsection 4.4 defines the simple factors of the pair $(\Lie(G^{\ad}),\phi)$. 

\medskip
\noindent
{\bf 4.1. Standard notations} 

Always the notations 
$$M=F^1\oplus F^0,\;\mu:\dbG_m\to G,\;\scrC_g=(M,g\phi,G),\;\dbL\leqslant G(W(k)),\;P,\;T\leqslant N_T\leqslant G,\;\vartheta,\;\bar M,$$
$$\;\bar{\phi},\bar\vartheta,\;\scrC_{\bar g}=(\bar M,\bar g\bar\phi,\bar \vartheta\bar g^{-1},G_k),\;l_0\in\End_{W(k)}(M),\;w\in W_G,\;g_w\in G(W(k)),\;W_P\leqslant W_G,$$
$$\dbI_{\infty}=\{<\scrC_g>|g\in G(W(k))\},\,\,\text{and}\,\,\dbI=\{<\scrC_{\bar g}>|\bar g\in G(k)\}$$
are as in Subsection 1.1. Let $r:=\text{rk}(M)$ and $d:=\text{rk}(F^1)$. Let $\{e_1,\ldots,e_r\}$ be a $W(k)$-basis for $M$ such that $\{e_1,\ldots,e_d\}\subseteq F^1$, $\{e_{d+1},\ldots,e_r\}\subseteq F^0$, and $T$ normalizes $W(k)e_i$ for all $i\in\{1,\ldots,r\}$. Let $\{e_{i_1,i_2}|1\le i_1,i_2\le r\}$ be the standard $W(k)$-basis for $\End_{W(k)}(M)$ such that $e_{i_1,i_2}(e_{i_3})\in\{0,e_{i_1}\}$ is $e_{i_1}$ if and only if $i_2=i_3$. The element $l_0$ is the image under $d\mu$ of the standard generator of $\Lie(\dbG_m)$; thus $l_0\in\Lie(T)$. For $x\in M\cup\End_{W(k)}(M)$, let $\bar x$ be the reduction mod $p$ of $x$. Thus $\bar e_i\in\bar M$, $\bar l_0\in\Lie(T_k)$, etc. If $x\in\End_{W(k)}(M)$ and $g\in G(W(k))$, let $g(x):=gxg^{-1}$. Let $L$, $\Lie(G)=\Lie(T)\oplus\bigoplus_{\alpha\in\Phi} \grg_{\alpha}$, and $\dbG_{a,\alpha}$, $N$, $U$, $Q$, $P$, $\Phi_N$, $\Phi_{L}$, $\Phi_{U}$, $\Phi_Q$, and $\Phi_P$ be as in Subsection 3.3. Let $B$ be a Borel subgroup scheme of $G$ such that $T\leqslant B\leqslant P$. Let $B^{\text{opp}}$ be the Borel subgroup scheme of $G$ that is the opposite of $B$ with respect to $T$; we have $B\cap B^{\text{opp}}=T$.  

\smallskip
\noindent
{\bf 4.1.1. Notations of combinatorial nature}

For $\alpha\in\Phi$ we define a number 
$$\eps(\alpha)\in\{-1,0,1\}$$ 
by the rule: $\eps(\alpha)=-1$ if $\alpha\in\Phi_N$, $\eps(\alpha)=0$ if $\alpha\in\Phi_{L}$, and $\eps(\alpha)=1$ if $\alpha\in\Phi_U$. As $g_w\phi(\Lie(T))=\Lie(T)$ (cf. (1)), there exists an automorphism of root systems 
$$\pi_w:\Phi\arrowsim\Phi$$ 
such that $g_w\phi(\grg_{\alpha}[{1\over p}])=\grg_{\pi_w(\alpha)}[{1\over p}]$. Let $O(w)$ be the set of orbits of $\pi_w$ (viewed as a permutation of $\Phi$). For $o\in O(w)$ let $|o|$ be the number of elements of $o$, let $\grg_o:=\oplus_{\alpha\in o} \grg_{\alpha}$, and let $\bar\grg_o:=\grg_o/p\grg_o$. Let $O(w)^0:=\{o\in O(w)|o\subseteq\Phi_L\}$.  
For $\alpha\in\Phi_N$ let 
$$w_{\alpha}\in\dbN$$ 
be the smallest number such that $\pi_w^{w_{\alpha}}(\alpha)\notin\Phi_{L}$ (i.e., such that $\pi_w^{w_{\alpha}}(\alpha)\in\Phi_N\cup\Phi_U$). Let 
$$\Phi_N^{+w}:=\{\alpha\in\Phi_N|\pi_w^{w_{\alpha}}(\alpha)\in\Phi_U\}=\{\alpha\in\Phi_N|-\pi_w^{w_{\alpha}}(\alpha)\in\Phi_N\}\;\;\;\text{and}$$
$$\Phi_N^{-w}:=\{\alpha\in\Phi_N|\pi_w^{w_{\alpha}}(\alpha)\in\Phi_N\}=\{\alpha\in\Phi_N|-\pi_w^{w_{\alpha}}(\alpha)\in\Phi_U\}=\Phi_N\setminus\Phi_N^{+w}.$$ 
The function $\dbS:W_G\to\dbN\cup\{0\}$ mentioned in Subsubsection 1.1.3 is defined by the rule 
$$\dbS(w):=\text{the}\;\;\text{number}\;\;\text{of}\;\;\text{elements}\;\;\text{of}\;\;\Phi_N^{+w}.$$ 
Let $o^+:=o\cap\Phi_N^{+w}$. Let $O(w)^+:=\{o\in O(w)|o^+\neq\emptyset\}$. For $o\in O(w)^+$ let 
$$\tilde o:=\cup_{\alpha\in o^+} \{\pi_w^i(\alpha)|1\le i\le w_{\alpha}\}\subseteq o\cap\Phi_P.$$ 

\medskip
\noindent
{\bf 4.2. Some $\dbZ_p$-structures} 

We start by not assuming that (1) holds. Let 
$$\sigma_g:=g\phi\mu(p),\;\; \text{where}\;\;g\in G(W(k)).$$ 
It is a $\sigma$-linear automorphism of $M$ that normalizes $\Lie(G)$. As $k$ is algebraically closed, 
$$M_{\dbZ_p}^g:=\{x\in M|\sigma_g(x)=x\}$$ 
is a $\dbZ_p$-structure of $M$ and $\Lie(G)$ is the tensorization with $W(k)$ of a Lie subalgebra of $\End_{\dbZ_p}(M^g_{\dbZ_p})$. Thus $G$ is the pull back to $\Spec(W(k))$ of a reductive, closed subgroup scheme $G_{\dbZ_p}^g$ of $\pmb{GL}_{M_{\dbZ_p}^g}$, cf. Proposition 3.2. For $i\in\dbZ$, $\sigma_g^i$ acts naturally on: (i) subgroup schemes and cocharacters of either $G_k$ or $G$, (ii) the groups $G(k)$ and $G(W(k))$, and (iii) $\Lie(G)$. For instance, $\sigma_g^i(\mu)$ is the $\sigma_g^i$-conjugate of $\mu$ i.e., it is the cocharacter of $G$ whose image has a Lie algebra generated by $\sigma_g^i(l_0)=\sigma_g^i l_0\sigma_g^{-i}\in\End_{W(k)}(M)$ and which acts on $M$ via the trivial and the inverse of the identical character of $\dbG_m$. We have $g\phi(\Lie(T))=\sigma_g(\Lie(T))=\Lie(\sigma_g(T))$. 

Let $\bar g_{\text{corr}}\in G(k)$ be such that $\bar g_{\text{corr}}\sigma_{1_M}(T_k)\bar g_{\text{corr}}^{-1}=T_k$. Let $g_{\text{corr}}\in G(W(k))$ be such that it lifts $\bar g_{\text{corr}}$ and we have $g_{\text{corr}}\sigma_{1_M}(T)g_{\text{corr}}^{-1}=\sigma_{g_{\text{corr}}}(T)=T$, cf. Subsubsection 3.3.2. Thus $g_{\text{corr}}\phi(\Lie(T))=\sigma_{g_{\text{corr}}}(\Lie(T))=\Lie(T)$. The two sets $\{\scrC_g|g\in G(W(k))\}$ and $\{\scrC_{gg_{\text{corr}}}|g\in G(W(k))\}$ are equal. Therefore for the study of $\scrC_g$'s and of their $D$-truncations mod $p$, we can replace $\phi$ by $g_{\text{corr}}\phi$. Thus to ease notations, until the end we will assume that $g_{\text{corr}}=1_M$. Therefore (1) holds i.e., we have $\phi(\Lie(T))=\Lie(T)$. 

\smallskip
Let $w\in W_G$. We have $\Lie(\sigma_{g_w}(T))=\sigma_{g_w}(\Lie(T))=g_w\phi(\Lie(T))=\Lie(T)$. Thus $\{x\in\Lie(T)|\sigma_{g_w}(x)=x\}$ is a direct summand of $\End_{\dbZ_p}(M_{\dbZ_p})$ and a $\dbZ_p$-structure of $\Lie(T)$. Thus $T$ is the pull back to $\Spec(W(k))$ of a maximal torus $T_{\dbZ_p}^{g_w}$ of $G^{g_w}_{\dbZ_p}$ (cf. Proposition 3.2); in general the isomorphism class of $T^{g_w}_{\dbZ_p}$ depends on $w$. For $t\in T(W(k))$ let $t_1\in T(W(k))$ be such that $t=t_1\sigma_{g_w}(t_1^{-1})$, cf. Lemma 3.1 applied to $T^{g_w}_{\dbZ_p}$. We have $tg_w\phi=t_1g_w\phi t_1^{-1}$. This implies that:

\smallskip
\noindent
{\bf Lemma 4.3.} {\it 
Let $w\in W_G$. The inner isomorphism class $<\scrC_{g_w}>$ depends only on $w$ and not on the choice of the representative $g_w\in N_T(W(k))$ of $w$. Also, if $g_1\in\Ker(G(W(k))\to G(k))$, then $<\scrC_{g_1tg_w}>=<\scrC_{g_2g_w}>$, where $g_2:=t_1^{-1}g_1t_1\in\Ker(G(W(k))\to G(k))$.}

\smallskip
\noindent
{\bf 4.2.1. Extra notations} 

Let 
$$G_{\dbZ_p}:=G_{\dbZ_p}^{1_M}\;\;\; {\text{and}}\;\;\;T_{\dbZ_p}:=T_{\dbZ_p}^{1_M}.$$ 
In what follows we will often denote $\sigma_{1_M}$ simply by $\sigma$. Thus we have 
$$\sigma_g=g\sigma\;\;\;\text{and}\;\;\;\phi=\sigma\mu({1\over p}).$$ 
If $h\in \dbL$, then $\mu({1\over p})h\mu(p)\in G(W(k))$ as $\mu({1\over p})h\mu(p)(M)=\mu({1\over p})h(M+{1\over p}F^1)=\mu({1\over p})(M+{1\over p}F^1)=M$. Thus 
$$\phi(h):=\phi h\phi^{-1}=\sigma\mu({1\over p})h\mu(p)\sigma^{-1}=\sigma(\mu({1\over p})h\mu(p))\in G(W(k)).$$ 
We have an abstract group action $\dbT_{G,\phi}:\dbL\times G(W(k))\to G(W(k))$ of $\dbL$ on the set $G(W(k))$, defined by
$$\dbT_{G,\phi}(h,g)=hg\phi(h^{-1}).$$ 
The set of orbits of this action is in one-to-one correspondence to $\dbI_{\infty}$. In this paper we will mainly study the mod $p$ version of this action. If $\alpha\in\Phi_N$ (resp. if $\alpha\in\Phi_U$) and if $x\in\grg_{\alpha}$, then for $h=1_M+px$ (resp. for $h=1_M+x$), we have $\phi(h)=1_M+\sigma(x)$ (resp. we have $\phi(h)=1_M+p\sigma(x)$). In particular, the reduction mod $p$ of each element of $\phi(U(W(k)))$ is $1_{\bar M}$ (i.e., we have $\phi(U(W(k)))\leqslant \Ker(G(W(k))\to G(k))$).

We consider the product decomposition into simple, adjoint groups over $\Spec(\dbZ_p)$ (cf. Subsection 3.1) 
$$G^{\ad}_{\dbZ_p}=\prod_{j\in J} G^j_{\dbZ_p}.\leqno (6)$$ 
If $G$ is a torus, then as a convention $J=\emptyset$. Let $J^c:=\{j\in J|\mu\;\text{has}\;\text{a}\;\text{trivial}\;\text{image}\;\text{in}\;G^j_{W(k)}\}$. 

\smallskip
\noindent
{\bf Lemma 4.2.} {\it  Suppose $\mu$ factors through $Z^0(G)$ (equivalently through $Z(G)$). Then $\dbI_{\infty}=\{<\scrC_{1_M}>\}$ i.e., for all $g\in G(W(k))$ we have $<\scrC_g>=<\scrC_{1_M}>$.}
\medskip

{\bf Proof:}
We have $G=L=P$. Thus $\dbL=G(W(k))$ and for $h\in G(W(k))$ we have $\phi(h)=\sigma(h)$. Let $h\in G(W(k))$ be such that $g=h\sigma(h^{-1})$, cf. Lemma 3.1 applied to $G_{\dbZ_p}$. Therefore $g\phi h$ is equal to $h\sigma(h^{-1})\phi h=h\phi(h^{-1})\phi h=h\phi$. Thus $h$ is an inner isomorphism between $\scrC_{1_M}$ and $\scrC_g$.\endproof

\smallskip
\noindent
{\bf Lemma 4.3.} {\it  If $j\in J\setminus J^c$, then each simple factor $G^j_{0}$ of $G^j_{W(k)}$ is of classical Lie type.}
\medskip

{\bf Proof:}
Let $i\in\dbZ$ be such that $\sigma^i(\mu_{B(k)})$ has a non-trivial image in $G^j_{0B(k)}$. Let $\tilde G^{j}_{0B(k)}$ be the reductive subgroup of $G_{B(k)}$ that is generated by $\im(\sigma^i(\mu_{B(k)}))$ and by the unique normal, semisimple subgroup of $G_{B(k)}$ whose adjoint is  $G^j_{0B(k)}$. The cocharacter $\sigma^i(\mu_{B(k)})$ factors through $\tilde G^{j}_{0B(k)}$; it is also a cocharacter of $\pmb{GL}_{M[{1\over p}]}$ of weights $\{-1,0\}$ whose $\tilde G^{j}_{0B(k)}(B(k))$-conjugates generate $\tilde G^{j}_{0B(k)}$. Therefore the simple, adjoint group $\tilde G^{j\ad}_{0B(k)}=G^j_{0B(k)}$ is of classical Lie type, cf. [30, Cor. 1, p. 182].\endproof

\smallskip
\noindent
{\bf 4.2.2. Epsilon strings and slopes} 

Let $w\in W_G$. Let $o\in O(w)$. Let $m_w^+(o)$ and $m_w^-(o)$ be the number of elements of the intersections $o\cap \Phi_U$ and $o\cap\Phi_N$ (respectively). We write $o=\{\alpha^1,\ldots,\alpha^{|o|}\}$, the numbering being such that $\pi_w(\alpha^i)=\alpha^{i+1}$ if $i\in\{1,\ldots,|o|\}$, where $\alpha^{|o|+1}:=\alpha^1$. By the epsilon string of $o$ with respect to $w$ we mean the $|o|$-tuple 
$$\eps_w(o):=(\eps(\alpha^1),\ldots,\eps(\alpha^{|o|}))\in\{-1,0,1\}^{|o|};$$ 
it is uniquely determined by $o$ and $w$ up to a circular rearrangement. 

As $\phi=\sigma\mu({1\over p})$, for $\alpha\in\Phi$ we have $g_w\phi(\grg_{\alpha})=\sigma_{g_w}(p^{\eps(\alpha)}\grg_{\alpha})=p^{\eps(\alpha)}\grg_{\pi_w(\alpha)}$. Thus for $i\in\{1,\ldots,|o|\}$ we have $(g_w\phi)^{|o|}(\grg_{\alpha^i})=p^{m_w^+(o)-m_w^-(o)}\grg_{\alpha^i}$. Therefore all (Newton polygon) slopes of $(\grg_o,g_w\phi)$ are equal to
$${{m_w^+(o)-m_w^-(o)}\over {|o|}}.\leqno (7)$$
If $o\in O(w)^0$, then $m_w^+(o)=m_w^-(o)=0$ and all slopes of $(\grg_o,g_w\phi)$ are $0$. As $\Phi_U=-\Phi_N$ and $\Phi_N$ have the same number of elements, we have
$$\sum_{o\in O(w)} [m_w^+(o)-m_w^-(o)]=0.\leqno (8)$$ 
 
\medskip
\noindent
{\bf 4.3. On Kraft's work  and $D$-truncations} 

A truncated Barsotti--Tate group of level 1 over $k$ is a finite, flat, commutative group scheme $\bar D_1$ over $k$ that is annihilated by $p$ and such that the complex $\bar D_1\operatornamewithlimits{\to}\limits^{\bar F_1} \bar D_1\times_k {}_{\sigma} k\operatornamewithlimits{\to}\limits^{\bar V_1} \bar D_1$ is exact, where $\bar F_1$ and $\bar V_1$ are the Frobenius and the Verschiebung (respectively) homomorphisms. Let $D_1$ be a Barsotti--Tate group over $k$ such that $D_1[p]=\bar D_1$, cf. [15, Thm. 4.4 e)]. Let $(M_1,\phi_1)$ be the Dieudonn\'e module of $D_1$. The association $D_1\mapsto (M_1,\phi_1)$ induces an antiequivalence from the category of Barsotti--Tate groups over $k$ and the category of Dieudonn\'e modules over $k$ (see [8, Ch. III, Prop. 6.1 iii)]). Let $\vartheta_1:M_1\to M_1$ be the Verschiebung map of $(M_1,\phi_1)$. Let $\bar{M}_1$, $\bar{F}_1$, and $\bar\vartheta_1$ be the reductions mod $p$ of $M_1$, $\phi_1$, and $\vartheta_1$ (respectively). The classical Dieudonn\'e theories (for instance, see [8, 5.1 of Ch. III and p. 160]) tell us that the triple $\dbD(\bar D_1):=(\bar{M}_1,\bar {\phi}_1,\bar\vartheta_1)$ depends only on $\bar D_1$ and we have $\Ker(\bar{\phi}_1)=\im(\bar\vartheta_1)$ (and thus also $\im(\bar{\phi}_1)=\Ker(\bar\vartheta_1)$); moreover, the association $\bar D_1\mapsto \dbD(\bar D_1)$ induces an antiequivalence of (the corresponding) categories. Thus for each $\bar g\in G(k)$, there exists a (unique up to isomorphism) truncated Barsotti--Tate group $\bar D_{\bar g}$ of level 1 over $k$ such that we have $\dbD(\bar D_{\bar g})=(\bar M,\bar g\bar\phi,\bar\vartheta\bar g^{-1})$. 

Suppose $\text{rk}(M_1)=r$ and $\dim_k(\Ker(\bar{\phi}_1))=d$ i.e., $\bar D_1$ is a truncated Barsotti--Tate group of level 1 that has height $r$ and dimension $d$. In [17] (see [26, Subsects. (2.3) and (2.4)] and [22, Subsect. 2.1]) it is shown that there exists a $k$-basis $\{\bar b_1,\ldots,\bar b_r\}$ for $\bar{M}_1$ and a permutation $\pi$ of $\{1,\ldots,r\}$ such that for $i\in\{1,\ldots,r\}$ we have:

\medskip\noindent
--  $\bar{\phi}_1(\bar b_i)=0$ if $i\le d$, $\bar{\phi}_1(\bar b_i)=\bar b_{\pi(i)}$ if $i>d$, $\bar\vartheta_1(\bar b_{\pi(i)})=0$ if $i>d$, and $\bar\vartheta_1(\bar b_{\pi(i)})=\bar b_i$ if $i\le d$.

\medskip
Let $g_w\in \pmb{GL}_M(W(k))$ be such that $g_w\phi$ takes $e_i$ to $p^{n_i}e_{\pi(i)}$ for all $i\in\{1,\ldots,r\}$, where $n_i\in\{0,1\}$ is $1$ if and only if $i\le d$. If $G=\pmb{GL}_M$, then based on Lemma 4.1 we can assume that $g_w$ is as in Subsection 1.1. By mapping $\bar b_i\to \bar e_i$ we get an isomorphism $(\bar{M}_1,\bar {\phi}_1,\bar\vartheta_1)\arrowsim (\bar M,\bar g_w\bar\phi,\bar \vartheta\bar g_w^{-1})$. Thus we got:

\smallskip
\noindent
{\bf Lemma 4.4.} {\it  Suppose $G=\pmb{GL}_M$. Each truncated Barsotti--Tate group of level 1 over $k$ of height $r$ and dimension $d$ is isomorphic to the truncated Barsotti--Tate group of level 1 defined by $(\bar M,\bar g_w\bar\phi,\bar\vartheta \bar g_w^{-1})$ for some $w\in W_G=W_{\pmb{GL}_M}$. Thus for each $\bar g\in G(k)$, there exists $w\in W_G$ such that $\scrC_{\bar g}$ and $\scrC_{\bar g_w}$ are (inner) isomorphic.}

\smallskip
\noindent
{\bf Lemma 4.5.} {\it   
Let $g_1$ and $g_2\in G(W(k))$. 

\medskip
{\bf (a)} We have an identity $\scrC_{\bar g_1}=\scrC_{\bar g_2}$ if and only if $\bar g_2^{-1}\bar g_1\in\sigma(N)(k)$ (equivalently  if and only if $\bar g_2\bar g_1^{-1}\in (g_1\sigma(N)g_1^{-1})(k)=\sigma_{g_1}(N)(k)$).

\smallskip
{\bf (b)} We have $<\scrC_{\bar g_1}>=<\scrC_{\bar g_2}>$ if and only if there exists $h\in\dbL$ such that $h_{12}:=hg_1\phi(h^{-1})g_2^{-1}$ is an element of $\Ker(G(W(k))\to G(k))$.

\smallskip
{\bf (c)} We have $<\scrC_{\bar g_1}>=<\scrC_{\bar g_2}>$ if and only if there exists $h_{12}\in\Ker(G(W(k))\to G(k))$ such that $<\scrC_{g_1}>=<\scrC_{h_{12}g_2}>$.

\smallskip
{\bf (d)} Statement 1.4 (a) holds if and only if we have $\dbI=\{<\scrC_{\bar g_w}>|w\in W_G\}$.

\smallskip
{\bf (e)} If statement 1.4 (a) holds, then statement 1.4 (b) holds if and only if the set $\dbI$ has $[W_G:W_P]$ elements.}
\medskip

{\bf Proof:}
Let $N^{\text{big}}$ be as in Subsubsection 3.3.3. As $N_k=N_k^{\text{big}}\cap G_k$ (see Subsubsection 3.3.3), it suffices to prove (a) under the extra assumption that $G=\pmb{GL}_M$. Thus $N=N^{\text{big}}$. As $\bar g_2\bar{\phi}=\bar g_1\bar{\phi}$, we have $\bar g_2^{-1}\bar g_1\bar\phi=\bar\phi$ and thus $\bar g_2^{-1}\bar g_1$ fixes $\bar\phi(\bar M)=\Ker(\bar\vartheta)$. Thus $\Ker(\bar\vartheta)=\Ker(\bar\vartheta\bar g_2^{-1}\bar g_1)$. Due to this and the equality $\bar \vartheta=\bar \vartheta\bar g_2^{-1}\bar g_1$, we get that $\bar g_2^{-1}\bar g_1$ also fixes $\bar M/\text{Ker}(\bar \vartheta)$. But we have $\text{Ker}(\bar \vartheta)=\bar\phi(\bar M)=\sigma(F^0/pF^0)$. Thus $\bar g_2^{-1}\bar g_1$ fixes $\sigma(F^0/pF^0)$ and $\bar M/\sigma(F^0/pF^0)$. This implies that $\bar g_2^{-1}\bar g_1\in\sigma(N^{\text{big}})(k)=\sigma(N)(k)$. This proves (a). 

We prove (b). The if part is trivial. We check the only if part. Let $\bar h\in P(k)$ be an inner isomorphism between $\scrC_{\bar g_1}$ and $\scrC_{\bar g_2}$. By replacing $g_1$ with $h_1g_1\phi(h_1^{-1})$, where $h_1\in\dbL$ lifts $\bar h$, we can assume $\bar h=1_{\bar M}$. Thus $\scrC_{\bar g_1}=\scrC_{\bar g_2}$. The element $g_3:=\sigma_{g_1}^{-1}g_2g_1^{-1}\sigma_{g_1}\in G(W(k))$ is such that $\bar g_3\in N(k)$, cf. (a). By multiplying $g_2$ with an element of $\Ker(G(W(k))\to G(k))$, we can also assume that $g_3\in N(W(k))$. Thus we can write $g_3=1_M+u$, where $u\in\Lie(N)$ (see Subsubsection 3.3.3). Let $h:=1_M-pu\in\Ker(N(W(k))\to N(k))$. We have $\phi=g_1^{-1}g_1\phi=g_1^{-1}\sigma_{g_1}\mu({1\over p})$ and $\mu({1\over p})h\mu(p)=1_M-u=g_3^{-1}$. We compute 
$$g_2\phi(h)=g_2g_1^{-1}\sigma_{g_1}\mu({1\over p})h\mu(p)\sigma_{g_1}^{-1}g_1=g_2g_1^{-1}\sigma_{g_1} g_3^{-1}\sigma_{g_1}^{-1} g_1=g_2g_1^{-1}g_1g_2^{-1}g_1=g_1.$$ 
Thus $h_{12}=hg_1\phi(h^{-1})g_2^{-1}=hg_1g_1^{-1}=h\in \Ker(G(W(k))\to G(k))$. Thus (b) holds. 

The if part of (c) is obvious. The only if part of (c) follows from (b), as $h$ is an inner isomorphism between $\scrC_{g_1}$ and $\scrC_{h_{12}g_2}$. Parts (d) and (e) follow from (c).\endproof

\smallskip
\noindent
{\bf Corollary 4.6.} {\it If $G=\pmb{GL}_M$, then the statement 1.4 (a) holds i.e., the set $R_{\pmb{GL}_M}$ exists.}
\medskip

{\bf Proof:}
This is only the translation of Lemma 4.4 and Lemma 4.5 (d), cf. also Lemma 4.1.\endproof

\medskip
\noindent
{\bf 4.4. Simple factors of $(\Lie(G^{\ad}),\phi)$}

As $\phi$ normalizes $\Lie(G_{B(k)})$, it also normalizes $\Lie(G^{\der}_{B(k)})=[\Lie(G_{B(k)}),\Lie(G_{B(k)})]$. We identify $\Lie(G^{\ad})$ with a $W(k)$-Lie subalgebra of $\Lie(G^{\ad}_{B(k)})=\Lie(G^{\der}_{B(k)})$. We have a direct sum decomposition (cf. (6))
$$(\Lie(G^{\ad}),\phi)=\oplus_{j\in J} (\Lie(G^j_{W(k)}),\phi).$$
As $\phi=\sigma\mu({1\over p})$, we have $\phi(\Lie(G^j_{W(k)}))=\Lie(G^j_{W(k)})$ if and only if $j\in J^c$. Accordingly, if $j\in J\setminus J^c$ (resp. if $j\in J^c$) we say $(\Lie(G^j_{W(k)}),\phi)$ is a non-trivial (resp. is a trivial) simple factor of $(\Lie(G^{\ad}),\phi)$. Let now $j\in J\setminus J^c$. The Lie algebras of the simple factors of $G^j_{B(k)}$ are permuted transitively by $\sigma$ and thus also by $\phi=\sigma\mu({1\over p})$. Therefore these factors have the same Lie type $\scrX_j$ which (cf. Lemma 4.3) is a classical Lie type.

Next we introduce types that extrapolate  [5, Table 2.3.8]. For $n\in\dbN$ we say $(\Lie(G^j_{W(k)}),\phi)$ is of $A_n$, $B_n$, or $C_n$ type if $\scrX_j$ is $A_n$, $B_n$, or $C_n$ (respectively). Suppose $\scrX_j=D_n$ with $n\ge 4$. We say $(\Lie(G^j_{W(k)}),\phi)$ is of $D_n^{\dbR}$ type if the non-trivial images of $\sigma^i(\mu)$'s with $i\in\dbZ$ in each simple factor $G^j_0$ of $G^j_{W(k)}$, are $G^j_0(W(k))$-conjugate and the centralizers of these images in $G^j_0$ have adjoint group schemes of $D_{n-1}$ Lie type (i.e., these images define $\varpi_1^{\vee}$ cocharacters of $\pmb{SO}_{2n}^{\ad}$; see [28, Table 4.2]). We say $(\Lie(G^j_{W(k)}),\phi)$ is of $D_n^{\dbH}$ type if it is not of $D_n^{\dbR}$ type.

\bigskip
\noindent
{\boldsectionfont \S5. Proof of Basic Theorem A}

\bigskip
Subsection 5.1 shows that there exists a group action $\dbT_{G_k,\sigma}$ as in Subsubsection 2.2.1 which is the mod $p$ version of the action $\dbT_{G,\phi}$ of Subsubsection 4.2.1. The orbits of $\dbT_{G_k,\sigma}$ will parameterize the set $\dbI$. Propositions 5.3 and 5.4 are key ingredients in the proof of the long form (see Subsection 5.2) of Basic Theorem A. From Subsubsection 5.1.1 until Remark 5.5 we study group schemes of inner automorphisms of $\scrC_{\bar g}$'s and stabilizer subgroup schemes of $\dbT_{G_k,\sigma}$. The Section ends with two examples and a remark. Let $\sigma:W_G\arrowsim W_G$ be the automorphism defined by the automorphism $\sigma$ of the pair $(G,T)$ that fixes $(G_{\dbZ_p},T_{\dbZ_p})$. We use the notations of Subsections 4.1, 4.1.1, 4.2, and 4.2.1. If $w\in W_G$, then $\sigma(w)\in W_G$ is the unique element such that the $\sigma$-linear automorphisms of $\Lie(T)$ defined by $g_{\sigma(w)}\phi$ and $\phi g_w$ coincide.

\medskip
\noindent
{\bf 5.1. A group action}

As the product $W(k)$-morphism $P\times_{\Spec(W(k))} N\to G$ is an open embedding (see Subsubsection 3.3.4), each element $h\in\dbL$ (i.e., each $h\in G(W(k))$ such that $\bar h\in P(k))$ can be written uniquely as a product $h=h_0h_4$, where $h_0\in P(W(k))$ and $h_4\in\Ker(N(W(k))\to N(k))$. We write $h_0=h_1h_2$, where $h_1\in U(W(k))$ and $h_2\in L(W(k))$. Thus $h=h_1h_2h_4$. We recall that for $\tilde h\in\dbL$ we defined $\phi(\tilde h)=\phi\tilde h\phi^{-1}\in G(W(k))$ (see Subsubsection 4.2.1). Let $g\in G(W(k))$. Let $g_1:=hg\phi(h^{-1})\in G(W(k))$. As $\phi=\sigma\mu({1\over p})$, we have
$$g_1=hg\phi(h_4^{-1})\phi(h_2^{-1})\phi(h_1^{-1})=hg\phi(h_4^{-1})\sigma(h_2^{-1})\phi(h_1^{-1}).$$
We have $\phi(h_1^{-1})\in\Ker(G(W(k))\to G(k))$ (see Subsubsection 4.2.1). As $L$ normalizes $N$, we can rewrite 
$$\phi(h_4^{-1})\sigma(h_2^{-1})=\sigma(h_2^{-1})\phi(\tilde h_3^{-1})=\sigma(h_2^{-1})\sigma(h_3^{-1}),$$ 
where $\tilde h_3:=h_2h_4h_2^{-1}\in\Ker(N(W(k))\to N(k))$ is $1_M+pn_3$  for some element $n_3\in\Lie(N)$ and where $h_3:=1_M+n_3\in N(W(k))$. We get
$$\bar g_1=\bar h_1\bar h_2\bar g\sigma(\bar h_2^{-1})\sigma(\bar h_3^{-1}).\leqno (9)$$
\indent
We introduce the group action that ``defines" (9). We endow $H_k:=U_k\times_k L_k\times_k N_k$ with the semidirect product group structure (see Subsubsection 2.2.1). Note that $H_k$ is the special fibre of the Bruhat--Tits affine, smooth group scheme over $\Spec(W(k))$ whose generic fibre is $G_{B(k)}$ and whose group of $W(k)$-valued points is $\dbL$ (equivalently of the dilatation in the sense of [3, Ch. 3, Subsect. 3.2] of $G$ centered on $P_k$). Let 
$$\dbT_{G_k,\sigma}:=\dbT_{G_k,P_k,T_k,\sigma}:H_k\times_k G_k\to G_k$$ 
be a $\dbT$-action as in Subsubsection 2.2.1. At the level of $k$-valued points we have $\dbT_{G_k,\sigma}((\bar h_1,\bar h_2,\bar h_3),\bar g)=\bar h_1\bar h_2\bar g\sigma(\bar h_2^{-1})\sigma(\bar h_3^{-1})$, cf. (5).
Let $o_{\bar g}$ 
be the orbit of $\bar g\in G_k(k)$ under $\dbT_{G_k,\sigma}$. We consider the set of orbits 
$$\dbO:=\{o_{\bar g}|\bar g\in G_k(k)\}$$
of $\dbT_{G_k,\sigma}$. We define an order relation $\le$ on $\dbO$ as follows: for $\bar g_1, \bar g_2\in G_k(k)$, we have $o_{\bar g_1}\le o_{\bar g_2}$ if and only if the Zariski closure of $o_{\bar g_1}$ in $G_k$ contains $o_{\bar g_2}$. This order defines a topology on $\dbO$ as follows. A subset $\dbU$ of $\dbO$ is open if and only if for all $\bar g_1, \bar g_2\in G_k(k)$ such that $o_{\bar g_1}\le o_{\bar g_2}$ and $o_{\bar g_2}\in\dbU$, we have $o_{\bar g_1}\in\dbU$.

\smallskip
\noindent
{\bf Lemma 5.1.} {\it  {\bf (a)} The association $<\scrC_{\bar g}>\mapsto o_{\bar g}$ defines a bijection $b_{\text{can}}:\dbI\arrowsim\dbO$. 

\smallskip
{\bf (b)} Conjecture 2.1 holds for the quadruple $(G_k,P_k,T_k,\sigma)$ if and only if the statements 1.4 (a) and (b) hold.}
\medskip

{\bf Proof:}
We prove (a). As Formula (9) says $\bar g_1=\dbT_{G_k,\sigma}((\bar h_1,\bar h_2,\bar h_3),\bar g)$, we have $o_{\bar g}=o_{\bar g_1}$. Thus the fact that the map $b_{\text{can}}$ is well defined follows from Lemma 4.5 (b). Let $g_2\in G(W(k))$ be such that $o_{\bar g_2}=o_{\bar g}$. Each element $(\bar h'_1,\bar h'_2,\bar h'_3)\in H_k(k)$ is associated to an element $h'\in\dbL$ in the same way we associated $(\bar h_1,\bar h_2,\bar h_3)$ to $h\in\dbL$. Thus if $(\bar h'_1,\bar h'_2,\bar h'_3)$ is such that $\bar g_2=\dbT_{G_k,\sigma}((\bar h'_1,\bar h'_2,\bar h'_3),\bar g)$, then we have $h_{12}' g_2=h' g\phi(h^{'-1})$ for some $h_{12}'\in\Ker(G(W(k))\to G(k))$. Thus $<\scrC_{\bar g_2}>=<\scrC_{\bar g}>$, cf. Lemma 4.5 (b). Thus the map $b_{\text{can}}$ is injective. As this map is obviously surjective, part (a) holds. Part (b) follows from (a) and Lemma 4.5 (d) and (e).\endproof

\smallskip
\noindent
{\bf 5.1.1. Some group schemes} 

Let $\Phi_R$ be the Frobenius endomorphism of a commutative $k$-algebra $R$. Let $\text{Aut}_{\bar g}$ be the subgroup scheme of $G_k$ of inner automorphisms of $\scrC_{\bar g}$. Thus $\text{Aut}_{\bar g}(R)$ is the subgroup of $G_k(R)$ formed by elements $\bar h_R\in G_k(R)$ that commute with the pulls back of $\bar g\bar{\phi}$ and $\bar \vartheta\bar g^{-1}$ to $R$. For instance, in the case of $\bar g\bar{\phi}$ by this commutation we mean that we have an identity $(\bar g\bar{\phi}\otimes 1_R)\circ (\Phi_R^*(\bar h_R))=\bar h_R\circ (\bar g\bar{\phi}\otimes 1_R):\bar M\otimes_{\sigma} R\to \bar M\otimes_k R$ of $R$-linear maps. Here we identify $\bar M\otimes_k {}_{\sigma} R=\bar M\otimes_k {}_{\sigma} k\otimes_k R=\bar M\otimes_k R\otimes_R {}_{\Phi_R} R$. 
As the mentioned commutations can be described by polynomial equations, we get that $\text{Aut}_{\bar g}$ is indeed a subgroup scheme of $G_k$. 

Let $\scrS_{\bar g}\leqslant H_k$ be the stabilizer subgroup scheme of $\bar g\in G(k)$ under $\dbT_{G_k,\sigma}$. Let $\text{Aut}_{\bar g}^{\text{red}}$ and $\scrS_{\bar g}^{\text{red}}$ be the reduced subgroup schemes of $\text{Aut}_{\bar g}$ and $\scrS_{\bar g}$. Let $\text{Aut}^{0\text{red}}_{\bar g}$ and $\scrS^{0\text{red}}_{\bar g}$ be the identity components of $\text{Aut}^{\text{red}}_{\bar g}$ and $\scrS_{\bar g}^{\text{red}}$.

\smallskip
\noindent
{\bf Lemma 5.2.} {\it  {\bf (a)} If $\bar h=(\bar h_1,\bar h_2,\bar h_3)\in\scrS^{\text{red}}_{\bar g}(k)$ is such that $\bar h_1$ and $\bar h_2$ are identity elements, then $\bar h$ is the identity element.

\smallskip
{\bf (b)} The group scheme $\text{Aut}_{\bar g}\cap L_k$ is finite.}
\medskip

{\bf Proof:}
As $\bar h_1$ and $\bar h_2$ are identity elements, we have $\bar g=\dbT_{G_k,\sigma}(\bar h,\bar g)=\bar g\sigma(\bar h_3^{-1})$. Thus $\sigma(\bar h_3^{-1})$ is the identity element $1_{\bar M}$ and therefore $\bar h_3$ is the identity element $1_{\bar M}$ of $N(k)$. Thus (a) holds.

It suffices to prove (b) under the extra assumption that $G=\pmb{GL}_M$. Not to introduce extra notations (by always replacing $g\in G(W(k))$ with some $g\tilde g\in G(W(k))$), to prove (b) we can also assume that $\sigma$ normalizes both $F^1$ and $F^0$. This implies that $\sigma(N)=N$. Let $\bar h_2\in (\text{Aut}_{\bar g}\cap L_k)(k)$. Let $h_2\in L(W(k))$ be a lift of $\bar h_2$. If $g_1:=h_2g\phi(h_2^{-1})$, then $\scrC_{\bar g_1}=\scrC_{\bar g}$. Thus $\bar h_5:=\bar g^{-1}\bar g_1\in \sigma(N)(k)=N(k)$ (cf. Lemma 4.5 (a); referring to (9) with $\bar h_1=1_{\bar M}$, we have $\bar h_5=\sigma(\bar h_3)$). As $\bar h_5=\bar g^{-1}\bar h_2\bar g\sigma(\bar h_2^{-1})$, we have
$$\bar g^{-1}\bar h_2\bar g=\bar h_5\sigma(\bar h_2).\leqno (10)$$ 
Let $I_0:=\{1,\ldots,d\}^2\cup \{d+1,\ldots,r\}^2$. We write $\bar h_2=1_{\bar M}+\sum_{(i_1,i_2)\in I_0}\bar x_{i_1,i_2}\bar e_{i_1,i_2}$ with $\bar x_{i_1,i_2}\in k$, cf. Subsection 4.1 for notations. Let $\bar X$ be the column vector formed by $\bar x_{i_1,i_2}$'s listed is some order. Let $\bar X^{[p]}$ be the column vector obtained from $\bar X$ by raising each entry to power $p$. 

Every element of $\pmb{GL}_M(k)$ can be uniquely written in the form 
$1_{\bar M}+\sum_{i_3,i_4\in\{1,\ldots,r\}} \bar y_{i_3,i_4}\sigma(\bar e_{i_3,i_4})$ with $\bar y_{i_3,i_4}\in k$.
Next for $(i_1,i_2)\in I_0$ we will identify the coefficients of $\sigma(\bar e_{i_1,i_2})$'s in Formula (10). 

As $\bar h_5\in\sigma(N)(k)=N(k)$ and as we have a short exact sequence $1\to N\to NL\to L\to 1$, the coefficient of $\sigma(\bar e_{i_1,i_2})$ for the right hand side of (10) is $1+\bar x_{i_1,i_2}^p$ if $i_1=i_2$ and is $\bar x_{i_1,i_2}^p$ if $i_1\neq i_2$. If $i_1=i_2$ (resp. $i_1\neq i_2$), then the coefficient of $\sigma(\bar e_{i_1,i_2})$ for the left hand side of (10) is $1$ plus a $k$-linear combination of the $\bar x_{i_1',i_2'}$'s with $(i_1',i_2')\in I_0$ (resp. is such a $k$-linear combination). Thus there exists a square matrix $\bar A$ with coefficients in $k$ such that we have $\bar A\bar X=\bar X^{[p]}$. The system $\bar A\bar X=\bar X^{[p]}$ defines a $p^{d^2+(r-d)^2}$ dimensional $k$-algebra and thus has a finite number of solutions. Thus $(\text{Aut}_{\bar g}\cap L_k)(k)$ is finite. Therefore the group scheme $\text{Aut}_{\bar g}\cap L_k$ is finite. 
\endproof

\smallskip
For $\alpha\in\Phi_N$ let $V_{a,\alpha}$ be the vector group scheme defined by $\grg_{\alpha}$; thus $V_{a,\alpha}(W(k))=\grg_{\alpha}$. The next two Propositions show that for each element $w\in W_G$, $\scrS^{0\text{red}}_{\bar g_w}$ and $o_{\bar g_w}$ are ``controlled'' by $\Phi_N^{+w}$ and $\Phi_N^{-w}$ (respectively).

\smallskip
\noindent
{\bf Proposition 5.3.} {\it  Let $w\in W_G$. Let $\alpha\in\Phi_N^{+w}$. Then there exists a curve of $\scrS^{0\text{red}}_{\bar g_w}$ which passes through the identity element and whose tangent space at this element is $\bar\grg_{\alpha}$.}
\medskip

{\bf Proof:}
Let $y_0\in\grg_{\alpha}$. For $i\in\{0,\ldots,1+w_{\alpha}\}$ let 
$$x_i:=(g_w\phi)^i(py_0).$$ 
We have $x_0=py_0$. As $y_0^2=0$, for $i\in\{0,\ldots,1+w_{\alpha}\}$ we have $x_i^2=0$. Thus $(1_M+x_i)(1_M-x_i)=1_M$. As $\alpha\in\Phi_N^{+w}$ we have $x_1,\ldots, x_{w_{\alpha}-1}\in\Lie(L)$ and $x_{w_{\alpha}}\in\grg_{\pi_w^{w_{\alpha}}(\alpha)}\subseteq\Lie(U)$. Also $x_{w_{\alpha}+1}=(g_w\phi)(x_{w_{\alpha}})\in p\Lie(G)$; thus $\bar x_{w_{\alpha}+1}=0$. Let $h_0:=\prod_{i=w_{\alpha}}^0 (1_{M}+x_i)$ be the product obtained using a decreasing order. 

We have $h_0^{-1}=\prod_{i=0}^{w_{\alpha}} (1_M+x_i)^{-1}=\prod_{i=0}^{w_{\alpha}} (1_M-x_i)$. Thus we can write $h_0g_w\phi(h_0^{-1})=g_1g_w$, where $g_1:=\prod_{i=w_{\alpha}}^0 (1_M+x_i)\prod_{i=1}^{1+w_{\alpha}} (1_M-x_i)\in G(W(k))$. As $\bar x_0=\bar x_{w_{\alpha}+1}=0$, we have $$\bar g_1=(\prod_{i=w_{\alpha}}^1 (1_{\bar M}+\bar x_i))(1_{\bar M}+ \bar x_0)(\prod_{i=1}^{w_{\alpha}} (1_{\bar M}-\bar x_i))(1_{\bar M}-\bar x_{w_{\alpha}+1})=1_{\bar M}.$$ 
\indent
Let $h_1:=1_M+x_{w_{\alpha}}\in U(W(k))$, $h_2:=\prod_{i=w_{\alpha}-1}^{1} (1_{M}+x_i)\in L(W(k))$, and $h_3:=h_2(1_M+y_0)h_2^{-1}\in N(W(k))$. We have $h_0=h_1h_2(1_M+x_0)$. Thus $h_0g_w\phi(h_0^{-1})$ mod $p$ is on one hand $\bar g_w$ (as $\bar g_1=1_{\bar M})$ and on the other hand it is 
$\bar h_1\bar h_2\bar g_w\sigma((1_{\bar M}+\bar y_0)^{-1})\sigma(\bar h_2^{-1})=\bar h_1\bar h_2\bar g_w \sigma(\bar h_2^{-1})\sigma(\bar h_3^{-1})=\dbT_{G_k,\sigma}((\bar h_1,\bar h_2,\bar h_3),\bar g_w)$.
Thus $(\bar h_1,\bar h_2,\bar h_3)\in\scrS_{\bar g_w}^{\text{red}}(k)$. 
Therefore the $k$-morphism 
$$c_{\alpha}:{V_{a,\alpha}}_k\to H_k$$ 
that maps $\bar y_0\in {V_{a,\alpha}}_k(k)$ to the element (defined using decreasing order products)
$$(\bar h_1,\bar h_2,\bar h_3)=(1_{\bar M}+\bar x_{w_{\alpha}},\prod_{i=w_{\alpha}-1}^{1} (1_{\bar M}+\bar x_i),(\prod_{i=w_{\alpha}-1}^{1} (1_{\bar M}+\bar x_i))(1_{\bar M}+\bar y_0)(\prod_{i=w_{\alpha}-1}^{1} (1_{\bar M}+\bar x_i))^{-1})\in H_k(k),\leqno (11)$$ factors through $\scrS^{\text{red}}_{\bar g_w}$ and thus (as ${V_{a,\alpha}}_k$ is connected) through $\scrS^{0\text{red}}_{\bar g_w}$. If $\bar y_0$ is multiplied by $\bar\gamma\in\dbG_m(k)$, then $\bar x_i$ gets multiplied by $\bar\gamma^{p^i}$. Thus $\im(c_{\alpha})$ is a curve of $\scrS^{0\text{red}}_{\bar g_w}$ whose tangent space at $1_{\bar M}$ is $\bar\grg_{\alpha}$.\endproof

\smallskip
\noindent
{\bf Proposition 5.4.} {\it   Let $w\in W_G$. Then $\dim(o_{\bar g_w})\ge\dim(G_k)-\dbS(w)$ (see Subsubsection 4.1.1 for $\dbS(w)$).}
\medskip

{\bf Proof:}
 Let $\alpha\in\Phi_N^{-w}$. Let $\tilde\alpha:=\pi_w^{w_{\alpha}}(\alpha)$. Let $y_0\in\grg_{\alpha}$. Let $x_0,x_1,\ldots,x_{w_{\alpha}}$ be obtained from $y_0$ as in the proof of Proposition 5.3. As $\alpha\in\Phi_N^{-w}$, we have $x_1,\ldots,x_{w_{\alpha}-1}\in \Lie(L)$ and $x_{w_{\alpha}}\in\grg_{\tilde\alpha}\subseteq\Lie(N)$. Let $r_{\bar g_w}:P_k\to o_{\bar g_w}\bar g_w^{-1}$ be the $k$-morphism that takes $\bar h_1\bar h_2\in U(k)L(k)=P(k)=P_k(k)$ to $\bar h_1\bar h_2\bar g_w\sigma(\bar h_2^{-1})\bar g_w^{-1}$. The tangent map of $r_{\bar g_w}$ at the identity element has image $\Lie(P_k)$. Thus $\im(r_{\bar g_w})$ has dimension $\dim(P_k)$ and its tangent space at the identity element is $\Lie(P_k)$. The following identity
$$(\prod_{i=0}^{w_{\alpha}-1}(1_M+x_i))(1_M+x_{w_{\alpha}})(g_w\phi)(\prod_{i=0}^{w_{\alpha}-1}(1_M+x_i))^{-1}=(1_M+x_0)g_w\phi=(1_M+py_0)g_w\phi$$ 
implies $<\scrC_{(1_{\bar M}+\bar x_{w_{\alpha}})\bar g_w}>=<\scrC_{\bar g_w}>\in\dbI$. From this and Lemma 5.1 (a) we get $(1_{\bar M}+\bar x_{w_{\alpha}})\bar g_w\in o_{\bar g_w}$. Thus we can speak about the $k$-morphism
$$c_{\alpha}:{V_{a,\tilde\alpha}}_k\to o_{\bar g_w}\bar g_w^{-1}$$ 
that takes $\bar x_{w_{\alpha}}\in {V_{a,\tilde\alpha}}_k(k)$ to $1_{\bar M}+\bar x_{w_{\alpha}}$. The image of the tangent map of $c_{\alpha}$ at the identity element is $\bar\grg_{\tilde\alpha}$. As $o_{\bar g_w}\bar g_w^{-1}$ is smooth, its dimension is the dimension of its tangent space at the identity element. As this tangent space contains $\Lie(P_k)\oplus\bigoplus_{\alpha\in\Phi_N^{-w}} \bar\grg_{\tilde\alpha}$ and as $\Phi_N^{-w}=\Phi_N\setminus\Phi_N^{+w}$ has $\dim(N_k)-\dbS(w)$ elements, we have the following inequality $\dim(o_{\bar g_w})\ge\dim(P_k)+\dim(N_k)-\dbS(w)=\dim(G_k)-\dbS(w)$.\endproof 

\medskip
\noindent
{\bf 5.2. Basic Theorem A (long form)} 

{\bf (a)} {\it There exists a natural finite epimorphism $a_{\bar  g}:\scrS^{\text{red}}_{\bar g}\twoheadrightarrow \text{Aut}^{\text{red}}_{\bar g}$ that induces an isomorphism at the level of $k$-valued points. Thus the four groups $\scrS^{\text{red}}_{\bar g}$, $\text{Aut}^{\text{red}}_{\bar g}$, $\scrS^{0\text{red}}_{\bar g}$, and $\text{Aut}^{0\text{red}}_{\bar g}$ have the same dimension. 

\smallskip
{\bf (b)} The group scheme $\scrS^{0\text{red}}_{\bar g}$ is unipotent and its Lie algebra is commutative. 

\smallskip
{\bf (c)} For $w\in W_G$ we have $\dim(\scrS^{0\text{red}}_{\bar g_w})=\dbS(w)$ and thus
$$\dim(o_{\bar g_w})=\dim(H_k)-\dbS(w)=\dim(G_k)-\dbS(w).\leqno (12)$$
\indent
{\bf (d)} The group $\text{Aut}^{0\text{red}}_{\bar g}$ is unipotent and for $w\in W_G$ we have $\dim(\text{Aut}^{\text{red}}_{\bar g_w})=\dbS(w)$.}

{\bf Proof:}
Let $\bar h=(\bar h_1,\bar h_2,\bar h_3)\in \scrS^{\text{red}}_{\bar g}(k)$. Let $h_{12}\in P(W(k))$ be a lift of  $\bar h_{12}:=\bar h_1\bar h_2\in P(k)$. Let $g_2:=h_{12}g\phi(h_{12}^{-1})\in G(W(k))$. As $\bar h_1\bar h_2\bar g\sigma(\bar h_2^{-1})\sigma(\bar h_3^{-1})=\bar g$, from (9) we get that $\bar g_2=\bar g\sigma(\bar h_3)$. Thus $\scrC_{\bar g}=\scrC_{\bar g_2}$, cf. Lemma 4.5 (a) and the relation $\bar h_3\in N(k)$. Thus $\bar h_{12}\in \text{Aut}^{\text{red}}_{\bar g}(k)$. Due to Formula (3), the rule
$$\bar h=(\bar h_1,\bar h_2,\bar h_3)\mapsto\bar h_{12}=\bar h_1\bar h_2\;\;\text{defines}\;\;\text{a}\;\;\text{homomorphism}\;\; a_{\bar  g}:\scrS^{\text{red}}_{\bar g}\to \text{Aut}^{\text{red}}_{\bar g}.$$ 
As $\Ker(a_{\bar  g})(k)$ is trivial (cf. Lemma 5.2 (a)), $\Ker(a_{\bar  g})$ is a finite, flat group scheme. Thus $a_{\bar  g}$ is finite. 

Next we work conversely. If $\bar h_{12}=\bar h_1\bar h_2\in\text{Aut}^{\text{red}}_{\bar g}(k)\leqslant P(k)$ and if $h_{12}$ and $g_2$ are as above, then $\scrC_{\bar g}=\scrC_{\bar g_2}$ and thus $\bar h_3:=\sigma^{-1}(\bar g^{-1}\bar g_2)\in N(k)$ (cf. Lemma 4.5 (a)). As $\bar g=\bar h_1\bar h_2\bar g\sigma(\bar h_2^{-1})\sigma(\bar h_3^{-1})$, we have $(\bar h_1,\bar h_2,\bar h_3)\in\scrS^{\text{red}}_{\bar g}(k)$ and thus $a_{\bar  g}(k)$ maps $(\bar h_1,\bar h_2,\bar h_3)$ to $\bar h_{12}$. Thus the homomorphism $a_{\bar  g}(k):\scrS^{\text{red}}_{\bar g}(k)\to \text{Aut}^{\text{red}}_{\bar g}(k)$ is onto. Thus $a_{\bar g}$ is an epimorphism. As $a_{\bar  g}(k)$ is also injective, it is an isomorphism. Thus (a) holds.

As $\dbT_{G_k,\sigma}$ is a restriction of the analogue action $\dbT_{\pmb{GL}_{\bar M},\sigma}$ associated to the triple $(M,\phi,\mu:\dbG_m\to \pmb{GL}_M)$, it suffices to prove (b) under the extra assumption $G=\pmb{GL}_M$. Based on Corollary 4.6, we can assume $g$ is a $g_w$ element for some $w\in W_G$. Thus for the rest of the proofs of (b) and (c) we will work with a $g_w$ element and we will not assume anymore that $G$ is $\pmb{GL}_M$. We recall $\dbS(w)$ is the number of elements of $\Phi_N^{+w}$. 

We have $\oplus_{\alpha\in\Phi_N^{+w}} \bar\grg_{\alpha}\subseteq \Lie(\scrS^{0\text{red}}_{\bar g_w})$, cf. Proposition 5.3. Thus $\dim_k(\Lie(\scrS^{0\text{red}}_{\bar g_w}))\ge \dbS(w)$. As the reduced group $\scrS^{0\text{red}}_{\bar g_w}$ is smooth, we have $\dim(\scrS^{0\text{red}}_{\bar g_w})=\dim_k(\Lie(\scrS^{0\text{red}}_{\bar g_w}))\ge \dbS(w)$. As we have
$$\dim(o_{\bar g_w})=\dim(H_k)-\dim(\scrS_{\bar g}^{\text{red}})=\dim(H_k)-\dim(\scrS^{0\text{red}}_{\bar g_w})=\dim(G_k)-\dim(\scrS^{0\text{red}}_{\bar g_w}),$$ 
from Proposition 5.4 we get $\dim(\scrS^{0\text{red}}_{\bar g_w})\le \dbS(w)$. Thus (12) holds. By reasons of dimensions we also get
 $$\Lie(\scrS^{0\text{red}}_{\bar g_w})=\oplus_{\alpha\in\Phi_N^{+w}} \bar\grg_{\alpha}\subseteq\Lie(N_k).\leqno (13)$$ 
Thus the Lie subalgebra $\Lie(\scrS^{0\text{red}}_{\bar g_w})$ of $\End_k(\bar M)$ is abelian and all its elements have square $0$. In particular, $\Lie(\scrS^{0\text{red}}_{\bar g_w})$ contains no non-zero semisimple element of $\End_k(\bar M)$. Thus $\scrS^{0\text{red}}_{\bar g_w}$ contains no non-trivial torus and therefore it is a unipotent group. Thus (b) and (c) hold. Part (d) is a direct consequence of (a), (b), and (c).\endproof 

\smallskip
\noindent
{\bf Remark 5.5.}
If $\alpha\in\Phi_N^{+w}$, then the tangent map of $a_{\bar  g_w}\circ  c_{\alpha}$ at the identity element is trivial. From this and (13) we get that $a_{\bar g_w}$ is the composite of the Frobenius endomorphism of $\scrS_{\bar g_w}^{\text{red}}$ with a finite homomorphism $a_{\bar g_w}':\scrS_{\bar g_w}^{\text{red}}\to\text{Aut}^{\text{red}}_{\bar g_w}$ which at the level of $k$-valued points induces an isomorphism. It is easy to see that $a_{\bar g_w}'$ induces at the level of Lie algebras an isomorphism. The last two sentences imply that $a_{\bar g_w}'$ is an isomorphism.

\smallskip
\noindent
{\bf Example 5.6.}
See Subsection 4.1 for $B$ and $B^{\text{opp}}$. As $\sigma(B)$ is a Borel subgroup scheme of $G$ that contains $T$, there exists a unique element $w\in W_G$ such that $g_w\sigma(B)g_w^{-1}=B$. Thus $g_w\phi$ normalizes $\Lie(B)[{1\over p}]$. As $g_w\phi$ also normalizes $\Lie(T)$, we get that $g_w\phi$ normalizes $\Lie(B^{\text{opp}})[{1\over p}]$ as well. Thus, as $N\leqslant B^{\text{opp}}$ and as $B^{\text{opp}}\cap U$ is the trivial group scheme, we conclude that the set $\Phi_N^{+w}$ is empty. Therefore $\dbS(w)=0$. This implies that $\dim(o_{\bar g_w})=\dim(G_k)$, cf. (12). Thus $o_{\bar g_w}$ is a dense, open subscheme of $G_k$.

\smallskip
\noindent
{\bf Example 5.7.} Let $w$ be as in Example 5.6. Let $G_0$ be a simple factor of $G^{\ad}$ such that the image of $\mu$ in it is non-trivial. Let $i_0\in\dbN$ be the smallest number such that the intersection $\sigma_{g_w}^{i_0}(N)\cap G_0$ is non-trivial (this intersection makes sense as $N$ is also a unipotent subgroup scheme of $G^{\ad}$). We assume that for each such simple factor $G_0$ of $G^{\ad}$, there exists $w_{00}\in W_{G_0}\vartriangleleft W_{G^{\ad}}=W_G$ such that $g_{w_{00}}(\sigma_{g_w}^{i_0}(N)\cap G_0)=U\cap G_0$. Let $w_0\in W_G$ be the product of all these $w_{00}$ elements. Let $w_1:=w_0w\in W_G$. We have $\sigma_{g_{w_1}}^{i_0}(N)\cap G_0=U\cap G_0$. Due to the minimal property of $i_0$, we have $\sigma_{g_{w_1}}^{i_1}(N)\cap G_0\leqslant L$ for $i_1\in\{1,\ldots,i_0-1\}$. From the last two sentences we get $\Phi_N^{+w_1}=\Phi_N$. Thus $\dbS(w_1)=\dim(N_k)$. From (12) we get 
$$\dim(o_{\bar g_{w_1}})=\dim(G_k)-\dim(N_k)=\dim(P_k).$$ 
\indent
One can check that the Weyl elements $w_{00}$ exist if all non-trivial simple factors of $(\Lie(G),\phi)$ are of $B_n$, $C_n$, or $D_n^{\dbR}$ type (see [2, plates II to IV, element $w_0$ of item (VIII)]).

\smallskip
\noindent
{\bf Remark 5.8.} Let $o'_{\bar g_w}$ be the integral, locally closed subscheme of $G_k$ whose $k$-valued points are of the form 
$$\bar h_1\bar h_2\bar h^{-w}\bar g_w\sigma(\bar h_2^{-1}),\;\;\text{where}\;\;\bar h_1\in U_k(k),\;\;\bar h_2\in L_k(k),\;\;\text{and}\;\;\bar h^{-w}\in N^{-w}_k(k);$$ 
here $N^{-w}:=\prod_{\alpha\in\Phi_N^{-w}} \dbG_{a,{\pi_w^{w_{\alpha}}(\alpha)}}\leqslant N$. The proof of Proposition 5.4 can be easily adapted to get that $o'_{\bar g_w}$ is a subscheme of $o_{\bar g_w}$ which has dimension $\dim(G_k)-\dbS(w)$. Thus $o'_{\bar g_w}$ is an open subscheme of $o_{\bar g_w}$, cf. (12).

\bigskip
\noindent
{\boldsectionfont \S6. On the number and the dimensions of orbits: a classical approach}

\bigskip
Subsection 6.1 reviews the classical analogue $\dbT^{\text{cl}}_{G_k}$ of $\dbT_{G_k,\sigma}$. Theorem 6.1 presents two results on the number of orbits of $\dbT_{G_k,\sigma}$ and of their dimensions. The Section ends with an example and three remarks. We use the notations of Subsections 4.1 and 5.1. In particular, we have $T\leqslant B\leqslant P\leqslant G$. 
 
\medskip
\noindent
{\bf 6.1. The classical context} 

For the Bruhat decompositions to be recalled here in the context of $G_k$, we refer to Subsection 2.1, [1, Ch. IV, Sect. 14], and to their axiomatized and abstract generalization formalized in [2, Ch. IV, Sect. 2]. Let $U^B$ be the unipotent radical of $B$; it contains the unipotent radical $U$ of $P$. Let $\Phi_{U^B}:=\{\alpha\in\Phi|\grg_{\alpha}\subseteq\Lie(U^B)\}$. Let $H_k^{\text{cl}}:=P_k\times_k U^B_k$. Next we recall properties of the classical action
$$\dbT_{G_k}^{\text{cl}}:H_k^{\text{cl}}\times_k G_k\to G_k$$
defined by the rule: if $\bar h_{12}\in P_k(k)$ and $\bar h_3\in U^B_k(k)$, then $\dbT^{\text{cl}}_{G_k}((\bar h_{12},\bar h_3),\bar g):=\bar h_{12}\bar g\bar h_3^{-1}$. 

\smallskip
{\bf (i)} The orbits of $\dbT_{G_k}^{\text{cl}}$ are parameterized by elements of $W_P\backslash W_G$ (cf. [2, Ch. IV, Sect. 2, Subsect. 6, Rm. of p. 22] or [6, Vol. III, Exp. XXVI, Subsubsect. 4.5.5]). Let $w\in W_G$. Let $U_k^{B,w,P}:=U_k^B\cap\bar g_w^{-1}P_k\bar g_w$ and $U_k^{B,w,N}:=U_k^B\cap\bar g_w^{-1}N_k\bar g_w$. These two group schemes are normalized by $T_k$, are directly span by the ${\dbG_{a,\alpha}}_k$'s with $\alpha\in S(w)^{\text{cl}}:=\Phi_{U^B}\cap w^{-1}\Phi_{P}w$ and with $\alpha\in \Phi_{U^B}\cap w^{-1}\Phi_{N}w$ (respectively), and their product $k$-morphism $U_k^{B,w,P}\times_k U_k^{B,w,N}\to U_k^B$ is an isomorphism (cf. [1, Ch. IV, Prop. 14.4 (2)]). In particular, the groups $U_k^{B,w,P}$ and $U_k^{B,w,N}$ are smooth and connected. The orbit of $\bar g_w$ under $\dbT_{G_k}^{\text{cl}}$ is the affine variety
$$P_k\bar g_wU^B_k=P_k(\bar g_wU^B_k\bar g_w^{-1})\bar g_w=(P_k\bar g_wU_k^{B,w,N}\bar g_w^{-1})\bar g_w.\leqno (14)$$
\indent
{\bf (ii)} The smooth group $\bar g_wU_k^{B,w,P}\bar g_w^{-1}$ is the stabilizer subgroup of $\bar g_w$ under $\dbT_{G_k}^{\text{cl}}$. The classical analogue $\dbS^{\text{cl}}:W_G\to\dbN\cup\{0\}$ of the function $\dbS:W_G\to\dbN\cup\{0\}$, is defined by:
$$\dbS^{\text{cl}}(w):=\dim(\bar g_wU_k^{B,w,P}\bar g_w^{-1})=\text{the}\;\text{number}\;\text{of}\;\text{elements}\;\text{of}\;\text{the}\;\text{set}\;S(w)^{\text{cl}}.$$ \indent
{\bf (iii)} We have $\dbS^{\text{cl}}(w)=\dim(U^B_k)$ if and only if $[w]\in W_P\backslash W_G$ is the class of the identity element of $W_G$. Also, we have $\dbS^{\text{cl}}(w)=\dim(U^B_k)$ if and only if the orbit of $\bar g_w$ under $\dbT_{G_k}^{\text{cl}}$ has dimension equal to (equivalently at most equal to) $\dim(H_k^{\text{cl}})-\dim(U^B_k)=\dim(P_k)$.

\smallskip
{\bf (iv)} Let $w_0\in W_G$ be such that $B^{\text{opp}}=g_{w_0}Bg_{w_0}^{-1}\leqslant Q$. Let
$$d_u:=\dim(U^B_k)-\dim(N_k)=\dim(U^B_k)-\dim(U_k).$$
We have $\dbS^{\text{cl}}(w)\ge d_u$. The equality holds if and only if $N_k\leqslant \bar g_wU^B_k\bar g_w^{-1}$ and thus if and only if the intersection $Q_k\cap \bar g_wB_k\bar g_w^{-1}$ is a Borel subgroup of $G_k$ (cf. [1, Ch. IV, Prop. 14.22 (i)]) and thus if and only if $\bar g_wB_k\bar g_w^{-1}\leqslant Q_k$. Therefore we have $\dbS^{\text{cl}}(w)=d_u$ if and only if $[w]=[w_0]\in W_P\backslash W_G$.  

\smallskip
{\bf (v)} Let $o_1^{\text{cl}},\ldots,o_{[W_G:W_P]}^{\text{cl}}$ be the orbits of $\dbT_{G_k}^{\text{cl}}$. For $i\in\{1,\ldots,[W_G:W_P]\}$ let $d_i^{\text{cl}}$ be the dimension of the stabilizer subgroup of any $k$-valued point of $o_i^{\text{cl}}$ under $\dbT_{G_k}^{\text{cl}}$.

\smallskip
\noindent
{\bf Theorem 6.1.} {\bf (a)} {\it Suppose $\dbT_{G_k,\sigma}$ has a finite number of orbits $o_1,\ldots,o_s$. For $i\in\{1,\ldots,s\}$ we fix a point $\bar g_i\in o_i(k)$; let $d_i:=\dim(\scrS_{\bar g_i})$. Then $s=[W_G:W_P]$ and up to a reindexing of the orbits $o_1,\ldots,o_s$, we have 
$$d_i=d_i^{\text{cl}}-d_u\;\;\;\forall i\in\{1,\ldots,s\}.\leqno (15)$$ 
Equality (15) is equivalent to the fact that for all $i\in\{1,\ldots,s\}$ we have $\dim(o_i)=\dim(o_i^{\text{cl}})$.

\smallskip
{\bf (b)} Suppose there exist distinct orbits $o_1,\ldots,o_{[W_G:W_P]}$ of $\dbT_{G_k,\sigma}$ such that for each $i\in\{1,\ldots,[W_G:W_P]\}$ we have $d_i=d_i^{\text{cl}}-d_u$. Then $\dbT_{G_k,\sigma}$ has precisely $[W_G:W_P]$ orbits.}
\medskip

{\bf Proof:}
We prove (a). We recall that the pair $(G,T)$ has a natural $\dbZ_p$-structure $(G_{\dbZ_p},T_{\dbZ_p})$, cf. Subsubsection 4.2.1. Let $q_0\in\dbN$ be such that the torus $T_{W(\dbF_{p^{q_0}})}$ is split. The cocharacter $\mu:\dbG_m\to T$ is definable over $\Spec(W(\dbF_{p^{q_0}}))$. Thus $P$, $U$, $Q$, and $N$ are definable over $\Spec(W(\dbF_{p^{q_0}}))$. As $T_{W(\dbF_{p^{q_0}})}$ is split and as $B$ contains $T$, $B$ is also definable over $\Spec(W(\dbF_{p^{q_0}}))$. Thus both groups $H_k^{\text{cl}}$ and $H_k$ are definable over $\dbF_{p^{q_0}}$. This implies that the two actions $\dbT_{G_k}^{\text{cl}}$ and $\dbT_{G_k,\sigma}$ are definable over $\dbF_{p^{q_0}}$; the resulting two actions over the algebraic closure $\dbF$ of $\dbF_{p^{q_0}}$ in $k$ will have also $[W_G:W_P]$ and $s$ (respectively) number of orbits.  

Thus there exists $q\in q_0\dbN$ such that $\dbT_{G_k}^{\text{cl}}$ and $\dbT_{G_k,\sigma}$ and all their orbits $o_1^{\text{cl}},\ldots,o_{[W_G:W_P]}^{\text{cl}}$, $o_1$,$\ldots$,$o_s$, are defined over $\dbF_{p^q}$. By enlarging $q$ we can assume that $\bar g_1,\ldots,\bar g_s$ are also defined over $\dbF_{p^q}$. Let $H_{\dbF_{p^q}}$, $\scrS_{\bar g_i\dbF_{p^q}}$, $\scrS^{\text{red}}_{\bar g_i\dbF_{p^q}}$, $\scrS^{0\text{red}}_{\bar g_i\dbF_{p^q}}$, and $o_{i\dbF_{p^q}}$ be the natural $\dbF_{p^q}$-structures of $H_k$, $\scrS_{\bar g_i}$, $\scrS^{\text{red}}_{\bar g_i}$, $\scrS^{0\text{red}}_{\bar g_i}$, and $o_i$ (respectively). 

Let $q_1\in q\dbN$. We have a finite epimorphism $H_{\dbF_{p^q}}/\scrS^{\text{red}}_{\bar g_i\dbF_{p^q}}\twoheadrightarrow o_{i\dbF_{p^q}}$ that is radicial (i.e., universal injective on geometric points). Thus the map $(H_{\dbF_{p^q}}/\scrS^{\text{red}}_{\bar g_i\dbF_{p^q}})(\dbF_{p^{q_1}})\to o_{i\dbF_{p^q}}(\dbF_{p^{q_1}})$ is a bijection. 
As $\scrS^{0\text{red}}_{\bar g_i\dbF_{p^q}}$ is a unipotent group (cf. (b) of Subsection 5.2), from Lemma 2.5 we get that $o_{i\dbF_{p^q}}(\dbF_{p^{q_1}})$ has the same number of elements $e_i(q_1)$ as $(H_{\dbF_{p^q}}/\scrS^{0\text{red}}_{\bar g_i\dbF_{p^q}})(\dbF_{p^{q_1}})$ an thus (cf. Lang theorem) as the quotient set $H_{\dbF_{p^q}}(\dbF_{p^{q_1}})/\scrS^{0\text{red}}_{\bar g_i\dbF_{p^q}}(\dbF_{p^{q_1}})$. 

Let $g(q_1)$, $h(q_1)$, and $h^{\text{cl}}(q_1)$ be the orders of the groups $G_{\dbF_{p^q}}(\dbF_{p^{q_1}})$, $H_{\dbF_{p^q}}(\dbF_{p^{q_1}})$, and $H_{\dbF_{p^q}}^{\text{cl}}(\dbF_{p^{q_1}})$ (respectively). As the unipotent group $\scrS^{0\text{red}}_{\bar g_i\dbF_{p^q}}$ has a normal series whose factors are $\dbG_a$ groups (cf. [6, Vol. II, Exp. XVII, Cor. 4.1.3]), the scheme $\scrS^{0\text{red}}_{\bar g_i\dbF_{p^q}}$ is the spectrum of a polynomial $\dbF_{p^q}$-algebra in $d_i$ variables; thus the order of $\scrS^{0\text{red}}_{\bar g_i\dbF_{p^q}}(\dbF_{p^{q_1}})$ is $p^{q_1d_i}$. Therefore $e_i(q_1)=h(q_1)p^{-q_1d_i}$. As $G_{i\dbF_{p^q}}(\dbF_{p^{q_1}})$ is a disjoint union $\sqcup_{i=1}^s o_{i\dbF_{p^q}}(\dbF_{p^{q_1}})$, we get 
$$g(q_1)=h(q_1)\sum_{i=1}^s p^{-q_1d_i}.\leqno (16)$$ 
\indent
The group $H_k$ (resp. $H_k^{\text{cl}}$) is the extension of $L_k$ by $U_k\times_k N_k$ (resp. by $U_k\times_k U^B_k$), cf. (4) (resp. cf. Subsection 6.1). As $U^B_k$ and $N_k$ are unipotent groups and as we have $\dim(U^B_k)=\dim(N_k)+d_u$, a simple computation shows that $h^{\text{cl}}(q_1)=p^{q_1d_u}h(q_1)$. Working with $\dbT_{G_k}^{\text{cl}}$ instead of $\dbT_{G_k,\sigma}$, we get the classical analogue of (16)
$$g(q_1)=h^{\text{cl}}(q_1)\sum_{i=1}^{[W_G:W_P]} p^{-q_1d_i^{\text{cl}}}.\leqno (17)$$ 
As $h^{\text{cl}}(q_1)=p^{q_1d_u}h(q_1)$, from (16) and (17) we get that
$$p^{q_1d_{\text{m}}}g(q_1)/h^{\text{cl}}(q_1)=\sum_{i=1}^{[W_G:W_P]} p^{q_1(d_{\text{m}}-d_i^{\text{cl}})}=\sum_{i=1}^s p^{q_1(d_{\text{m}}-d_i-d_u)}, \leqno (18)$$ 
where $d_{\text{m}}:=\max\{d_u+d_i,d_j^{\text{cl}}|1\le i\le s,1\le j\le [W_G:W_P]\}$. We take $q_1>\max\{s,[W_G:W_P]\}$. For $a\in\dbN$, $p^{q_1a}$ is greater than each sum of $q_1-1$ terms of the form $p^{q_1b}$ with $b\in\{0,\ldots, a-1\}$. From this and (18), we get that $s=[W_G:W_P]$ and that (15) holds after a reindexing of the orbits $o_1,\ldots,o_s$. We compute
$$\dim(o_i)=\dim(G_k)-d_i=\dim(P_k)+\dim(N_k)+d_u-d_i^{\text{cl}}=\dim(P_k)+\dim(U_k^B)-d_i^{\text{cl}}=$$
$$\dim(H_k^{\text{cl}})-d_i^{\text{cl}}=\dim(o_i^{\text{cl}}).$$ 
\indent
Thus (a) holds. We prove (b). Our hypotheses imply that by taking $s=[W_G:W_P]$, the right hand sides of (16) and (17) are equal. Thus identity (16) holds. Therefore for each $q_1\in q\dbN$ we have a disjoint union $G_{\dbF_{p^q}}(\dbF_{p^{q_1}})=\sqcup_{i=1}^{[W_G:W_P]} o_{i\dbF_{p^q}}(\dbF_{p^{q_1}})$. 
This implies that $\dbT_{G_k,\sigma}$ has precisely $[W_G:W_P]$ orbits.\endproof

\smallskip
\noindent
{\bf Corollary 6.2.} {\it 
If the statement 1.4 (a) holds, then statement 1.4 (b) also holds.}
\medskip

{\bf Proof:}
We have $\dbI=\{<\scrC_{\bar g_w}>|w\in W_G\}$, cf. Lemma 4.5 (d). Thus the set $\dbI$ is finite. Therefore the set $\dbO$ is also finite, cf. Lemma 5.1 (a). Thus $\dbO$ has $[W_G:W_P]$ elements, cf. Theorem 6.1 (a). Therefore $\dbI$ also has  $[W_G:W_P]$ elements, cf. Lemma 5.1 (a). Thus the statement 1.4 (b) holds, cf. Lemma 4.5 (e).\endproof

\smallskip
\noindent
{\bf Example 6.3.}
Let $n\in\dbN$. Let $G$ be such that $G^{\ad}$ is absolutely simple of $B_n$ Lie type. Thus $G^{\ad}_{\dbZ_p}$ is the $\pmb{SO}$-group of the quadratic form $x_0^2+x_1x_{1+n}+x_2x_{2+n}+\cdots+x_{n}x_{2n}$ on $\dbZ_p^{2n+1}$ (if $p=2$, then $G^{\ad}_{\dbF_p}$ is semisimple cf. [1, Ch. V, Subsect. 23.6]). Let $\{b_0,\ldots,b_{2n}\}$ be the standard $\dbZ_p$-basis for $\dbZ_p^{2n+1}$. To ease notations we will assume that: (i) $T_{\dbZ_p}$ is split, (ii) the image $T'_{\dbZ_p}$ of $T_{\dbZ_p}$ in $G^{\ad}_{\dbZ_p}$ normalizes each $\dbZ_pb_i$ with $i\in\{1,\ldots, 2n\}$ and fixes $\dbZ_pb_0$, and (iii) the cocharacter of $T'_{\dbZ_p}$ defined by $\mu$ fixes $\{b_i|2\le i\le 2n,i\neq n+1\}$, acts as the inverse of the identical character of $\dbG_m$ on $\dbZ_pb_{n+1}$, and acts as the identical character of $\dbG_m$ on $\dbZ_pb_1$. For $w\in W_G$ we denote also by $w$ the permutation of $\{1,\ldots,2n\}$ such that $g_w(b_i)\in W(k)b_{w(i)}$. The resulting map $W_G\to S_{2n}$ is a monomorphism that identifies $W_G$ with the subgroup of the symmetric group $S_{2n}$ of $\{1,\ldots,2n\}$  which permutes the subsets $\{1,n+1\}$, $\{2,n+2\},\ldots,\{n,2n\}$ (see [2, plate II]). Thus $W_G$ has order $2^nn!$. Moreover $W_P$ is the  subgroup of $W_G$ that fixes $1$ and $n+1$ and thus its order is $2^{n-1}(n-1)!$. Thus $[W_G:W_P]=2n$.

Let $\{b_{i_1,i_2}|0\le i_1,i_2\le 2n\}$ be the standard $\dbZ_p$-basis for $\End_{\dbZ_p}(\dbZ_p^{2n+1})$. It is known (cf. [18, Ch. VI, Sect. 25, Type B of p. 358]) that $\Lie(N)$ and $\Lie(U)$ are $W(k)$-generated by 
$$\{b_{1,i}-b_{n+i,n+1},b_{1,n+i}-b_{i,n+1}|2\le i\le n\}\cup\{b_{0,n+1}-2b_{1,0}\}\;\;\text{and}$$
$$\{b_{n+1,i}-b_{i+n,1},b_{n+1,n+i}-b_{i,1}|2\le i\le n\}\cup\{b_{0,1}-2b_{n+1,0}\}$$ 
(respectively) and thus have ranks $2n-1$. To each non-zero element $b_{i_1,j_1}-\eta b_{i_2,j_2}$ (with $\eta\in\{1,2\}$) we associate the pair $(i_1,j_1)$ and its attached quadruple $(i_1,j_1,i_2,j_2)$. Accordingly we define 
$$I_N:=\{(1,i),(1,n+i)|2\le i\le n\}\cup\{(0,n+1)\}\;\;\text{and}\;\;I_U:=\{(n+1,i),(n+1,n+i)|2\le i\le n\}\cup\{(0,1)\}.$$ 
For $w\in W_G$, $\dbS(w)$ is the number of elements of the subset $S(w)$ of $I_N$ formed by pairs $(a,b)\in I_N$ whose attached quadruples $(a,b,c,d)$ have the property that the smallest number $l\in\dbN$ for which the intersection 
$$\{(w^l(a),w^l(b)),(w^l(c),w^l(d))\}\cap (I_N\cup I_U)$$ is non-empty, is such that $\{(w^l(a),w^l(b)),(w^l(c),w^l(d))\}\cap I_U\neq\emptyset$.  

For $(j,\varepsilon)\in\{1,\ldots,n\}\times\{-1,1\}$, let $w_{j,\varepsilon}\in W_G\subset S_{2n}$ be such that: it fixes $j+1,\ldots,n$, $n+j+1,\ldots,2n$, maps $s$ to $s+1$ if $1\le s\le j-1$, and $(j,n+j)$ to $(1,n+1)$ (resp. $(n+1,1)$) if $\varepsilon=1$ (resp. $\varepsilon=-1$). We have 
$$S(w_{j,1})=\{(1,n+i)|2\le i\le j\}\;\;\text{and}\;\;S(w_{j,-1})=\{(1,i)|2\le i\le n\}\cup\{(1,n+i)|2\le i\le j\}\cup\{(0,n+1)\}.$$ 
Thus $\dbS(w_{j,1})=j-1$ and $\dbS(w_{j,-1})=(n-1)+(j-1)+1=n+j-1$. 

The reduced stabilizer subgroup of a point of $o_{\bar g_{w_{j,\varepsilon}}}$ under $\dbT_{G_k,\sigma}$ has dimension $\dbS(w_{j,\varepsilon})$ (cf. (12)); this dimension $\dbS(w_{j,\varepsilon})$ determines uniquely $w_{j,\varepsilon}$. Thus the orbits $o_{\bar g_{w_{j,\varepsilon}}}$'s are distinct and we can index them $o_1,\ldots,o_{2n}$ in such a way that for $i\in\{1,\ldots,2n\}$ we have $d_i=i-1$. As $\Phi$ has $n^2$ positive roots (cf. [2, plate II]), we have $d_u=n^2-(2n-1)=(n-1)^2$ and $\dim(G_k)=2n^2+n$. It is well known that we can assume that for $i\in\{1,\ldots,2n\}$ we have $d_i^{\text{cl}}=d_u+d_i=(n-1)^2+i-1$ (i.e.,  $\dim(G_k)-\dim(o_i^{\text{cl}})=i-1$ and thus that $\dim(o_i^{\text{cl}})=\dim(o_i)=2n^2+n-i+1$; this is implicitly reproved in this paper once Basic Theorem C is proved, cf. Theorem 6.1 (a)). Thus $o_1,\ldots,o_{2n}$ are the only orbits of $\dbT_{G_k,\sigma}$, cf. Theorem 6.1 (b).

\smallskip
\noindent
{\bf Remark 6.4.} {\bf (a)} See Remark 5.8 for $o'_{\bar g_w}$. For $w_1$, $w_2\in W_G$, we have $o_{\bar g_{w_1}}\le o_{\bar g_{w_2}}$ if and only if the Zariski closure of $o_{\bar g_{w_1}}'$ in $G_k$ contains $o_{\bar g_{w_2}}'$. Even more, it seems to us that the counting arguments of the proof of Theorem 6.1 can be adapted to show that the following presumed analogue $o'_{\bar g_w}=o_{\bar g_w}$ 
of (14) holds. 

\smallskip
{\bf (b)} We assume that: (i) $\dbO=\{o_{\bar g_w}|w\in W_G\}$ and (ii) the association $w\mapsto o_{\bar g_w}$ defines a surjective map $W_G\twoheadrightarrow\dbO$ whose pulls back of elements of $\dbO$, have a number of elements equal to the order of $W_P$ (in Corollary 11.1 (b) and Lemma 9.4 we check that these two assumptions always hold). Then Formula (15) says: there exists a natural bijection $\pi_G:W_G\arrowsim W_G$ such that the difference of functions $\dbS^{\text{cl}}\circ \pi_G-\dbS$ is the constant function $d_u$ on $W_G$. From this we get that $\dbS(w)=\dim(N_k)-l([\pi_G(w)])$. Here $l([w]):=\text{min}\{l(\tilde w)|\tilde w\in [w]\}$, where $l(\tilde w)$ is the length of $\tilde w\in W_G$ (with respect to the bases of $\Phi$ associated to $B$) as defined in the theory of Coxeter groups (see [2, Ch. IV, Sect. 1]); we have $\dbS^{\text{cl}}(w)=\dim(U_k^B)-l([w])$. 

\smallskip
{\bf (c)} Suppose $\dbT_{G_k,\sigma}$ has a finite number of orbits. From the property 6.1 (iii) and Theorem 6.1 (a) we get that $\dbT_{G_k,\sigma}$ has a unique orbit of dimension $\dim(P_k)$ (equivalently of dimension at most $\dim(P_k)$).

\bigskip
\noindent
{\boldsectionfont \S7. Zero spaces}

\bigskip
In Subsection 7.1 we define the zero space of $\scrC_{\bar g}$; examples and elementary properties of them are listed until Theorem 7.5. Theorem 7.5 lists formulas that pertains to the zero spaces of $\scrC_{\bar g_w}$'s. Theorems 7.6 and 7.8 present two results on the existence of connected, smooth subgroups of $G_k$ whose Lie algebras are related to zero spaces of $\scrC_{\bar g_w}$'s. We will use the notations of Subsections 4.1, 4.1.1, (beginning of) 4.2, and 4.2.1. 

\medskip
\noindent
{\bf 7.1. Basic constructions}
 
Let $g\in G(W(k))$. Though $F^1$ and $P$ depend on $\mu:\dbG_m\to G$, the direct summand $F^1/pF^1$ of $\bar M$ is equal to $(g\phi)^{-1}(pM)/pM$ and therefore both it and $P_k$ are uniquely determined by $\bar{\phi}$. Thus the following Lie subalgebra $\Lie(P)+p\Lie(G)$ of $\Lie(G)$ does not depend on either $\mu$ or $g$. We consider the $\sigma$-linear Lie monomorphism
$$L_g:\Lie(P)+p\Lie(G)\hookrightarrow\Lie(G)$$ 
that maps $x\in \Lie(P)+p\Lie(G)$ to $(g\phi)(x)=g\circ\phi\circ x\circ\phi^{-1}\circ g^{-1}\in\Lie(G)$. Let 
$$I_G:\Lie(P)+p\Lie(G)\hookrightarrow\Lie(G)$$
be the natural inclusion. Let $\bar L_g$ and $\bar I_G$ be the reductions mod $p$ of $L_g$ and $I_G$. Let 
$$\tilde\grz_{\bar g}:=\{x\in (\Lie(P)+p\Lie(G))/p(\Lie(P)+p\Lie(G))|\bar L_g(x)=\bar I_G(x)\}$$
and $$\grz_{\bar g}:=\im(\tilde\grz_{\bar g}\to\Lie(G_k))\subseteq \Lie(P_k).$$ 
\noindent
Let $\grw_{\bar g}$ be the $k$-span of $\grz_{\bar g}$ inside $\Lie(P_k)$. Let $\gry^0_{\bar g}$ be the maximal $k$-linear subspace of $\grz_{\bar g}$. We have inclusions
$$\gry^0_{\bar g}\subseteq\grz_{\bar g}\subseteq\grw_{\bar g}\subseteq\Lie(P_k).$$  
 
\smallskip
\noindent
{\bf Lemma 7.1.} {\it  {\bf (a)} When $g$ is replaced by $hg\phi(h^{-1})$ with $h\in\dbL$, $\grz_{\bar g}$ is replaced by its inner conjugate under $\bar h$. 

\smallskip
{\bf (b)} If $g_1\in G(W(k))$ is such that $\scrC_{\bar g}=\scrC_{\bar g_1}$, then $\grz_{\bar g}=\grz_{\bar g_1}$. 

\smallskip
{\bf (c)} The group $\text{Aut}_{\bar g}(k)$ normalizes $\grz_{\bar g}$ (and thus also $\grw_{\bar g}$ and $\gry^0_{\bar g}$).}
\medskip

{\bf Proof:}
Part (a) is obvious as $L_{hg\phi(h^{-1})}=hL_g h^{-1}$. We check (b). The subset $\grz_{\bar g}$ of $\Lie(P_k)$ depends only on $\bar g$ as $L_{\bar g}$ does. Up to a multiplication of $g_1$ by an element in $\Ker(G(W(k))\to G(k))$, we can assume $g_1$ is $hg\phi(h^{-1})$, where $h\in\dbL$ is such that $\bar h=1_{\bar M}$ (cf. proof of Lemma 4.5 (b)). Thus $\grz_{\bar g}=\bar h(\grz_{\bar g})=\grz_{\bar g_1}$, cf. (a) for the last identity. Part (c) follows from (a) and (b).\endproof

\smallskip
The definition of $\tilde\grz_{\bar g}$ and Lemma 7.1 (b) explain the following definition.

\smallskip
\noindent
{\bf Definition 7.2.} We call $\grz_{\bar g}$, $\grw_{\bar g}$, and $\gry^0_{\bar g}$ the {\it zero space}, the {\it spanned zero space}, and the {\it linear zero space} (respectively) of $\scrC_{\bar g}$.

\smallskip
\noindent
{\bf Lemma 7.3.} {\it  {\bf (a)} The subset $\grz_{\bar g}$ of $\Lie(P_k)$ is an $\dbF_p$-Lie subalgebra.

\smallskip
{\bf (b)} The subsets $\gry^0_{\bar g}$  and $\grw_{\bar g}$ of $\Lie(P_k)$ are Lie subalgebras. 

\smallskip
{\bf (c)} The natural map $\tilde\grz_{\bar g}\to\grz_{\bar g}$ is a bijection.

\smallskip
{\bf (d)} The set $\grz_{\bar g}$ is the image in $\Lie(P_k)$ of the set 
$\{x\in\Lie(P)|\exists y\in\Lie(N)\,\text{with}\,g\phi(x+py)-x\in p\Lie(G)\}$.}
\medskip

{\bf Proof:} 
As $L_g$ and $I_G$ are Lie homomorphisms, part (a) follows. Part (b) follows from (a). To prove (c) it suffices to show that there exists no $y\in \tilde\grz_{\bar g}\setminus\{0\}$ such that $\bar I_G(y)=0$. Each such $y$ is a non-zero element of the image of $p\Lie(N)$ in $(\Lie(P)+p\Lie(G))/p(\Lie(P)+p\Lie(G))$. Thus $\bar L_g(y)=\bar g(\overline{\phi(y)})\in\Lie(G_k)\setminus\{0\}$ and therefore  $\bar I_G(y)=0\neq\bar L_g(y)$. Thus $y\notin \tilde\grz_{\bar g}$ and therefore such an element $y$ does not exist. Thus (c) holds. 

Part (d) is only a translation of the definition of $\grz_{\bar g}$ in terms of elements of $\Lie(G)$.
\endproof

\smallskip
\noindent
{\bf Example 7.4.} Suppose $(r,d)=(2,1)$, $G=\pmb{GL}_M$, $\phi(e_1)=pe_1$, and $\phi(e_2)=e_2$. Then $\{\bar e_{1,1},\bar e_{2,2}\}$ is an $\dbF_p$-basis for $\grz_{1_{\bar M}}$. Let $w\in W_{\pmb{GL}_M}$ be the non-identity element. We can assume $g_w$ permutes $e_1$ and $e_2$. Then $g_w\phi$ permutes $e_{1,1}$ and $e_{2,2}$ and maps $e_{1,2}$ to $pe_{2,1}$. If $\gamma\in W(k)$, then $g_w\phi(p\sigma^{-1}(\gamma)e_{2,1}+\gamma e_{1,2})-\gamma e_{1,2}\in p\Lie(G)$. Thus $k\bar e_{1,2}\subseteq\grz_{\bar g_w}$. Also $\grx_{\bar g_w}:=\{\bar\gamma\bar e_{1,1}+\sigma(\bar\gamma)\bar e_{2,2}|\bar\gamma\in k,\,\bar\gamma^{p^2}=\bar\gamma\}\subseteq\grz_{\bar g_w}$. It is easy to see that $\grz_{\bar g_w}=\grx_{\bar g_w}\oplus k\bar e_{1,2}$. Thus $\grw_{\bar g_w}=\Lie(P_k)$ and $\gry^0_{\bar g_w}=\Lie(U_k)=k\bar e_{1,2}$.

\smallskip
\noindent
{\bf Theorem 7.5}  {\it Let $w\in W_G$ be fixed. We use the notations of Subsubsection 4.1.1. We consider the direct sum $\gry_{\bar g_w}:=\oplus_{o\in O(w)^+} (\grz_{\bar g_w}\cap \bar\grg_o)$.
Then the following eight properties hold:

\medskip
{\bf (a)} We have $\grz_{\bar g_w}=\Lie(T^{g_w}_{\dbF_p})\oplus\bigoplus_{o\in O(w)} \grz_{\bar g_w}\cap\bar\grg_o$.

\smallskip
{\bf (b)} Suppose $o\in O(w)^0$. Then $\grz_{\bar g_w}\cap\bar\grg_o=\bar\grg_{o\dbZ_p}$, where $\grg_{o\dbZ_p}:=\{x\in\grg_o|g_w\phi(x)=x\}$ is a $\dbZ_p$-structure of $\grg_o$ and where $\bar\grg_{o\dbZ_p}:=\grg_{o\dbZ_p}/p\grg_{o\dbZ_p}$. Thus $\grz_{\bar g_w}\cap\bar\grg_o$ is an $\dbF_p$-vector space of dimension $\dim_k(\bar\grg_o)$. 

\smallskip
{\bf (c)} If $o\notin (O(w)^+\cup O(w)^0)$, then $\grz_{\bar g_w}\cap\bar\grg_o=0$. Also the restriction to $\gry_{\bar g_w}$ of the projection of $\Lie(G_k)$ on $\oplus_{\alpha\in\Phi_N^{+w}} \bar\grg_{\pi_w^{w_{\alpha}}(\alpha)}$ along $\Lie(T_k)\oplus\bigoplus_{\alpha\in\Phi\setminus\{\pi_w^{w_{\alpha}}(\alpha)|\alpha\in\Phi_N^{+w}\}} \bar\grg_{\alpha}$ is a bijection and thus $\gry_{\bar g_w}$ has a natural structure of a $k$-vector space of dimension $\dbS(w)$ but which in general does not define a $k$-vector subspace of $\Lie(G_k)$. 

\smallskip
{\bf (d)} If $o\in O(w)^+$, then $\grw_{\bar g_w}\cap\bar\grg_o=\oplus_{\beta\in \tilde o}\bar\grg_{\beta}$. Thus $\dim_k(\grw_{\bar g_w}\cap\bar\grg_o)=\sum_{\alpha\in o^+} w_{\alpha}$.

\smallskip
{\bf (e)} We have $\grz_{\bar g_w}=\grx_{\bar g_w}^{\dbF_p}\oplus\gry_{\bar g_w}$, where $\grx_{\bar g_w}^{\dbF_p}:=\Lie(T^{g_w}_{\dbF_p})\oplus\bigoplus_{o\in O(w)^0} \bar\grg_{o\dbZ_p}$. Therefore we have identities $\grw_{\bar g_w}=(\grx_{\bar g_w}^{\dbF_p}\otimes_{\dbF_p} k)\oplus\bigoplus_{o\in O(w)^+,\;\beta\in \tilde o} \bar\grg_{\beta}=\Lie(T_k)\oplus\bigoplus_{o\in O(w)^0} \bar\grg_o\oplus\bigoplus_{o\in O(w)^+,\;\beta\in \tilde o} \bar\grg_{\beta}$.

\smallskip
{\bf (f)} We have $\gry^0_{\bar g_w}=\oplus_{\alpha\in\Phi_N^{+w},\,w_{\alpha}=1} \bar\grg_{\pi_w(\alpha)}=\Lie(U_k)\cap\grz_{\bar g_w}$.

\smallskip
{\bf (g)} We have $\grx^{\dbF_p}_{\bar g_w}=\Lie(L_k)\cap\grz_{\bar g_w}$.

\smallskip
{\bf (h)} We have $\grw_{\bar g_w}\cap\Lie(U_k)=\oplus_{\alpha\in\Phi_N^{+w}} \bar\grg_{\pi_w^{w_{\alpha}}(\alpha)}$.}
\medskip

{\bf Proof:}
For $z\in\Lie(G)$ we write 
$z=z_t+\sum_{o\in O(w)} z_o=z_t+\sum_{\alpha\in\Phi} z_{\alpha}$, 
where $z_t\in\Lie(T)$, $z_o\in\grg_o$, and $z_{\alpha}\in\grg_{\alpha}$. Let $\bar x\in\grz_{\bar g_w}$. Let $x\in\Lie(P)$ lifting $\bar x$ and $y\in\Lie(N)$ be such that $g_w\phi(x+py)-x\in p\Lie(G)$, cf. Lemma 7.3 (d). As the decomposition $\Lie(G_{B(k)})=\Lie(T_{B(k)})\oplus\bigoplus_{o\in O(w)} \grg_o[{1\over p}]$ is normalized by $g_w\phi$, we get 
$$g_w\phi(x_t+py_t)-x_t\in p\Lie(T)\;\;\text{and}\;\;g_w\phi(x_o+py_o)-x_o\in p\grg_o\;\;\forall o\in O(w).$$ 
From Lemma 7.3 (d) and the fact that $y_t=0$, we get that $\bar x_t\in\Lie(T^{g_w}_{\dbF_p})$ and $\bar x_o\in \grz_{\bar g_w}\cap\bar\grg_o$. Thus $\grz_{\bar g_w}\subseteq \Lie(T^{g_w}_{\dbF_p})\oplus\bigoplus_{o\in O(w)} \grz_{\bar g_w}\cap\bar\grg_o$. As $\Lie(T^{g_w}_{\dbZ_p})$ is fixed by $g_w\phi$, $\grz_{\bar g_w}$ contains $\Lie(T^{g_w}_{\dbF_p})$. Therefore (a) holds. 

If $o\in O(w)^0$, then $y_o=0$ and thus the relation $g_w\phi(x_o)-x_o\in p\grg_o$ is equivalent to the relation $\bar x_o\in\grg_{o\dbZ_p}/p\grg_{o\dbZ_p}$. Therefore (b) holds. 

Let $o\notin O(w)^0$. We list $o=\{\alpha^1,\ldots,\alpha^{|o|}\}$, where $\pi_w(\alpha^i)=\alpha^{i+1}$ with $\alpha^{|o|+j}:=\alpha^j$. If $o^+\neq\emptyset$ we assume $\alpha^{|o|}\in\Phi_N^{+w}$. As $g_w\phi(x_o+py_o)-x_o\in p\grg_o$, we have $g_w\phi(x_{\alpha^i}+py_{\alpha^i})-x_{\alpha^{i+1}}\in p\grg_{\alpha^{i+1}}$ for all $i\in\{1,\ldots,|o|\}$. But $g_w\phi(\grg_{\alpha})=p^{\eps(\alpha)}\grg_{\pi_w(\alpha)}$ (see Subsubsection 4.2.2). As $y\in\Lie(N)$ and $x\in\Lie(P)$, we have $y_{\alpha^i}=0$ if $\alpha^i\notin\Phi_N$ and $x_{\alpha^i}=0$ if $\alpha^i\in\Phi_N$. From the last three sentences we easily get that: 

\medskip
{\bf (i)} if $\alpha^i\not\in\tilde o$, then $x_{\alpha^i}\in p\grg_{\alpha^i}$;

\smallskip
{\bf (ii)} if $\alpha^i\in \tilde o$ is of the form $\pi_w^{j}(\alpha^{i_0})$ with $\alpha^{i_0}\in\Phi_N^{+w}$ and $j\in\{1,\ldots,w_{\alpha^{i_0}}\}$, then $\bar x_{\alpha^i}$ can be any element of $\bar\grg_{\alpha^i}$ and it is uniquely determined by $\bar y_{\alpha^{i_0}}$ in a $\sigma^{j}$-linear way.
\medskip

Thus (c) follows from properties (i) and (ii).

We check (d). We have $\grw_{\bar g_w}\cap\bar\grg_o\subseteq\oplus_{\beta\in \tilde o}\bar\grg_{\beta}$, cf. properties (i) and (ii). Referring to (ii), if we multiply $\bar y_{\alpha^{i_0}}$ by $\bar\gamma\in \dbG_m(k)$, then for each $j\in\{1,\ldots,w_{\alpha^{i_0}}\}$ the element $\bar x_{\alpha^{i_0+j}}$ gets multiplied by $\bar\gamma^{p^j}$. Let $\bar\gamma_1,\ldots,\bar\gamma_{w_{\alpha^{i_0}}}$ be elements of $k$ such that their $p$-th powers are linearly independent over $\dbF_p$. The Moore determinant of the $w_{\alpha^{i_0}}\times w_{\alpha^{i_0}}$ matrix whose rows indexed by $s\in\{1,\ldots,w_{\alpha^{i_0}}\}$ are $(\bar\gamma_s^p,\ldots,\bar\gamma_s^{p^{w_{\alpha^{i_0}}}})$ is invertible, cf. [10, Def. 1.3.2 and Lemma 1.3.3]. Thus if for $s\in\{1,\ldots,w_{\alpha^{i_0}}\}$ we consider an element $\bar x_s\in \grz_{\bar g_w}\cap \bar\grg_o$ such that its component in $\bar\grg_{\alpha^{i_0+1}}$ is $\bar\gamma_s^p$ times a fixed generator of $\bar\grg_{\alpha^{i_0+1}}$ and its component in $\bar\grg_{\alpha^j}$ is $0$ if $j$ mod $|o|$ does not belong to $\{i_0+1,\ldots,i_0+w_{\alpha^{i_0}}\}$ mod $|o|$, then the $k$-span of $x_1,\ldots,x_{w_{\alpha^{i_0}}}$ is on one hand $\oplus_{j=1}^{w_{\alpha^{i_0}}} \bar\grg_{\alpha^{i_0+j}}$ (cf. the mentioned invertibility) and on the other hand it is included in $\grw_{\bar g_w}$. As $\oplus_{\beta\in \tilde o} \bar\grg_{\beta}$ is a direct sum of $k$-vector spaces of the form $\oplus_{j=1}^{w_{\alpha^{i_0}}} \bar\grg_{\alpha^{i_0+j}}$, we have $\oplus_{\beta\in \tilde o}\bar\grg_{\beta}\subseteq \grw_{\bar g_w}\cap\bar\grg_o$. Thus $\grw_{\bar g_w}\cap\bar\grg_o=\oplus_{\beta\in \tilde o}\bar\grg_{\beta}$ i.e., (d) holds. 

Parts (e) to (h) are left as an easy exercise to the reader.\endproof

\smallskip
\noindent
{\bf Theorem 7.6.}  {\it Let $w\in W_G$. There exists a unique connected, smooth subgroup $X_{\bar g_w}$ of $G_k$ whose Lie algebra is $\grx_{\bar g_w}:=\grx_{\bar g_w}^{\dbF_p}\otimes_{\dbF_p} k$ and which contains a maximal torus of $G_k$. The group $X_{\bar g_w}$ has a natural $\dbF_p$-structure $X^{\dbF_p}_{\bar g_w}$ such that $\Lie(X^{\dbF_p}_{\bar g_w})=\grx_{\bar g_w}^{\dbF_p}$ and $T^{g_w}_{\dbF_p}$ is a maximal torus of $X^{\dbF_p}_{\bar g_w}$. Moreover, the group $X^{\dbF_p}_{\bar g_w}$ is reductive.}
\medskip

{\bf Proof:}
For $i\in\dbZ$ let $\grl_i:=\Lie(T_k)\oplus\bigoplus_{\alpha\in\Phi_L} \bar\grg_{\pi^i_w(\alpha)}$. Thus $\grl_i$ is the Lie algebra of the centralizer of $\sigma_{g_w}^i(\mu_k)$ in $G_k$. But $\grx_{\bar g_w}=\cap_{i\in\dbZ} \grl_i$. Therefore $\grx_{\bar g_w}$ is the Lie algebra of the centralizer $X_{\bar g_w}$ in $G_k$ of the subtorus $T^{1\bar g_w}_k$ of $T_k$ generated by the images of $\sigma_{g_w}^i(\mu_k)$'s, where $i\in\dbZ$. The group $X_{\bar g_w}$ is reductive (cf. Proposition 3.3 (a)) and we have $T_k\leqslant X_{\bar g_w}$. The torus $T^{1\bar g_w}_k$ is normalized by $\sigma_{g_w}$ and thus it is the pull back to $\Spec(k)$ of a subtorus $T_{\dbF_p}^{1\bar g_w}$ of $T^{\bar g_w}_{\dbF_p}$. Let $X_{\bar g_w}^{\dbF_p}$ be the centralizer of $T_{\dbF_p}^{1\bar g_w}$ in $G_{\dbF_p}^{\bar g_w}$; it is an $\dbF_p$-structure of $X_{\bar g_w}$ and thus it is reductive. 

We define an isogeny decomposition of $X_{\bar g_w}$. Let $A_k'$ and $A_k''$ be reductive subgroups of $G_k$ which are normal and for which the following two properties hold: (i) $A_k'$ is the subgroup of $G_k$ generated by cocharacters of $G_k$ that are cocharacters of $\pmb{GL}_{\bar M}$ of weights $\{-1,0\}$, and (ii) we have a natural central isogeny $A_k'\times_k A_k''\to G_k$. As $T^{1\bar g_w}_k\leqslant A_k'$, we have $A_k''\vartriangleleft X_{\bar g_w}$. Let $Y_{\bar g_w}:=X_{\bar g_w}\cap A_k'$. As $T_k\leqslant X_{\bar g_w}$, the inverse image $X_{\bar g_w}'$ of $X_{\bar g_w}$ through the isogeny $A_k'\times_k A_k''\to G_k$ contains a maximal torus of $A_k'\times_k A_k''$. As $\Ker(A_k'\times_k A_k''\to G_k)$ is contained in each maximal torus of $A_k'\times_k A_k''$, we get that $X_{\bar g_w}'$ is generated by a smooth, connected group and by a maximal torus of $A_k'\times_k A_k''$ and therefore it is smooth and connected. As $A_k''\vartriangleleft X_{\bar g_w}$, $X_{\bar g_w}'$ is isogenous to $X_{\bar g_w}$ and is equal to $Y_{\bar g_w}\times_k A_k''$. Thus $Y_{\bar g_w}$ is a connected, smooth group that contains a maximal torus of $A_k'$.

Let $\tilde X_{\bar g_w}$ be a connected, smooth subgroup of $G_k$ that contains a maximal torus of $G_k$ and we have $\Lie(\tilde X_{\bar g_w})=\Lie(X_{\bar g_w})$. As above, the inverse image of $\tilde X_{\bar g_w}$ in $A_k'\times_k A_k''$ is a smooth, connected group that contains a maximal torus $S_k'\times_k S_k''$ of $A_k'\times_k A_k''$; here $S_k'\leqslant A_k'$ and $S_k''\leqslant A_k''$. We check that $\tilde X_{\bar g_w}=X_{\bar g_w}$ i.e., $X_{\bar g_w}$ is unique.

Let $\tilde S_k'$ be the subtorus of $S_k'$ generated by cocharacters of $S_k'$ that are cocharacters of $\pmb{GL}_{\bar M}$ of weights $\{-1,0\}$. We claim that $S_k'=\tilde S_k'$. As $\tilde S_k'$ is a normal subgroup of the normalizer of $S_k'$ in $A_k'$ (i.e., it is normalized by the Weyl group of $A_k'$), to prove this claim it suffices to show that $\tilde S_k'$ has a non-trivial image in each simple factor of $A_k^{'\ad}$ and surjects onto $A_k^{'\ab}$. But these two properties follow from the definition of $A_k'$.

As $\Lie(S_k')\subseteq\Lie(A_k')\cap \Lie(X_{\bar g_w}^0)=\Lie(A_k')\cap \Lie(X_{\bar g_w})=\Lie(Y_{\bar g_w}),$
the torus $S_k'$ centralizes a maximal torus of $Y_{\bar g_w}$ (cf. Proposition 2.4 (b) applied with $\scrG=A_k'$ based on the previous paragraph). Thus $S_k'$ is this maximal torus of $Y_{\bar g_w}$. Therefore $\im(S_k'\times_k S_k''\to G_k)$ is a maximal torus of both $X_{\bar g_w}$ and $\tilde X_{\bar g_w}$. Thus we have $\tilde X_{\bar g_w}=X_{\bar g_w}$, cf. Proposition 3.3 (c).\endproof

\smallskip
\noindent
{\bf Remark 7.7.}
One checks that $\Lie(S_k')$ is a Cartan Lie subalgebra of $\Lie(A_k')$ in the sense of [6, Vol. II, Exp. XIII, Sect. 4]. Thus the uniqueness of $Y_{\bar g_w}$ (and thus of $X_{\bar g_w}$) also follows from [6, Vol. II, Exp. XIII, Cor. 5.3]. This approach to prove the uniqueness part of Theorem 7.6, is more suited in connection to Subsection 1.8.

\smallskip
\noindent
{\bf Theorem 7.8}  {\it 
Let $w\in W_G$. With the notations of Subsubsections 4.1.1 and 7.1 we have:

\medskip
{\bf (a)} Let $\Phi_w:=\{\alpha\in\Phi|\alpha\in O(w)^0\cup\bigcup_{o\in O(w)^+} \tilde o\}$; thus $\grw_{\bar g_w}=\Lie(T_k)\oplus\bigoplus_{\alpha\in\Phi_w} \bar\grg_{\alpha}$ (cf. Theorem 7.5 (e)). Then $\Phi_w$ is a closed subset of $\Phi$.

\smallskip
{\bf (b)} There exists a unique connected, smooth subgroup $W_{\bar g_w}$ of $G_k$ whose Lie algebra is $\grw_{\bar g_w}$ and which contains a maximal torus of $G_k$. We have $T_k\leqslant X_{\bar g_w}\leqslant W_{\bar g_w}\leqslant P_k$. 

\smallskip
{\bf (c)} The group $\text{Aut}_{\bar g_w}^{\text{red}}$ normalizes $W_{\bar g_w}$.}
\medskip

{\bf Proof:}
Let $\beta_1,\beta_2\in\Phi_w$ be such that $\beta_3:=\beta_1+\beta_2\in\Phi$. To prove (a) we need to check that $\beta_3\in\Phi_w$. For $i\in\{1,2,3\}$ let $o_i\in O(w)$ be such that $\beta_i\in o_i$. We consider three cases: (i) $o_1,o_2\in O(w)^0$, (ii) $o_1\in O(w)^+$ and $o_2\in O(w)^0$, and (iii) $o_1, o_2\in O(w)^+$. If $o_1,o_2\in O(w)^0$, then $o_3\in O(w)^0\subseteq\Phi_w$. 

If $o_1\in O(w)^+$ and $o_2\in O(w)^0$, then let $\alpha_1\in\Phi_N^{+w}$ be such that there exists $i\in\{1,\ldots,w_{\alpha_1}\}$ with the property that $\beta_1=\pi_w^i(\alpha_1)\in \tilde o_1$. As $o_2\in O(w)^0$, we have $\pi_w^{-j}(\beta_2)\in\Phi_{L}$ for all $j\in\dbZ$. Thus 
$$\alpha_3:=\pi_w^{-i}(\beta_3)=\alpha_1+\pi_w^{-i}(\beta_2)\in (\Phi_N+\Phi_{L})\cap \Phi\subseteq\Phi_N.$$ 
Moreover $\pi_w^j(\alpha_3)\in(\Phi_{L}+\Phi_{L})\cap \Phi\subseteq\Phi_{L}$ if $j\in\{1,\ldots,w_{\alpha_1}-1\}$ and $\pi_w^{w_{\alpha_1}}(\alpha_3)\in (\Phi_{L}+\Phi_U)\cap\Phi\subseteq\Phi_U$. Thus $\alpha_3\in\Phi_N^{+w}$ and $w_{\alpha_3}=w_{\alpha_1}$. As $\beta_3=\pi_w^i(\alpha_3)$ and as $1\le i\le w_{\alpha_1}=w_{\alpha_3}$, we get $\beta_3\in\tilde o_3\subseteq\Phi_w$.

Let now $o_1, o_2\in O(w)^+$. Thus $\beta_1\in\tilde o_1$ and $\beta_2\in\tilde o_2$. Let $\alpha_1, \alpha_2\in\Phi_N^{+w}$ be such that $\beta_1=\pi_w^{i_1}(\alpha_1)$ and $\beta_2=\pi_w^{i_2}(\alpha_2)$, where $i_1\in \{1,\ldots,w_{\alpha_1}\}$ and $i_2\in \{1,\ldots,w_{\alpha_2}\}$. We can assume that $i_1\le i_2$. Let $\alpha_3:=\alpha_1+\pi_w^{i_2-i_1}(\alpha_2)=\pi_w^{-i_1}(\beta_3)\in\Phi$. As in the previous case we argue that we have $\alpha_3\in\Phi_N^{+w}$ (with $w_{\alpha_3}=\min\{w_{\alpha_1},w_{\alpha_2}+i_1-i_2\}$) and that $\beta_3\in\tilde o_3\subseteq\Phi_w$. Thus (a) holds.

We prove (b). Due to (a), from the end of Subsubsection 3.3.1 we get that there exists a unique connected, smooth  subgroup $W_{\grw_{\bar g_w}}$ of $G_k$ which contains $T_k$ and whose Lie algebra is $\grw_{\bar g_w}$. As $\cup_{o\in O(w)^0} o\subseteq\Phi_w\subseteq \Phi_P$, we have $T_k\leqslant X_{\bar g_w}\leqslant W_{\grw_{\bar g_w}}\leqslant P_k$. Thus $W_{\grw_{\bar g_w}}$ contains the subgroup $A_k''$ of $X_{\bar g_w}$ of the proof of Theorem 7.6. Thus the uniqueness of $W_{\bar g_w}$ is proved similarly to the uniqueness of $X_{\bar g_w}$ of Theorem 7.6. Thus (b) holds.

As $\text{Aut}_{\bar g_w}^{\text{red}}$ normalizes $\grw_{\bar g_w}$ (cf. Lemma 7.1 (c)), from the uniqueness part of (b) we get that each $\bar h\in \text{Aut}_{\bar g_w}^{\text{red}}(k)$ normalizes $W_{\bar g_w}$. Thus $\text{Aut}_{\bar g_w}^{\text{red}}$ normalizes $W_{\bar g_w}$. Therefore (c) holds.\endproof

\bigskip
\noindent
{\boldsectionfont \S8. Two applications} 

\bigskip
In this Section we include two applications. In Subsection 8.1 we include a complement to Proposition 5.4 that pertains to rational inner isomorphisms. Proposition 8.1 and Theorem 8.3 generalize the well known fact that the Barsotti--Tate group of a finite product of supersingular elliptic curves over $k$ is uniquely determined by its truncation of level $1$. We use the notations of Subsections 4.1, 4.1.1, and 4.2.2.

\medskip
\noindent
{\bf 8.1. A complement to Proposition 5.4} 

We assume that the cocharacter $\mu:\dbG_m\to G$ does not factor through $Z(G)$; thus the sets $\Phi_N$ and $\Phi_U$ are non-empty. Let $w\in W_G$ be such that not all slopes of $(\Lie(G_{B(k)}),g_w\phi)$ are $0$. Let $o\in O(w)$ be such that $m_w^+(o)\neq m_w^-(o)$, cf. (7). Based on  this and (8), we can choose $o\in O(w)$ such that  $m_w^+(o)< m_w^-(o)$; thus all slopes of $(\grg_o,g_w\phi)$ are negative (cf. (7)). As the set $o\cap\Phi_N$ has more elements than the set $o\cap \Phi_U$, there exists $\alpha\in o$ such that $\alpha\in\Phi_N^{-w}$. Let $x_{w_{\alpha}}\in\grg_{\pi_w^{w_{\alpha}}(\alpha)}\subseteq\Lie(N)$ be such that $x_{w_{\alpha}}\notin p\grg_{\pi_w^{w_{\alpha}}(\alpha)}$. We have $<\scrC_{(1_{\bar M}+\bar x_{w_{\alpha}})\bar g_w}>=<\scrC_{\bar g_w}>$, cf. proof of Proposition 5.4. 

Next we check that $\scrC_{(1_M+x_{w_{\alpha}})g_w}$ and $\scrC_{g_w}$ are rational inner isomorphic i.e.,  there exists $h\in G(B(k))$ such that $h(1_M+x_{w_{\alpha}})g_w\phi=g_w\phi h$ (equivalently such that $h(1_M+x_{w_{\alpha}})[(g_w\phi)(h^{-1})]=1_M$). 

For $i\in\dbN\cup\{0\}$ let $x_{w_{\alpha}}(i):=(g_w\phi)^{-i}(x_{w_{\alpha}})$; thus $x_{w_{\alpha}}(0)=x_{w_{\alpha}}$. As all slopes of $(\grg_o,(g_w\phi)^{-1})$ are positive, the elements $x_{w_{\alpha}}(i)\in\Lie(G_{B(k)})$ converge to $0$ in the $p$-adic topology. As $x_{w_{\alpha}}(i)^2=0$ for all $i\in\dbN$, $h^{-1}:=\text{lim}_{j\to\infty}\prod_{i=1}^{j} (1_M-x_{w_{\alpha}}(i))\in \pmb{GL}_M(B(k))$
is well defined and belongs to $G(B(k))$. We have
$$h(1_M+x_{w_{\alpha}})[(g_w\phi)(h^{-1})]=[\operatornamewithlimits{\lim}\limits_{j\to\infty}\prod_{i=j}^1 (1_M+x_{w_{\alpha}}(i))](1_M+x_{w_{\alpha}})\prod_{i=0}^{j-1} (1_M-x_{w_{\alpha}}(i))=$$
$$\operatornamewithlimits{\lim}\limits_{j\to\infty} 1_M+x_{w_{\alpha}}(j)=1_M.$$
\smallskip
\noindent
{\bf Proposition 8.1.} {\it  Let $w\in W_G$. Then the following three statements are equivalent:

\medskip
{\bf (i)}  we have  $\dim(o_{\bar g_w})=\dim(P_k)$;

\smallskip
{\bf (ii)} we have $\dbS(w)=\dim(N_k)$;

\smallskip
{\bf (iii)} if $\grd:=\{x\in\Lie(G)|g_w(\phi(x))=x\}$, then we have 
$$p\Lie(G)\subseteq \grd\otimes_{\dbZ_p} W(k)\subseteq \Lie(P)+p\Lie(G).\leqno (19)$$
\indent
Moreover if statements (i) to (iii) hold, then all slopes of $(\Lie(G_{B(k)}),g_w\phi)$ are $0$.}
\medskip

{\bf Proof:} 
As $\dim(G_k)=\dim(P_k)+\dim(N_k)$, from (12) we get that statements (i) and (ii) are equivalent. We check that the statement (ii) implies the statement (iii). As $\dbS(w)=\dim(N_k)$, we have $\Phi_N=\Phi_N^{+w}$. As $\Phi_N=\Phi_N^{+w}$, for each orbit $o\in O(w)\setminus O(w)^0$, the epsilon string $\eps_w(o)$ is special in the sense that after disregarding the $0$'s, each $-1$ is followed by $+1$ (in a circular reading). Thus $m_w^+(o)\ge m_w^-(o)$ for all $o\in O(w)\setminus O(w)^0$. If $o\in O(w)^0$, then  we have $m_w^+(o)=0\ge 0=m_w^-(o)$ (cf. Subsubsection 4.2.2). Thus we have $m_w^+(o)\ge m_w^-(o)$ for all $o\in O(w)$. From this and (8) we get that we have 
$$m_w^+(o)=m_w^-(o)\;\;\forall o\in O(w).$$ 
Thus all slopes of $(\grg_o,g_w\phi)$ are $0$, cf. (7). As all slopes of $(\Lie(T),g_w\phi)$ are also $0$, we get that all slopes of $(\Lie(G_{B(k)}),g_w\phi)$ are $0$. Thus $\grd\otimes_{\dbZ_p} B(k)=\Lie(G_{B(k)})$. 

If $o\in O(w)\setminus O(w)^0$, then $m_w^+(o)=m_w^-(o)>0$. From this and the identity $\Phi_N=\Phi_N^{+w}$ we get that $\tilde o\neq\emptyset$ and $(o\setminus\tilde o)\neq\emptyset$. In particular, we have $O(w)^+=O(w)\setminus O(w)^0$. 

As $\Lie(G)=\Lie(T)\oplus\bigoplus_{o\in O(w)} \grg_o$ and $\Lie(T_{\dbZ_p}^{g_w})\subseteq \grd\subseteq \Lie(G)\cap (g_w\phi)^{-1}(\Lie(G))=\Lie(P)+p\Lie(G)$, 
to prove (19) it suffices to show for all $o\in O(w)$ we have $p\grg_o\subseteq \grd\otimes_{\dbZ_p} W(k)$. If $o\in O(w)^0$, then $\grg_{o\dbZ_p}\subseteq \grd$ (see Theorem 7.5 (b) for $\grg_{o\dbZ_p}$) and thus $\grg_o=\grg_{o\dbZ_p}\otimes_{\dbZ_p} W(k)\subseteq \grd\otimes_{\dbZ_p} W(k)$. For $o\in O(w)^+=O(w)\setminus O(w)^0$ let 
$$\grg_o^-:=\bigoplus_{\alpha\in\tilde o} \grg_{\alpha}\oplus\bigoplus_{\alpha\in o\setminus\tilde o} p\grg_{\alpha}.$$ 
\noindent
As $\Phi_N^{+w}=\Phi_N$ and $m_w^+(o)=m_w^-(o)$, the epsilon string $\eps_w(o)$ is very special in the sense that after disregarding the $0$'s, the $-1$ and $+1$ are alternating each others (in a circular reading). For $\alpha\in o$ we have $g_w\phi(\grg_{\alpha})=p^{\eps(\alpha)}\grg_{\pi_w(\alpha)}$, cf. Subsubsection 4.2.2. From the last two sentences we get the identities of the below four cases:

\medskip
-- if $\alpha,\pi_w(\alpha)\in \tilde o$, then $\alpha\in\Phi_L$ and thus $g_w\phi(\grg_{\alpha})=\grg_{\pi_w(\alpha)}$;

\smallskip
-- if $(\alpha,\pi_w(\alpha))\in \tilde o\times (o\setminus\tilde o)$, then $\alpha\in\Phi_U$ and thus $g_w\phi(\grg_{\alpha})=p\grg_{\pi_w(\alpha)}$;

\smallskip
-- if $\alpha,\pi_w(\alpha)\in o\setminus\tilde o$, then $\alpha\in\Phi_L$ and thus $g_w\phi(p\grg_{\alpha})=p\grg_{\pi_w(\alpha)}$;

\smallskip
-- if $(\alpha,\pi_w(\alpha))\in (o\setminus\tilde o)\times\tilde o$, then $\alpha\in\Phi_N$ and thus $g_w\phi(p\grg_{\alpha})=\grg_{\pi_w(\alpha)}$.

\medskip\noindent
These cases imply $g_w\phi(\grg_o^-)=\grg_o^-$. Thus $\grg_o^-$ is $W(k)$-generated by elements of $\grd$; therefore $\grg_o^-\subseteq\grd\otimes_{\dbZ_p} W(k)$. As $p\grg_o\subseteq\grg_o^-$ we get $p\grg_o\subseteq \grd\otimes_{\dbZ_p} W(k)$. Thus (19) holds. Thus the statement (ii) implies the statement (iii). 

As for $o\in O(w)^+$ we have $\tilde o\neq\emptyset$ and $o\setminus \tilde o\neq\emptyset$, it is easy to see that $\grg_o^-$ is a direct summand of $\grd\otimes_{\dbZ_p} W(k)$. Thus there exists a function $\eta:\Phi\to\{0,1\}$ such that we have an identity
$$\grd\otimes_{\dbZ_p} W(k)=\Lie(T)\oplus\bigoplus_{\alpha\in\Phi} p^{\eta(\alpha)}\grg_{\alpha}.\leqno (20)$$
We have $\eta(\alpha)=0$ if and only if $\alpha\in\Phi_w=\bigcup_{o\in O(w)^0} o\cup\bigcup_{o\in O(w)^+} \tilde o$. 

We check that the statement (iii) implies the statement (ii); thus (19) holds. As $g_w\phi$ normalizes $T(W(k))$ and $\grd\otimes_{\dbZ_p} W(k)$, it also normalizes the $W(k)$-submodule $\sum_{t\in T(W(k))} t(\grd\otimes_{\dbZ_p} W(k))$ of $\Lie(G)$. Thus this $W(k)$-submodule is $W(k)$-generated by $\grd$ and therefore $T(W(k))$ normalizes $\grd\otimes_{\dbZ_p} W(k)$. From this and (19) we get that there exists a function $\eta:\Phi\to\{0,1\}$ such that Formula (20) holds. As $\im(\eta)\subseteq\{0,1\}$, for $o\in O(w)$ the epsilon string  $\eps_w(o)$ is special. Therefore $\Phi_N^{+w}=\Phi_N$. Thus $\dbS(w)=\dim(N_k)$ i.e., the statement (ii) holds.
\endproof

\smallskip
\noindent
{\bf Definition 8.2.}  Let $w\in W_G$. We say $\scrC_{g_w}$ is {\it pivotal} if statements (i) to (iii) of Proposition 8.1 hold.

\smallskip
\noindent
{\bf Theorem 8.3.}  {\it  Suppose $\scrC_{g_w}$ is pivotal. Let $g\in G(W(k))$. If $<\scrC_{\bar g}>=<\scrC_{\bar g_w}>$, then $<\scrC_g>=<\scrC_{g_w}>$.}
\medskip

{\bf Proof:}
Let the function $\eta:\Phi\to\{0,1\}$ be as in the proof of Proposition 8.1. For $\alpha\in\Phi$ we have $\eta(\alpha)=0$ if and only if $\alpha\in\Phi_w$. This implies that (see Theorem 7.8 (a) and (20))
$$\im(\grd\otimes_{\dbZ_p} k\to \Lie(G_k))=\grw_{\bar g_w}.\leqno (21)$$ 
For the sake of clarity, we divide the proof into numbered and bold-faced paragraphs. 

\smallskip
{\bf (i) A dilatation.} Let $W_{\bar g_w}\leqslant P_k$ be as in Theorem 7.8 (b). Let $D$ be the group scheme over $\Spec(W(k))$ that is the dilatation of $W_{\bar g_w}$ centered on $G_k$. If $G=\Spec(O_G)$ and $I_{W_{\bar g_w}}$ is the ideal of $O_G$ that defines $W_{\bar g_w}$, then $D=\Spec(O_D)$ is defined by the $O_G$-algebra $O_D$ generated by all ${x\over p}$ with $x\in I_{W_{\bar g_w}}$. It is known that $D$ is smooth (see [3, Ch. 3, Sect. 3.2, Prop. 3]) and that we have a natural homomorphism $D\to G$ defined by the $W(k)$-monomorphism $O_G\hookrightarrow O_D$ (see [3, Ch. 3, Sect. 3.2, Prop. 2 (d)]). Obviously $D_{B(k)}=G_{B(k)}$ and
$$D(W(k))=\{h\in G(W(k))|\bar h\in W_{\bar g_w}(k)\}\leqslant\dbL.$$ 
\indent
Let $h_x:O_G\to W(k)[\eps]/(\eps^2)$ be a $W(k)$-homomorphism that defines an element $x\in\Lie(G)$. We have $x\in\Lie(D)$ (i.e., $h_x$ factors through the monomorphism $O_G\hookrightarrow O_D$) if and only if $h_x(I_{W_{\bar g_w}})\subseteq pW(k)[\eps]/(\eps^2)$ and thus if and only if $\bar x\in\Lie(W_{\bar g_w})=\grw_{\bar g_w}$. Therefore $\Lie(D)$ is the inverse image of $\grw_{\bar g_w}\subseteq\Lie(G_k)$ via the reduction mod $p$ epimorphism $\Lie(G)\twoheadrightarrow \Lie(G_k)$. Thus we have $\Lie(D)=\grd\otimes_{\dbZ_p} W(k)$, cf. (19) and (21). 

\smallskip
{\bf (ii) The $g_w\phi$-action on $W(k)$-valued points.} Based on (20), for $\alpha\in\Phi$ we define 
$$\scrD_{\alpha}:=p^{\eta(\alpha)}\dbG_{a,\alpha}(W(k))=1_M+p^{\eta(\alpha)}\grg_{\alpha}\leqslant \dbG_{a,\alpha}(W(k))\leqslant G(W(k)).$$ 
Let $\scrD$ be the subgroup of $G(W(k))$ generated by $T(W(k))$ and $\scrD_{\alpha}$'s with $\alpha\in\Phi$. As $\eta(\alpha)\in\{0,1\}$ and as $\Lie(G)=\Lie(T)\oplus\bigoplus_{\alpha\in\Phi} \grg_{\alpha}$, we have $\Ker(G(W(k))\to G(k))\leqslant\scrD$. The group $W_{\bar g_w}$ is the subgroup of $G_k$ generated by $T_k$ and by ${\dbG_{a,\alpha}}_k$'s with $\alpha\in\Phi_w$, cf. Subsubsection 3.3.1. From the last two sentences we get that $\scrD=\{h\in G(W(k))|\bar h\in W_{\bar g_w}(k)\}=D(W(k))$. As $g_w\phi$-normalizes $\Lie(D)=\grd\otimes_{\dbZ_p} W(k)$ and $\Lie(T)$, we have $T(W(k))=\{g_w\phi(h)g_w^{-1}|h\in T(W(k))\}$ and $\scrD_{\pi_w(\alpha)}=1_M+p^{\eta(\pi_w(\alpha))}\grg_{\pi_w(\alpha)}=1_M+g_w\phi(p^{\eta(\alpha)}\grg_{\alpha})=\{g_w\phi(h)g_w^{-1}|h\in\scrD_{\alpha}\}$. Thus
$D(W(k))=\scrD=\{g_w\phi(h)g_w^{-1}|h\in\scrD\}=\{g_w\phi(h)g_w^{-1}|h\in D(W(k))\}$.

\smallskip
{\bf (iii) The $\sigma$-automorphism.} We first recall with full details the well known fact that the smooth, affine group scheme $D$ over $\Spec(W(k))$ is uniquely determined up to unique isomorphism by $D(W(k))$ and $D_{B(k)}=G_{B(k)}$. Let $\tilde D=\Spec(O_{\tilde D})$ be a smooth, affine group scheme over $\Spec(W(k))$ such that $\tilde D_{B(k)}=D_{B(k)}$ and $\tilde D(W(k))=D(W(k))$. Let $D_0=\Spec(O_0)$ be the Zariski closure of the diagonal subgroup $D_{B(k)}$ of $D_{B(k)}\times_{B(k)} \tilde D_{B(k)}$ in $D\times_{\Spec(W(k))} \tilde D$. We check that the two natural projection homomorphisms $q_0:D_0\to D$ and $\tilde q_0:D_0\to\tilde D$ are isomorphisms; this will imply the existence of a unique isomorphism $\tilde D\arrowsim D$ that extends the identity $\tilde D_{B(k)}=D_{B(k)}$. We will only check that $\tilde q_0$ is an isomorphism (as the arguments for $q_0$ are the same). As $D(W(k))=\tilde D(W(k))$, we have an isomorphism $\tilde q_0(W(k)):D_0(W(k))\arrowsim \tilde D(W(k))$. Let $O_{\tilde D}\hookrightarrow O_0$ be the $W(k)$-monomorphism that defines $\tilde q_0$. We identify $O_0[{1\over p}]=O_{\tilde D}[{1\over p}]$ and thus we can view $O_{\tilde D}$ as a $W(k)$-subalgebra of $O_0$. We show that the assumption that $O_{\tilde D}\neq O_0$ leads to a contradiction. This assumption implies that there exists $y\in O_0\setminus O_{\tilde D}$ such that $py\in O_{\tilde D}\setminus pO_{\tilde D}$. As $\tilde D$ is smooth, there exists a $W(k)$-homomorphism $h_y:O_{\tilde D}\to W(k)$ that takes $py$ to an invertible element $i_y$ of $W(k)$. As ${{i_y}\over p}\notin W(k)$, such a homomorphism $h_y$ does not factor through a  $W(k)$-homomorphism $O_0\to W(k)$. This contradicts the fact that $\tilde q_0(W(k))$ is an isomorphism. Thus $O_0=O_{\tilde D}$ and therefore $\tilde q_0$ is an isomorphism.

As $D(W(k))=\{g_w\phi(h)g_w^{-1}|h\in D(W(k))\}$, the uniqueness part of the previous paragraph implies that $g_w\phi$ induces naturally a $\sigma$-automorphism $\sigma_D$ of $D$. In other words, we have a unique isomorphism $\sigma_D:D\times_{\Spec(W(k))} {}_{\sigma} \Spec(W(k))\arrowsim D$ with the property that for $h\in D(W(k))$ we have $\sigma_D(h\times\sigma^{-1})=g_w\phi(h)g_w^{-1}\in D(W(k))$.

\smallskip
{\bf (iv) The $\dbZ_p$-structure.} We now show that $\sigma_D$ defines naturally a $\dbZ_p$-structure $D_{\dbZ_p}$ of $D$ with the property that $\Lie(D_{\dbZ_p})=\grd$. The adjoint group $D^{\ad}_{\dbQ_p}$ is the identity component of the group of Lie automorphisms of the semisimple Lie algebra $[\grd[{1\over p}],\grd[{1\over p}]]$ over $\dbQ_p$. Accordingly we will take $D_{\dbQ_p}$ to be a natural isogeny cover of $D^{\ad}_{\dbQ_p}\times_{\dbQ_p} G_{\dbQ_p}^{\ab}$  (as $g_w$ acts trivially on $Z^0(G)$, we have $G_{\dbQ_p}^{\ab}={G^{g_w}_{\dbQ_p}}^{\ab}$). This isogeny cover is the one defined naturally by the isogeny $T^{g_w}_{\dbQ_p}\to T^{0g_w}_{\dbQ_p}\times_{\dbQ_p} G^{\ab}_{\dbQ_p}$, where $T^{0g_w}_{\dbQ_p}:=\im(T^{g_w}_{\dbQ_p}\to D^{\ad}_{\dbQ_p})$ (obviously $T^{g_w}_{\dbQ_p}$ acts via inner conjugation on $\grd[{1\over p}]$ and thus also on $[\grd[{1\over p}],\grd[{1\over p}]]$). Using Galois descent, it suffices to show that there exists a finite field $k_0=\dbF_{p^q}$ such that the iterate isomorphism $D\times_{\Spec(W(k))} {}_{\sigma^q} \Spec(W(k))\arrowsim D$ induced by $\sigma_D$, defines naturally a $W(k_0)$-structure $D_{W(k_0)}$ of $D$ with the property that $\Lie(D_{W(k_0)})=\grd\otimes_{\dbZ_p} W(k_0)$. We pick $k_0$ such that the torus $T^{g_w}_{W(k_0)}$ is split. This implies that the group $D_{B(k_0)}$ is also split and that there exists a Lie subalgebra $\tilde\grd_{W(k_0)}$ of the Lie algebra $\Lie(D_{B(k_0)})$ (viewed over $W(k_0)$) such that $\Lie(G)=\tilde\grd_{W(k_0)}\otimes_{W(k_0)} W(k)$. Let $\tilde D_{W(k_0)}'$ be a reductive group scheme over $\Spec(W(k_0))$ that has $T^{g_w}_{W(k_0)}$ as a maximal torus and has $D_{B(k_0)}$ as its generic fibre, cf. [34]. We check that there exists an element $\tilde h\in T^{0g_w}_{B(k_0)}(B(k_0))$ such that $\tilde D_{W(k_0)}:=\tilde h\tilde D_{W(k_0)}'\tilde h^{-1}$ has $\tilde\grd_{W(k_0)}$ as its Lie algebra. Based on the $G^{\ad}(B(k))$-conjugacy of reductive group schemes over $\Spec(W(k))$ that have $G_{B(k)}$ as their generic fibres (see [34, end of Subsect. 2.5]) and on the $G(W(k))$-conjugacy of maximal tori of $G$ (see Subsubsection 3.3.2), such an element exists over $B(k)$. As we have $T^{0g_w}_{B(k)}(B(k))/T^{0g_w}_{W(k)}(W(k))=T^{0g_w}_{B(k_0)}(B(k_0))/T^{0g_w}_{W(k_0)}(W(k_0))$, we easily get that $\tilde h$ exists. 

Both groups $G(W(k))$ and $\tilde D_{W(k_0)}(W(k))$ are maximal bounded subgroups of $G(B(k))$ (cf. [34, Subsect. 3.2]) and thus are the maximal bounded subgroup of $G(B(k))$ that normalizes $\Lie(G)=\Lie(\tilde D_{W(k_0)})\otimes_{W(k_0)} W(k)$. As $G(W(k))=\tilde D_{W(k_0)}(W(k))$, as in (iii) we argue that we can identify $G=\tilde D_{W(k)}$. It is easy to see that $W_{\bar g_w}$ is the pull back of a closed subgroup $W_{\bar g_w}^{k_0}$ of $\tilde D_{k_0}$ and that $D$ is the pull back of the dilatation $D_{W(k_0)}$ of $W_{\bar g_w}^{k_0}$ centered on $\tilde D_{W(k)}$. Thus $D_{W(k_0)}$ is the desired $W(k_0)$-structure of $D$. Therefore $D_{\dbZ_p}$ exists. 

As $<\scrC_{\bar g}>=<\scrC_{\bar g_w}>$,  we can assume that $g=g_0g_w$ with $g_0\in\Ker(G(W(k))\to G(k))\leqslant D(W(k))$ (cf. Lemma 4.5 (c)). Let $h\in D(W(k))\leqslant\dbL$ be such that $g_0=h\sigma_D(h\times\sigma^{-1})^{-1}=hg_w\phi(h)^{-1}g_w$, cf. Lemma 3.1. Thus $hg_w\phi=g_0g_w\phi h=g\phi h$ i.e., $h$ is an inner isomorphism between $\scrC_{g_w}$ and $\scrC_g$. Thus $<\scrC_g>=<\scrC_{g_w}>$.
\endproof

\bigskip
\noindent
{\boldsectionfont \S9. The proof of Basic Theorem B}

\bigskip
Theorem 9.1 proves a technical result on $P(k)$-conjugates of $\mu_k:\dbG_m\to G_k$. We use Theorem 9.1 to prove Basic Theorem B in Subsection 9.1. Lemma 9.4 and Remark 9.5 present a result on the relation $\scrR$ on $W_G$ (see Basic Theorem B) and two remarks. We use the notations of Subsections 4.1, 4.1.1, 5.1.1, and 7.1.

\smallskip
\noindent
{\bf Theorem 9.1.}  {\it 
Let $w\in W_G$. Let $\mu_{1k}:\dbG_m\to G_k$ be a $P(k)$-conjugate of $\mu_k$. Let $\bar l_1\in\Lie(G_k)$ be the image under $d\mu_{1k}$ of the standard generator of $\Lie(\dbG_m)$. Then the following three statements are equivalent:

\medskip
{\bf (i)} the element $\bar l_1$ is contained in the $k$-span $\grw_{\bar g_w}$ of $\grz_{\bar g_w}$;

\smallskip
{\bf (ii)} the two cocharacters $\mu_{k}$ and $\mu_{1k}$ (equivalently $\bar l_0$ and $\bar l_1$) are $\text{Aut}_{\bar g_w}(k)=\text{Aut}^{\text{red}}_{\bar g_w}(k)$-conjugate;

\smallskip
{\bf (iii)} the two cocharacters $\mu_{k}$ and $\mu_{1k}$ (equivalently $\bar l_0$ and $\bar l_1$) are $\text{Aut}^{0\text{red}}_{\bar g_w}(k)$-conjugate.}
\medskip

{\bf Proof:}
Obviously the statement (iii) implies the statement (ii). We check that the statement (ii) implies the statement (i). We have $\bar l_0\in \Lie(T_k)\subseteq\grw_{\bar g_w}$ (cf. Theorem 7.5 (e)) and $\text{Aut}_{\bar g_w}(k)$ normalizes $\grw_{\bar g_w}$ (cf. Lemma 7.1 (c)). Thus, if the statement (ii) holds, then we have $\bar l_1\in \grw_{\bar g_w}$ and therefore the statement (i) holds. 

We prove that the statement (i) implies the statement (iii); thus $\bar l_1\in\grw_{\bar g_w}$. Let 
$$U^{+w}:=\prod_{\alpha\in\Phi_N^{+w}} \dbG_{a,{\pi_w^{w_{\alpha}}(\alpha)}}\leqslant U.$$ 
\noindent
We have $\Lie(U^{+w}_k)=\oplus_{\alpha\in\Phi_N^{+w}} \bar\grg_{\pi_w^{w_{\alpha}}(\alpha)}=\grw_{\bar g_w}\cap\Lie(U_k)$, cf. Theorem 7.5 (h). The group $\text{Aut}^{\text{red}}_{\bar g_w}$ normalizes the $k$-vector spaces $\grw_{\bar g_w}$ (cf. Lemma 7.1 (c)) and $\Lie(U_k)$ (as $\text{Aut}_{\bar g_w}\leqslant P_k$). Thus $\text{Aut}_{\bar g_w}^{\text{red}}$ normalizes $\Lie(U^{+w}_k)$. As $U^{+w}_k(k)=\{1_{\bar M}+\bar u|\bar u\in\Lie(U^{+w}_k)\}$ (see Subsubsection 3.3.3 for $U(k)$), we get that the group $\text{Aut}^{\text{red}}_{\bar g_w}$ normalizes $U^{+w}_k$ (another way to argue this is to use Theorem 7.8 (c) and to remark that $U^{+w}_k=W_{\bar g_w}\cap U_k$). 

For $\alpha\in\Phi_N^{+w}$, let $c_{\alpha}:{V_{a,\alpha}}_k\to H_k$ be as in the proof of Proposition 5.3. The direct sum of the tangent spaces at the identity element of the curves $\im(c_{\alpha})$'s with $\alpha\in\Phi_N^{+w}$ is $\Lie(\scrS^{0\text{red}}_{\bar g_w})$, cf. (13) and the proof of Proposition 5.3. Thus the group $\scrS^{0\text{red}}_{\bar g_w}$ is generated by the $\im(c_{\alpha})$'s with $\alpha\in\Phi_N^{+w}$. Thus $\text{Aut}^{0\text{red}}_{\bar g_w}$ is generated by the images of $a_{\bar g_w}\circ c_{\alpha}$'s with $\alpha\in\Phi_N^{+w}$, cf. Subsection 5.2 (a). Therefore the group $\text{Aut}^{0\text{red}}_{\bar g_w}(k)$ is generated by its elements which are products of the form $\bar h_1\bar h_2$ with $\bar h_1\in U^{+w}(k)$ and $\bar h_2\in L(k)$, cf. (11) and the proof of Subsection 5.2 (a). The product of two such elements $\bar h_1\bar h_2$ and $\bar h'_1\bar h'_2$ of $\text{Aut}^{0\text{red}}_{\bar g_w}(k)$ is $\bar r_1\bar r_2$, where 
$$\bar r_1:=\bar h_1\bar h_2\bar h'_1(\bar h_1\bar h_2)^{-1}\bar h_1\;\; \text{and}\;\;\bar r_2:=\bar h_2\bar h'_2\in L(k).$$ 
As $\text{Aut}^{0\text{red}}_{\bar g_w}$ normalizes $U^{+w}_k$, we have $\bar r_1\in U^{+w}(k)$. By induction, each element $\bar h\in \text{Aut}^{0\text{red}}_{\bar g_w}(k)$ is a product $\bar h=\bar h_1\bar h_2$, where $\bar h_1\in U^{+w}(k)$ and $\bar h_2\in L(k)$; as $\text{Aut}^{0\text{red}}_{\bar g_w}\leqslant P_k$, $\bar h$ is uniquely written as such a product $\bar h_1\bar h_2$. 

We check that the rule $\bar h\mapsto\bar h_1$ defines a finite, surjective $k$-morphism (in general is not a homomorphism)
$$p_{\bar g_w}:\text{Aut}^{0\text{red}}_{\bar g_w}\to U^{+w}_k.$$ 
Each non-empty $k$-fibre of $p_{\bar g_w}$ has a number of $k$-valued points equal to the order of the group 
$$B_{\bar g_w}:=(\text{Aut}^{0\text{red}}_{\bar g_w}\cap L_k)(k).$$ 
The group $B_{\bar g_w}$ is finite, cf. Lemma 5.2 (b). Thus $p_{\bar g_w}$ is quasi-finite. As $\dim(U^{+w}_k)=\dbS(w)=\dim(\text{Aut}^{0\text{red}}_{\bar g_w})$ (cf. Subsection 5.2 (d) for the last equality), the quasi-finite morphism $p_{\bar g_w}$ is dominant. Thus $\text{Aut}^{0\text{red}}_{\bar g_w}$ is an open, dense subscheme of the normalization $V^{+w}$ of $U^{+w}_k$ in the field of fractions of $\text{Aut}^{0\text{red}}_{\bar g_w}$, cf. Zariski main theorem of [29, Ch. IV]. The action of $B_{\bar g_w}$ on $\text{Aut}^{0\text{red}}_{\bar g_w}$ via right translations extends to an action on $V^{+w}$. The quotient $k$-morphism $\text{Aut}^{0\text{red}}_{\bar g_w}\to \text{Aut}^{0\text{red}}_{\bar g_w}/B_{\bar g_w}$ is a Galois cover. The natural finite factorization $V^{+w}/B_{\bar g_w}\to U^{+w}_k$ is radicial (i.e., it is universal injective on geometric points) above an open, dense subscheme of $U^{+w}_k$. Thus the field extension between the fields of fractions of $V^{+w}/B_{\bar g_w}$ and $U^{+w}_k$ is purely inseparable. Therefore there exists a bijection between prime divisors of $U^{+w}_k$ and $V^{+w}/B_{\bar g_w}$ defined by taking reduced structures on pulls back. 

We show that the assumption that the open, affine subscheme $\text{Aut}^{0\text{red}}_{\bar g_w}$ of $V^{+w}$ is not the whole affine scheme $V^{+w}$, leads to a contradiction. This assumption implies that the complement $C^{+w}$ of $\text{Aut}^{0\text{red}}_{\bar g_w}$ in $V^{+w}$ contains a divisor $I^{+w}$ (i.e., $C^{+w}$ is not of codimension $\ge 2$ in $V^{+w}$), cf. [21, Thm. 38]. The image $I^{+w}$ in $U^{+w}_k$ of $V^{+w}$ is a prime divisor of $U^{+w}_k$ whose pull back to $V^{+w}$ is contained in $C^{+w}$. As $U_k^{+w}\arrowsim\dbA_k^{\dbS(w)}$, there exists a global function $f$ of $U^{+w}_k$ whose zero locus is $I^{+w}$. We have $f\notin k$ and moreover $f$ is an invertible element of the $k$-algebra $A^0_{\bar g_w}$ of global functions of $\text{Aut}^{0\text{red}}_{\bar g_w}$. As $\text{Aut}^{0\text{red}}_{\bar g_w}$ is a connected, smooth, unipotent group (cf. Subsection 5.2 (d)), $A^0_{\bar g_w}$ is a polynomial ring over $k$ (to be compared with the case of $\scrS^{\text{red}}_{\bar g_i\dbF_{p^q}}$ in the proof of Theorem 6.1). Thus each invertible element of $A^0_{\bar g_w}$ belongs to $k$. Therefore such an element $f$ does not exist. Contradiction. 

Thus $V^{+w}=\text{Aut}^{0\text{red}}_{\bar g_w}$. Therefore the $k$-morphism $p_{\bar g_w}$ is finite and surjective.

As $P(k)=U(k)L(k)$ and $L_k$ fixes $\mu_k$, each $P(k)$-conjugate of $\mu_k$ is a $U(k)$-conjugate of $\mu_k$. Let $\bar u\in U(k)$ be such that $\mu_{1k}=\bar u\mu_k\bar u^{-1}$. We write $\bar u=1_{\bar M}+\bar y$, where $\bar y\in\Lie(U_k)$. As $\bar y\bar l_0=0$ and $\bar l_0\bar y=-\bar y$, we have 
$$\bar l_1=(1_{\bar M}+\bar y)(\bar l_0)(1_{\bar M}-\bar y)=\bar l_0+[\bar y,\bar l_0]-\bar y\bar l_0\bar y=\bar l_0+[\bar y,\bar l_0]=\bar l_0+\bar y.$$ 
As $\bar l_0\in\Lie(T_k)\subseteq\grw_{\bar g_w}$, we get that $\bar y=\bar l_1-\bar l_0\in\grw_{\bar g_w}\cap\Lie(U_k)=\Lie(U_k^{+w})$. Thus $\bar u\in U^{+w}(k)$. Let $\bar a\in \text{Aut}^{0\text{red}}_{\bar g_w}(k)$ be such that $p_{\bar g_w}(\bar a)=\bar u$. As $\bar u^{-1}\bar a\in L(k)$ centralizes $\mu_k$, the inner conjugate of $\mu_k$ through $\bar a$ is $\bar u(\mu_k)\bar u^{-1}=\mu_{1k}$. We conclude that the statement (i) implies the statement (iii).\endproof 

\smallskip
\noindent
{\bf Corollary 9.2.} {\it  Let $w_1, w_2\in W_G$. If $<\scrC_{\bar g_{w_1}}>=<\scrC_{\bar g_{w_2}}>$, then there exists $\bar h_2\in L(k)$ which is an inner isomorphism between $\scrC_{\bar g_{w_1}}$ and $\scrC_{\bar g_{w_2}}$.}
\medskip

{\bf Proof:}
Let $\bar h\in P(k)$ be an inner isomorphism between $\scrC_{\bar g_{w_1}}$ and $\scrC_{\bar g_{w_2}}$. As $\bar h(\grz_{\bar g_{w_1}})=\grz_{\bar g_{w_2}}$ (cf. Lemma 7.1 (a) and (b)), we have $\bar h(\grw_{\bar g_{w_1}})=\grw_{\bar g_{w_2}}$. Moreover $\bar l_0\in\Lie(T_k)\subseteq \grw_{\bar g_{w_s}}\cap \grw_{\bar g_{w_s}}$. Therefore $\bar h(\bar l_0),\bar l_0\in\grw_{\bar g_{w_2}}$. Let $\bar h_1\in \text{Aut}^{\text{red}}_{\bar g_{w_2}}(k)$ be such that $\bar h_1(\bar h(\bar l_0))=\bar l_0$, cf. Theorem 9.1. Let $\bar h_2:=\bar h_1\bar h$; it is an inner isomorphism between $\scrC_{\bar g_{w_1}}$ and $\scrC_{\bar g_{w_2}}$. As $\bar h_2$ centralizes $\bar l_0$, we have $\bar h_2\in L(k)$.
\endproof

\smallskip
\noindent
{\bf Corollary 9.3.} {\it  Each connected component of $\text{Aut}^{\text{red}}_{\bar g_w}$ intersects the centralizer $L_k$ of $\mu_k$ in $G_k$.}
\medskip

{\bf Proof:}
This follows from the equivalence of statements (ii) and (iii) of Theorem 9.1.
\endproof 

\medskip
\noindent
{\bf 9.1. End of the proof of Basic Theorem B} 

See Subsubsection 1.1.3 for the function $\scrS:W_G\to\scrS_{W_P}$. For $s\in\{1,2\}$ let $X_{g_{w_s}}$  be the centralizer in $G$ of the subtorus of $T$ generated by the images of the cocharacters in the set $\{\sigma_{g_{w_s}}^i(\mu)|i\in\dbZ\}$. It is a reductive group scheme (cf. Proposition 3.3 (a)) that contains $T$ and that depends only on $w_s\in W_G$. For $w\in W_G$ we have $w\in \scrS(w_s)$ if and only if $g_w\in X_{g_{w_s}}(W(k))$. Thus $\scrS(w_s)=((N_T\cap X_{g_{w_s}})/T)(W(k))\leqslant W_P$. The special fibre of $X_{g_w}$ is the subgroup $X_{\bar g_w}$ of $G_k$ we introduced in Theorem 7.6, cf. proof of Theorem 7.6.

We prove the property 1.3 (a). The only if part is trivial. We prove the if part; thus $<\scrC_{\bar g_{w_1}}>=<\scrC_{\bar g_{w_2}}>$. Let $\bar h_2\in L(k)$ be as in Corollary 9.2. We have $\bar h_2(\Lie(L_k))=\Lie(L_k)$ and (cf. Lemma 7.1 (a)) $\bar h_2(\grz_{\bar g_{w_1}})=\grz_{\bar g_{w_2}}$. Thus from Theorem 7.5 (g) we get $\bar h_2(\grx^{\dbF_p}_{\bar g_{w_1}})=\grx^{\dbF_p}_{\bar g_{w_2}}$. We compute
$$\Lie(\bar h_2X_{\bar g_{w_1}}\bar h_2^{-1})=\bar h_2(\Lie(X_{\bar g_{w_1}}))=\bar h_2(\grx^{\dbF_p}_{\bar g_{w_1}}\otimes_{\dbF_p} k)=\grx^{\dbF_p}_{\bar g_{w_2}}\otimes_{\dbF_p} k=\Lie(X_{\bar g_{w_2}}).$$ 
This implies that $\bar h_2X_{\bar g_{w_1}}\bar h_2^{-1}=X_{\bar g_{w_2}}$, cf. the uniqueness part of Theorem 7.6. 

Let $\bar h_2'\in X_{\bar g_{w_2}}(k)\leqslant L(k)$ be such that $\bar h_2'\bar h_2T_k\bar h_2^{-1}\bar h_2^{'-1}=T_k$. The element $\bar g_3:=\bar h_2'\bar h_2\in L(k)$ normalizes $T_k$. Thus $\bar g_3\in (N_T\cap L)(k)$. Let $w_3\in W_P$ be such that there exists $g_3\in (N_T\cap L)(W(k))$ which lifts $\bar g_3$ and which is a representative of $w_3^{-1}$. Based on Lemma 4.1, we can assume that $g_{w_3^{-1}w_1\sigma(w_3)}=g_{w_3^{-1}}g_{w_1}\sigma(g_{w_3})=g_{w_3^{-1}}g_{w_1}\phi(g_{w_3})$. Thus $\scrC_{g_{w_1}}$ is inner isomorphic to $\scrC_{g_{w_3^{-1}w_1\sigma(w_3)}}$. By replacing $w_1$ with $w_3^{-1}w_1\sigma(w_3)$ and thus implicitly by replacing $g_{w_1}$ with $g_{w_3^{-1}w_1\sigma(w_3)}=g_{w_3}^{-1}g_{w_1}\sigma(g_{w_3})$, we can assume $w_3$ is the identity element of $W_P$. Thus $\bar g_3\in T(k)$. This implies that $\bar h_2=\bar g_3\bar h_2^{'-1}\in X_{\bar g_{w_2}}(k)$. Thus we have $X_{\bar g_{w_1}}=\bar h_2^{-1}X_{\bar g_{w_2}}\bar h_2=X_{\bar g_{w_2}}$. 

The group scheme $X_{g_{w_s}}$ contains $T$ and its special fibre is $X_{\bar g_{w_1}}=X_{\bar g_{w_2}}$. Thus the group scheme $X_{g_{w_s}}$ does not depend on $s\in\{1,2\}$ (cf. Proposition 3.3 (d)); we will denote it simply by $X$. Let $X^{\text{n}}$ be the normalizer of $X$ in $G$, cf. Proposition 3.3 (b). As $g_{w_s}\phi$ permutes the cocharacters in the set $\{\sigma_{g_{w_s}}^i(\mu)|i\in\dbZ\}$ and as $\mu$ factors through $Z^0(X)$, we have $g_{w_s}\phi(\Lie(X))=\Lie(X)$. Thus the triples $(M,g_{w_1}\phi,X)$ and $(M,g_{w_2}\phi,X)$ are Shimura $F$-crystals over $k$ and the element $h_{12}:=g_{w_2}g_{w_1}^{-1}=(g_{w_2}\phi)(g_{w_1}\phi)^{-1}$ of $N_T(W(k))$ normalizes $\Lie(X)$. From the uniqueness part of Lemma 3.4 we get that $h_{12}$ normalizes $X$. Thus $h_{12}\in X^{\text{n}}(W(k))$. Let $h_2\in X(W(k))$ be a lift of $\bar h_2\in X(k)$. We have identities $h_2g_{w_1}\phi h_2^{-1}=g_0g_{w_2}\phi=g_0h_{12}g_{w_1}\phi$, where $g_0\in G(W(k))$ is such that $\bar g_0\in \sigma_{g_{w_2}}(N)(k)$ (cf. Lemma 4.5 (a)). As $(M,g_{w_1}\phi,X)$ is a Shimura $F$-crystal over $k$, we have $g_0h_{12}=h_2(g_{w_1}\phi)h_2^{-1}(g_{w_1}\phi)^{-1}\in X(W(k))$. Thus $\bar g_0\bar h_{12}\in X(k)$. We get $\bar g_0\in X(k)\bar h_{12}^{-1}\subseteq X^{\text{n}}(k)$. The group scheme $X=\sigma_{g_{w_2}}(X)\leqslant \sigma_{g_{w_2}}(L)$ normalizes $\sigma_{g_{w_2}}(N)$ and the intersection $X\cap \sigma_{g_{w_2}}(N)$ is the trivial group scheme over $\Spec(W(k))$. The last two sentences imply that the commutator of the element $\bar g_0\in (\sigma_{g_{w_2}}(N)_k\cap X^{\text{n}}_k)(k)$ with any element of $X_k(k)$, belongs to $\sigma_{g_{w_2}}(N)(k)\cap X_k(k)$ and therefore it is $1_{\bar M}$. Thus $\bar g_0$ centralizes $X_k(k)$ and therefore also $T_k$. Thus $\bar g_0\in T(k)$. As $\bar g_0$ is a unipotent element, we get $\bar g_0=1_{\bar M}$. Thus $\bar h_{12}=\bar g_0\bar h_{12}\in X(k)$. As $X$ is the identity component of $X^{\text{n}}$ (cf. Proposition 3.3 (b)), as $h_{12}\in X^{\text{n}}(W(k))$, and as $\bar h_{12}\in X(k)$, we have $h_{12}\in X(W(k))$. As $\mu$ factors through $Z^0(X)$, there exists an inner isomorphism between $(M,g_{w_1}\phi,X)$ and $(M,g_{w_2}\phi,X)=(M,h_{12}g_{w_1}\phi,X)$ (cf. Lemma 4.2). Thus $<\scrC_{g_{w_1}}>=<\scrC_{g_{w_2}}>$.

We prove the property 1.3 (b) The if part follows from the above proof of the property 1.3 (a) (the role of $h_{12}\in (X\cap N_T)(W(k))$ above being that of $g_{w^{-1}_4}$). To prove the only if part of the property 1.3 (b), we can assume that $w_3$ is the identity element and (cf. Lemma 4.1) that $g_{w_1}=g_{w_4}g_{w_2}$. We have $g_{w_4}\in X_{g_{w_2}}(W(k))$. The Shimura $F$-crystals $(M,g_{w_1}\phi,X_{g_{w_2}})=(M,g_{w_4}g_{w_2}\phi,X_{g_{w_2}})$ and $(M,g_{w_2}\phi,X_{g_{w_2}})$ are inner isomorphic, cf. Lemma 4.2. Thus $<\scrC_{g_{w_1}}>=<\scrC_{g_{w_2}}>$ and therefore $<\scrC_{\bar g_{w_1}}>=<\scrC_{\bar g_{w_2}}>$. Thus the property 1.3 (b) holds. Property 1.3 (c) follows from the property 1.3 (b).\endproof

\smallskip
\noindent
{\bf Lemma 9.4.} {\it  Let the relation $\scrR$ on $W_G$ be as in Subsection 1.3. Let $w\in W_G$. Then the equivalence class $[w]\in \scrR\backslash W_G$ has a number of elements equal to the order $|W_P|$ of $W_P$. Thus $\scrR\backslash W_G$ has $[W_G:W_P]$ elements.}
\medskip

{\bf Proof:}
Let $\sigma_w:W_G\arrowsim W_G$ be the composite of $\sigma:W_G\arrowsim W_G$ with the inner automorphism of $W_G$ defined by $w$; thus $\sigma_w$ is the automorphism of $W_G$ defined by the $\dbZ_p$-structure $(G_{\dbZ_p}^{g_w},T_{\dbZ_p}^{g_w})$ of $(G,T)$. We have $[w]=\{w_3w_4\sigma_w(w_3^{-1})w|w_3\in W_P,w_4\in \scrS(w)\}$. Let $w_{3,1},\ldots,w_{3,[W_P:\scrS(w)]}$ be a set of representatives for $W_P/\scrS(w)$. From very definitions, we have $\scrS(w):=\cap_{i\in\dbZ} \sigma_w^i(W_P)\leqslant W_P$; thus $\sigma_w(\scrS(w))=\scrS(w)$. Therefore 
$$[w]=\{w_{3,j}w_4\sigma_w(w_{3,j}^{-1})w|1\le j\le [W_P:\scrS(w)],\;w_4\in \scrS(w)\}.$$ 
Suppose $w_{3,j}w_4\sigma_w(w_{3,j}^{-1})=w_{3,j'}w_4'\sigma_w(w_{3,j'}^{-1})$, where $j,j'\in\{1,\ldots,[W_P:\scrS(w)]\}$ and $w_4,w_4'\in\scrS(w)$. Let $w_5:=w_{3,j'}^{-1}w_{3,j}$. We have $w_5=w_4'\sigma_w(w_5)w_4^{-1}$. As $w_4, w_4'\in \scrS(w):=\cap_{i\in\dbZ} \sigma_w^i(W_P)$, by induction on $s\ge 0$ we get $w_5\in \sigma_w^s(W_P)$. Thus $w_5\in\cap_{i\in\dbN\cup\{0\}} \sigma_w^i(W_P)=\cap_{i\in\dbZ} \sigma_w^i(W_P)=\scrS(w)$; the first equality holds as $\sigma_w\in\text{Aut}(W_G)$ has finite order. As $w_5\in\scrS(w)$, we have $j=j'$. Thus $(w_{3,j},w_4)=(w_{3,j'},w_4')$ i.e., the elements of $[w]$ are listed above without repetitions. Thus $[w]$ has $|W_P|$ elements.
\endproof

\smallskip
\noindent
{\bf Remark 9.5.} {\bf (a)} As in Example 6.3, we can use Theorem 6.1 (b) to get a short proof of Corollary 4.6 that does not rely on [17]. This is so as Basic Theorem B tells us how to get distinct orbits of $\dbT_{G_k,\sigma}$ that intersect $N_T(k)$. 

\smallskip
{\bf (b)} Corollary 6.2 also follows from Lemma 4.5 (d) and (e), the property 1.3 (b), and Lemma 9.4.

\bigskip
\noindent
{\boldsectionfont \S10. The inductive step}

\bigskip
We present an inductive approach that play a key role in the proof of the statement 1.4 (a) in Subsection 11. We use the notations of Subsections 4.1, (beginning of) 4.2, 4.2.1, and 7.1. Let $G_1$ be a reductive, closed subgroup scheme of $\pmb{GL}_M$ that contains $G$ and such that the triple $(M,\phi,G_1)$ is a Shimura $F$-crystal over $k$. Let $G_{1\dbZ_p}$ be the subgroup scheme of $\pmb{GL}_{M_{\dbZ_p}}$ whose pull back to $\Spec(W(k))$ is $G_1$, cf. Subsection 4.2 applied to $(M,\phi,G_1)$. Let $T_{1\dbZ_p}$ be a maximal torus of the centralizer of $T_{\dbZ_p}$ in $G_{1\dbZ_p}$, cf. Proposition 3.3 (a) and (e). Let $T_1:=T_{1W(k)}$; we have $\phi(\Lie(T_1))=\Lie(T_1)$. Let $N_{T_1}$ be the normalizer of $T_1$ in $G_1$ and let $W_{G_1}:=(N_{T_1}/T_1)(W(k))$. Let $P_1$, $U_1$, and $N_1$ be for $G_1$ and $\mu$ what $P$, $U$, and $N$ are for $G$ and $\mu$; thus $P_1$ is the normalizer of $F^1$ in $G_1$, etc. Let $G^{0\ad}$ be the maximal direct factor of $G^{\ad}$ in which $\mu$ has a trivial image.   

\smallskip
\noindent
{\bf Theorem 10.1.}  {\it  We assume that the group scheme $G^{\ad}$ is non-trivial and that each simple factor of either $(\Lie(G_1^{\ad}),\phi)$ or $(\Lie(G^{\ad}),\phi)$ is non-trivial; thus $\mu$ does not factor through either $Z^0(G)$ or $Z^0(G_1)$. We also assume that the following four conditions hold:

\medskip
{\bf (i)} the maximal tori of $G_k$ are generated by cocharacters of $\pmb{GL}_{\bar M}$ of weights $\{-1,0\}$;

\smallskip
{\bf (ii)} there exists a subset $R_{G_1}$ of $W_{G_1}$ such that the statement 1.4 (a) holds in the context of the family $\{(M,g_1\phi,G_1)|g_1\in G_1(W(k))\}$ of Shimura $F$-crystals over $k$;

\smallskip
{\bf (iii)} each $G_1(k)$-conjugate $\mu_{2k}:\dbG_m\to G_{1k}$ of a cocharacter in the set $\{\sigma^i(\mu_k)|i\in\dbZ\}$ with the property that $\text{Im}(d\mu_{2k})\subseteq\Lie(G_k)$, factors through $G_k$;

\smallskip
{\bf (iv)} either (iv.a) the quotient finite group scheme $Z(G)/Z^0(G)$ is \'etale and there exists a direct sum decomposition $\Lie(G_1)=\Lie(G)\oplus\Lie(G)^\perp$ of $G$-modules which after inverting $p$ is preserved by $\phi$ or (iv.b) $Z^0(G)=Z(G)=Z(G_1)$ and there exists a direct sum decomposition $\Lie(G_1^{\ad})=\Lie(G^{\ad})\oplus \Lie(G^{\ad})^\perp$ of $G^{\ad}$-modules which after inverting $p$ is preserved by $\phi$.

\medskip
Then a subset $R_G$ of $W_G$ as in the statement 1.4 (a) also exists.}
\medskip
 
{\bf Proof:}
Let $g\in G(W(k))$. To prove the Theorem it suffices to show that there exists $w\in W_G$ such that $<\scrC_{\bar g}>=<\scrC_{\bar g_w}>$, cf. Lemma 4.5 (d). To achieve this, we will often replace $g$ by another element $\dag\in G(W(k))$ such that $<\scrC_{\bar g}>=<\scrC_{\bar\dag}>$. As the proof is quite long, its main parts are bold-faced and numbered.

\smallskip
{\bf Part I. Translating (ii) and $\dbF_p$-structures.} 
For $g_1\in G_1(W(k))$, let $\scrC_{1\bar g_1}$ be the $D$-truncation mod $p$ of $(M,g_1\phi,G_1)$. Due to (ii), there exists $w_1\in W_{G_1}$ such that $\scrC_{1\bar g}$ and $\scrC_{1\bar g_{w_1}}$ are inner isomorphic; here $g_{w_1}\in N_{T_1}(W(k))$ represents $w_1$. Let $h\in G_1(W(k))$ be an element such that $\bar h$ belongs to $P_1(k)$ and is an inner isomorphism between $\scrC_{1\bar g}$ and $\scrC_{1\bar g_{w_1}}$. Let $g_0\in G_1(W(k))$ be such that $\bar g_0\in \sigma_{g_{w_1}}(N_1)(k)$ and we have
$$g_{w_1}=g_0hg\phi(h^{-1})\in N_{T_1}(W(k)),$$ 
cf. Lemma 4.5 (a) and (b) applied to $(M,g\phi,G_1)$ and $(M,g_{w_1}\phi,G_1)$. Here $\sigma_{g_{w_1}}:=g_{w_1}\sigma:M\arrowsim M$. For $\tilde h\in\dbL$, we have $<\scrC_g>=<\scrC_{\tilde hg\phi(\tilde h)^{-1}}>$ (cf. Subsubsection 4.2.1). Thus in what follows we will often replace the pair $(g,h)$ by another pair of the form $(\tilde hg\phi(\tilde h)^{-1},h\tilde h^{-1})$.

Let $X_{1\bar g_{w_1}}$ be the reductive subgroup of $G_{1k}$ that is the analogue of $X_{\bar g_w}$ of Theorem 7.6 but obtained working with $(M,g_{w_1}\phi,G_1,\mu)$. Let $X^{\dbF_p}_{1\bar g_{w_1}}$ be the natural $\dbF_p$-structure of $X_{1\bar g_{w_1}}$, cf. Theorem 7.6. Let $T^{g_{w_1}}_{1\dbZ_p}$ be the $\dbZ_p$-structure of $T_1$ defined by $g_{w_1}$, cf. Subsection 4.2 applied to $(M,g_{w_1}\phi,G_1)$. The special fibre $T^{g_{w_1}}_{1\dbF_p}$ of $T^{g_{w_1}}_{1\dbZ_p}$ is a maximal torus of $X^{\dbF_p}_{1\bar g_{w_1}}$. Let $\tilde T^{g_{w_1}}_{\dbZ_p}$ be the smallest subtorus of $T^{g_{w_1}}_{1\dbZ_p}$ such that $\mu$ factors through $\tilde T^{g_{w_1}}:=\tilde T^{g_{w_1}}_{W(k)}$. As the subtorus $\tilde T^{g_{w_1}}_{\dbZ_p}$ of $T^{g_{w_1}}_{1\dbZ_p}$ is uniquely determined by its special fibre $\tilde T^{g_{w_1}}_{\dbF_p}$, $\tilde T^{g_{w_1}}_{\dbF_p}$ is the smallest subtorus of $T^{g_{w_1}}_{1\dbF_p}$ such that $\mu_k$ factors through $\tilde T^{g_{w_1}}_k$; thus we have $\tilde T^{g_{w_1}}_{\dbF_p}\leqslant Z^0(X_{1\bar g_{w_1}})$ (cf. proof of Theorem 7.6). 

Let $X^{\dbF_p}_{1\bar g}:=\bar h^{-1}(X^{\dbF_p}_{1\bar g_{w_1}})\bar h$; it is an $\dbF_p$-structure of the subgroup $X_{1\bar g}:=\bar h^{-1}(X_{1\bar g_{w_1}})\bar h$ of $G_{1k}$. Let 
$$\mu_1:=h^{-1}\mu h:\dbG_m\to G_1.$$ 
Let $\tilde T_{\dbF_p}:=\bar h^{-1}\tilde T^{g_{w_1}}_{\dbF_p}\bar h\leqslant Z^0(X_{1\bar g}^{\dbF_p})$. It is the smallest subtorus of $Z^0(X^{\dbF_p}_{1\bar g})$ such that $\mu_{1k}$ factors through $\tilde T_k$ (i.e., the $\dbF_p$-structure $\tilde T_{\dbF_p}$ of $\tilde T_k$ is defined naturally with respect to $\mu_{1k}$). 

Let $\bar l_1\in\Lie(Z^0(X_{1\bar g}))$ be the image under $d\mu_{1k}$ of the standard generator of $\Lie(\dbG_m)$. We have $\bar l_1=\bar h^{-1}(\bar l_0)$ (see Subsection 4.1 for $\bar l_0$). Let $\grz_{1\bar g_{w_1}}$ be the zero space of $\scrC_{1\bar g_{w_1}}$. The zero space $\grz_{1\bar g}$ of $\scrC_{1\bar g}$ is $\bar h^{-1}(\grz_{1\bar g_0^{-1}\bar g_{w_1}})=\bar h^{-1}(\grz_{1\bar g_{w_1}})$, cf. Lemma 7.1 (a) and (b). As $\Lie(T_{1\dbF_p}^{g_{w_1}})$ (and thus also $\bar l_0$) is included in the $k$-span of $\grz_{1\bar g_{w_1}}$ (cf. Theorem 7.5 (a) applied to $\scrC_{1\bar g_{w_1}}$), the element $\bar l_1=\bar h^{-1}(\bar l_0)$ belongs to the $k$-span of $\grz_{1\bar g}=\bar h^{-1}(\grz_{1\bar g_{w_1}})$.

\smallskip
{\bf Part II. Applying (iv) to study $\bar l_1$.} In this Part II we will show that we can assume that $\bar l_1$ is $\bar l_0$. We first consider the case when (iv.a) holds. Let $\Lie(G_k)^{\perp}:=\Lie(G)^{\perp}/p\Lie(G)^{\perp}$. As $\mu$ factors through $G$ and due to the $G$-module part of (iv.a), we have a direct sum decomposition
$$\Lie(U_{1k})=\Lie(U_k)\oplus (\Lie(U_{1k})\cap \Lie(G_k)^{\perp}).$$ 
The injective map $L_{1g}:\Lie(P_1)+p\Lie(G_1)\hookrightarrow\Lie(G_1)$ that is the analogue of the map $L_g$ of Subsection 7.1 but for $(M,g\phi,G_1)$, is the direct sum of $L_g$ and of an injective map $L_g^{\perp}:(\Lie(P_1)\cap\Lie(G)^{\perp})+p\Lie(G)^{\perp}\hookrightarrow \Lie(G)^{\perp}$.
Similarly, the inclusion map $I_{G_1}:\Lie(P_1)+p\Lie(G_1)\hookrightarrow \Lie(G_1)$ is a direct sum of the inclusion map $I_G$ of Subsection 7.1 and of the inclusion map $I_G^{\perp}:(\Lie(P_1)\cap\Lie(G)^{\perp})+p\Lie(G)^{\perp}\hookrightarrow \Lie(G)^{\perp}$.
From the very definition of zero spaces we get that we have a direct sum decomposition
$$\grz_{1\bar g}=\grz_{\bar g}\oplus (\grz_{1\bar g}\cap \Lie(G_k)^{\perp}).$$ 
\indent
Thus we can write $\bar l_1=\bar l+\bar l^{\perp}$, where $\bar l$ and $\bar l^{\perp}$ belong to $\grw_{\bar g}$ and to the $k$-span of the intersection $\grz_{1\bar g}\cap \Lie(G_k)^{\perp}$ (respectively). For each element $\bar y\in\Lie(G_k)$ we have $[\bar l_1-\bar l,\bar y]\in\Lie(G_k)^\perp$; thus, if $[\bar l_1,\bar y]\in\Lie(G_k)$, then the element $[\bar l_1-\bar l,\bar y]=[\bar l^{\perp},\bar y]$ belongs to $\Lie(G_k)$ as well as to $\Lie(G_k)^{\perp}$ and therefore it is $0$. By applying this to $\bar x\in\Lie(U_k)$, as $[\bar l_1,\bar x]=-\bar x\in\Lie(G_k)$, we get that  $[\bar l,\bar x]=[\bar l_1,\bar x]=-\bar x=[\bar l_0,\bar x]$. 

Let $P'$ be the image in $G^{\ad}$ of $P$. We identify $U$ with the unipotent radical of the parabolic subgroup scheme $P'$ of $G^{\ad}$. Let $I_k$ be the centralizer of $\Lie(U_k)\oplus\Lie(G_k^{0\ad})$ in $\Lie(P'_k)$.

In this paragraph we check that $I_k=\Lie(U_k)$. It suffices to show that: (v) for each simple factor $G_0$ of $G^{\ad}$, the centralizer of $(\Lie(U_k)\oplus\Lie(G_k^{0\ad}))\cap\Lie(G_{0k})$ in $\Lie(P_k'\cap G_{0k})=\Lie(P'_k)\cap\Lie(G_{0k})$ is $\Lie(U_k)\cap\Lie(G_{0k})$. If $G_0\leqslant G^{0\ad}$, then $(\Lie(U_k)\oplus\Lie(G_k^{0\ad}))\cap\Lie(G_{0k})$ and $\Lie(P'_k)\cap\Lie(G_{0k})$ are both $\Lie(G_{0k})$ and $\Lie(U_k)\cap\Lie(G_{0k})=0$; thus (v) is implied by Lemma 2.3 (a). Suppose $G_0\nleqslant G^{0\ad}$ i.e., $\mu$ has a non-trivial image in $G_0$. From this and the fact that $U$ is commutative, we get that $U_k\cap G_{0k}$ is a commutative unipotent radical of the maximal parabolic subgroup $P'_k\cap G_{0k}$ of $G_{0k}$. Thus (v) holds, cf. Lemma 2.3 (b). Thus $I_k=\Lie(U_k)$. 

Let $\bar l_0':=\bar l-\bar l_0\in\Lie(P_k)$. The element $\bar l_0'$ commutes with $\Lie(U_k)$. The image $\bar l_0''$ of $\bar l_0'$ in $\Lie(G^{\ad}_k)$ has a trivial component in $\Lie(G_k^{0\ad})$ (cf. the definition of $G_k^{0\ad}$) and thus it commutes with $\Lie(G_k^{0\ad})$. Therefore $\bar l_0''\in\Lie(P_k')$ commutes with $\Lie(U_k)\oplus \Lie(G_k^{0\ad})$. Thus $l_0''\in I_k=\Lie(U_k)$. Therefore $\bar l_0'\in \Ker(\Lie(G_k)\to\Lie(G_k^{\ad}))+\Lie(U_k)=\Lie(Z(G_k))+\Lie(U_k)$. 

We write $\bar l_0'=\bar l-\bar l_0=\bar x+\bar y$, where $\bar x\in\Lie(Z(G_k))$ and $\bar y\in\Lie(U_k)$; thus $\bar l=\bar l_0+\bar x+\bar y$. Obviously $\Lie(Z(G_{\dbF_p}))=\Lie(Z^0(G_{\dbF_p}))\subseteq\grz_{\bar g}$. Thus $\bar x\in \grw_{\bar g}$. Therefore $\bar l_0+\bar y=\bar l-\bar x$ belongs to $\grw_{\bar g}$ and is the conjugate of $\bar l_1=\bar h^{-1}(\bar l_0)$ through $(1_{\bar M}+\bar y)\bar h\in P_1(k)$. Thus $\bar h(\bar l_0+\bar y)$ and $\bar l_0=\bar h(\bar l_1)$ are $P_1(k)$-conjugate and belong to the $k$-span of $\bar h(\grz_{1\bar g})=\grz_{1\bar g_{w_1}}$. Let $h_1\in P_1(W(k))$ be such that it lifts a $k$-valued point of the group of automorphisms of $\scrC_{1\bar g_{w_1}}$ and we have $\bar h_1\bar h(\bar l_0+\bar y)=\bar l_0$, cf. Theorem 9.1 applied to $(M,g_{w_1}\phi,G_1)$. As $\bar h_1\bar h$ is also an inner isomorphism between $\scrC_{1\bar g}$ and $\scrC_{1\bar g_{w_1}}$, by replacing $h$ with $h_1h$ and $g_0$ by a multiple of it by an element of $G_1(W(k))$ that lifts a $k$-valued point of $\sigma_{g_{w_1}}(N_1)$, we can assume that $\bar l_1=(\bar h)^{-1}(\bar l_0)$ is $\bar l_0+\bar y\in\Lie(P_k)$. This implies that $\bar l^{\perp}=0$ and $\bar l_1=\bar l=\bar l_0+\bar y$. 

Up to a replacement of the triple $(g,h,\bar l_1)$ with the triple $((1_M-y)g\phi(1_M+y),h(1_M+y),(1_{\bar M}-\bar y)(\bar l_1))$, where $y\in\Lie(U)$ lifts $\bar y$, we can assume that $\bar l_1=\bar l$ is $(1_{\bar M}-\bar y)(\bar l_0+\bar y)=\bar l_0+\bar y-\bar y=\bar l_0$.  

If (iv.b) holds, then by working with the images of $\bar l_1$ and $\bar l_0$ in $\Lie(G_{1k}^{\ad})$ instead of $\bar l_1$ and $\bar l_0$, as above we argue that $\bar l=\bar l_0+\bar x+\bar y$, where $\bar x\in\Lie(Z(G_k))$ and $\bar y\in\Lie(U_k)$, and that we can assume $\bar l_1=\bar l_0=\bar l$.  

\smallskip
{\bf Part III. Applying (iii).} As $\bar l_1=\bar l_0$ we have $\mu_{1k}=\mu_k$. Thus $\bar h$ centralizes $\mu_k$. The centralizer of $\bar l_1=\bar l_0$ in $\grz_{1\bar g}$ is $\Lie(X_{1\bar g}^{\dbF_p})$ (after inner conjugation with $\bar h$, this follows from Theorem 7.5 (g) applied to $(M,g_{w_1}\phi,G_1)$) and it is also the direct sum of the centralizers of $\bar l_0$ in $\grz_{\bar g}$ and in $\grz_{1\bar g}\cap\Lie(G_k)^{\perp}$. Thus the smallest $\dbF_p$-vector subspace $\grs$ of $\Lie(X_{1\bar g}^{\dbF_p})$ whose tensorization with $k$ contains $\bar l_0$, is contained in $\grw_{\bar g}$. As $\tilde T_{\dbF_p}\leqslant Z^0(X_{1\bar g}^{\dbF_p})$ and as $\tilde T_k=\bar h^{-1}\tilde T^{g_{w_1}}_k\bar h$ is generated by the cocharacters $\bar h^{-1}\sigma_{g_{w_1}}^i(\mu_{k}) \bar h=\bar h^{-1}\sigma_{g_{w_1}}^i(\mu_{1k}) \bar h$ with $i\in\dbZ$, we get that the Lie algebras of the images of these cocharacters are contained in $\grs$ and thus also in $\Lie(G_k)$. Therefore these images are tori of $G_k$, cf. (iii). Thus $\tilde T_k$ is a torus of $G_k$. As $\tilde T_k$ centralizes $\bar h^{-1}\mu_k\bar h=\mu_k$, it is in fact a torus of $L_k$. As the maximal tori of $L_k$ are $L(k)$-conjugate, up to a replacement of the pair $(g,h)$ with $(h_0g\phi(h_0^{-1}),hh_0^{-1})$, where $h_0\in L(W(k))$ is such that $\bar h_0\tilde T_k\bar h_0^{-1}\leqslant T_k$, we can assume $\tilde T_k\leqslant T_k$. Let $\tilde T$ be the unique subtorus of $T$ that lifts $\tilde T_k$. As $\mu_k=\mu_{1k}$ factors through $\tilde T_k$, the cocharacter $\mu$ of $T$ factors through $\tilde T$.

\smallskip
{\bf Part IV. Applying (i) and toric properties.} Let 
$$\tilde T':=\sigma_g(\tilde T),\;\;\tilde S':=\sigma_g(\tilde TN),\;\; \text{and}\;\;S':=\sigma_g(TN).$$ 
Let $\bar x_0\in\Lie(\tilde T_{\dbF_p})\subseteq\grz_{\bar g}\cap\Lie(\tilde T_k)$. Let $x_0\in\Lie(\tilde T)$ be a lift of $\bar x_0$. There exists $y\in\Lie(N)$ such that $g\phi(x_0+py)-x_0\in p\Lie(G)$, cf. Lemma 7.3 (d). Thus $\bar x_0=\sigma_g(\bar x_0+\bar y)\in\Lie(\tilde S'_k)$. This implies that $\Lie(\tilde T_k)\subseteq\Lie(\tilde S'_k)\subseteq\Lie(S'_k)$. As $\Lie(\tilde T_k)\subseteq\Lie(S'_k)$, from (i) and Proposition 2.4 (b) we get that $\tilde T_k$ centralizes a maximal torus $\tilde T''_k$ of $S'_k$. Thus $\tilde T_k$ is a subtorus of $\tilde T''_k$. Let $\bar g_2\in\sigma_g(N)(k)$ be such that $\bar g_2\tilde T'_k\bar g_2^{-1}$ is a subtorus of $\tilde T''_k$, cf. [1, Ch. III, Thm. 10.6 (4)]. We have $\scrC_{\bar g}=\scrC_{\bar g_2\bar g}$ (cf. Lemma 4.5 (a)) and thus $\grz_{\bar g}=\grz_{\bar g_2\bar g}$ (cf. Lemma 7.1 (b)). We get that $\scrC_{1\bar g}=\scrC_{1\bar g_2\bar g}$. Thus $\bar h$ is an inner isomorphism between $\scrC_{1\bar g_2\bar g}$ and $\scrC_{1\bar g_{w_1}}$. By replacing $(g,g_0)$ with $(g_2g,g_0hg_2h^{-1})$, where $g_2\in G(W(k))$ lifts $\bar g_2$, we can assume $\tilde T_k$ and $\tilde T_k'$ commute.  

By replacing the pair $(g,h)$ with $(h_2 g\phi(h_2^{-1}),hh_2^{-1})$, where $h_2\in L(W(k))$ lifts a $k$-valued point of the centralizer of $\tilde T_k$ in $L_k$ and we have $\bar h_2\tilde T_k' \bar h_2^{-1}\leqslant T_k$, we can assume that $\tilde T_k'\leqslant T_k$. By replacing the pair $(g,g_0)$ with $(g_3g,g_0hg_3h^{-1})$, where $g_3\in\Ker(G(W(k))\to G(k))$ is such that $g_3\tilde T' g_3^{-1}\leqslant T$ (cf. Subsubsection 3.3.2), we can assume $\tilde T'\leqslant T$. As $\bar h \tilde T_k \bar h^{-1}=\tilde T^{g_{w_1}}_k$, from Subsubsection 3.3.2 we also get that there exists $g_4\in\Ker(G_1(W(k))\to G_1(k))$ such that $h\tilde T h^{-1}=g_4\tilde T^{g_{w_1}}g_4^{-1}$. As $g_{w_1}\phi(\Lie(\tilde T^{g_{w_1}}))=\Lie(\tilde T^{g_{w_1}})$ and $g_0^{-1}g_{w_1}\phi(\Lie(h\tilde T h^{-1}))=h\sigma_g\mu({1\over p})h^{-1}(\Lie(h\tilde T h^{-1}))$ is $\Lie(h\sigma_g(\tilde T)h^{-1})=\Lie(h\tilde T' h^{-1})$, we get
$$g_0(\Lie(h\tilde T' h^{-1}))=g_{w_1}\phi(\Lie(h\tilde T h^{-1}))=g_{w_1}\phi(g_4(\Lie(\tilde T^{g_{w_1}})))$$
$$=g_{w_1}\phi g_4\phi^{-1}g_{w_1}^{-1}g_{w_1}\phi(\Lie(\tilde T^{g_{w_1}}))=g_{w_1}\phi(g_4)g_{w_1}^{-1}(\Lie(\tilde T^{g_{w_1}}))=g_{w_1}\phi(g_4)g_{w_1}^{-1}g_4^{-1}(\Lie(h\tilde T h^{-1})).$$ 
Thus $h\tilde T h^{-1}$ and $h\tilde T' h^{-1}$ are inner conjugate under $g_5:=g_0^{-1}g_4'\in G_1(W(k))$, with $g_4':=g_{w_1}\phi(g_4)g_{w_1}^{-1}g_4^{-1}\in G_1(W(k))$. 
As $g_4\in\Ker(G_1(W(k))\to G_1(k))$ and as the reduction mod $p$ of $\phi(g_4)$ belongs to $N_1(k)$, we have $\bar g_4'\in\sigma_{g_{w_1}}(N_1)(k)$. Thus $\bar g_5=\bar g_0^{-1}\bar g_4'\in\sigma_{g_{w_1}}(N_1)(k)$. As $\bar g_5\bar h\tilde T_k \bar h^{-1}\bar g_5^{-1}=\bar h\tilde T_k'\bar h^{-1}$, we get that $\tilde T_k$ and $\tilde T_k'$ are $(h^{-1}\sigma_{g_{w_1}}(N_1)h)(k)$-conjugate. As $\bar h^{-1}$ is an inner isomorphism between $\scrC_{1\bar g_{w_1}}$ and $\scrC_{1\bar g}$, the group of elements $\bar *\in G_1(k)$ such that $\scrC_{1\bar g}=\scrC_{1\bar *\bar g}$ is on one hand $[h^{-1}\sigma_{g_{w_1}}(N_1)h](k)$ and on the other hand is $\sigma_g(N_1)(k)$ (see Lemma 4.5 (a)). Thus the two tori $\tilde T_k$ and $\tilde T_k'$ are $\sigma_g(N_1)(k)$-conjugate. Therefore $\tilde T_k$ and $\tilde T_k'$ are tori of $\sigma_g(TN_1)_k$; moreover their images in $(\sigma_g(TN_1)/\sigma(N_1))_k$ coincide. The tori $\tilde T_k$ and $\tilde T_k'$ commute (being subtori of $T_k$). From the last two sentences we get that $\tilde T_k$ and $\tilde T_k'$ are the same torus of $\sigma_g(TN_1)_k$ and thus also of $T_k$. Thus the two tori $\tilde T$ and $\tilde T'$ of $T$ coincide. Therefore $\sigma_g$ normalizes $\Lie(\tilde T)$.

\smallskip
{\bf Part V. Applying Lemma 4.2.} Let $\tilde G$ be the centralizer of $\tilde T$ in $G$. As $\tilde T\leqslant T$ and as $\mu$ factors through $\tilde T$, $\mu$ factors also through $Z^0(\tilde G)$. As $\sigma_g$ normalizes $\Lie(\tilde T)$, it also normalizes $\Lie(\tilde G)$. Thus $g\phi=\sigma_g\mu({1\over p})$ normalizes $\Lie(\tilde G)$ and therefore the triple $(M,g\phi,\tilde G)$ is a Shimura $F$-crystal over $k$. As $\mu$ factors through $Z^0(\tilde G)$, for each $\tilde g\in \tilde G(W(k))$ there exists an inner isomorphism between $(M,g\phi,\tilde G)$ and $(M,\tilde gg\phi,\tilde G)$ (cf. Lemma 4.2) and thus we have $<\scrC_{g}>=<\scrC_{\tilde gg}>$. We take $\tilde g$ such that $\tilde gg\phi$ normalizes $\Lie(T)$ (to be compared with Subsection 4.2). Thus $\tilde gg\in N_T(W(k))$. We can assume there exists $w\in W_G$ such that $\tilde gg=g_w$. We have $<\scrC_{g}>=<\scrC_{\tilde gg}>=<\scrC_{g_w}>$. This ends the argument for the existence of $w\in W_G$.\endproof

\bigskip
\noindent
{\boldsectionfont \S11. On Basic Theorem C}

\bigskip
In Subsection 11.1 we prove Basic Theorem C. A fast way to achieve this, is to first use Basic Theorems A and B and Corollary 6.2 to prove the shifting process which says that it suffices to prove Basic Theorem C for $p>>0$ and then to show that Theorem 10.1 applies if $p>>0$. For the sake of completeness, we will formalize and use the shifting process only when it is truly required (i.e., only in the Case 5 below). See Corollary 11.1 and Subsection 11.2 for a Corollary and three functorial complements. 

We use the notations of Subsections 4.1, 4.2.1, and 5.1. Let $\mu':\dbG_m\to G^{\ad}$ be the composite of $\mu:\dbG_m\to G$ with $G\twoheadrightarrow G^{\ad}$. Let $P'$, $T'$, and $N_{T'}$ be the images of $P$, $T$, and $N_T$ in $G^{\ad}$. Let $T'_{\dbZ_p}:=\im(T_{\dbZ_p}\to G^{\ad}_{\dbZ_p})$; it is a $\dbZ_p$-structure of the maximal torus $T'$ of $G^{\ad}$. Let $\Tr:\End_{W(k)}(M)\times \End_{W(k)}(M)\to W(k)$ be the trace map. 

\medskip
\noindent
{\bf 11.1. Proof of Basic Theorem C} 

We now prove Basic Theorem C. Let $\dbT_{G_k^{\ad},P'_k,T_k',\sigma}$ be the action defined as in Subsubsection 2.2.1 for the Frobenius endomorphism $\sigma$ of $G^{\ad}_k=G^{\ad}_{\dbF_p}\times_{\dbF_p} k$. We have a product decomposition (cf. (6))
$$\dbT_{G_k^{\ad},P'_k,T_k',\sigma}=\prod_{j\in J} \dbT_{G_k^j,P'_k\cap G^j_k,T'_k\cap G^j_k,\sigma}.$$ 
To prove Basic Theorem C is equivalent to checking that Conjecture 2.1 holds for the quadruple $(G_k,P_k,T_k,\sigma)$, cf. Lemma 5.1 (b). Based on this and Proposition 2.2 (d), to prove Basic Theorem C is equivalent to checking that for each $j\in J$ the Conjecture 2.1 holds for $(G_k^j,P'_k\cap G^j_k,T'_k\cap G^j_k,\sigma)$. For this we can assume that $J\neq\emptyset$. Let $\tilde G^j$ be the reductive, closed subgroup scheme of $G$ that is the inverse image of $G^j_{W(k)}\times_{\Spec(W(k))}\prod_{j'\in J\setminus\{j\}} T'\cap G^{j'}_{W(k)}$ through $G\twoheadrightarrow G^{\ad}$. The Lie algebra $\Lie(\tilde G^j_{B(k)})=\Lie(T_{B(k)})+\Lie(G^j_{B(k)})$ is normalized by $\phi$. As $\mu$ factors through $T$, it also factors through $\tilde G^j$. Thus the triple $(M,\phi,\tilde G^j)$ is a Shimura $F$-crystal over $k$. The Conjecture 2.1 holds for $(G_k^j,P'_k\cap G^j_k,T'_k\cap G^j_k,\sigma)$ if and only if Basic Theorem C holds for the family $\{(M,g_j\phi,\tilde G^j)|g_j\in\tilde G^j(W(k))\}$ of Shimura $F$-crystals over $k$, cf. Proposition 2.2 (d) and Lemma 5.1 (b). Thus by replacing $G$ with $\tilde G^j$ for some $j\in J$, it suffices to prove Basic Theorem C under the extra assumption that:

\medskip
{\bf (i)} the set $J$ has only one element (i.e., the adjoint group scheme $G^{\ad}_{\dbZ_p}$ is simple). 

\medskip
Moreover we have the following replacement process. 

\medskip
{\bf (ii)} To prove Basic Theorem C (i.e., to prove Conjecture 2.1 for $(G_k^{\ad},P'_k,T'_k,\sigma)$), we can replace $(M,\phi,G,\mu)$ by any other analogue quadruple $(M_*,\phi_*,G_*,\mu_*)$ for which we can identity $(G^{\ad}_{\dbZ_p},T'_{\dbZ_p},\mu')=(G^{\ad}_{*\dbZ_p},T'_{*\dbZ_p},\mu_*')$.

\medskip
Based on Lemma 4.2 we can also assume that $\mu$ does not factor through $Z(G)$; thus both $\mu':\dbG_m\to G^{\ad}$ and $(\Lie(G^{\ad}),\phi)$ are non-trivial. As $\mu$ is a cocharacter of $\pmb{GL}_M$ of weights $\{-1,0\}$, for each $i\in\dbZ$ the image of $\sigma^i(\mu)$ in a simple factor $G_0$ of $G^{\ad}$ is either a minuscule cocharacter of $G_0$ (i.e., it acts on $\Lie(G_0)$ via precisely the trivial, the identical, and the inverse of the identical character of $\dbG_m$) or trivial. 

To prove Basic Theorem C it suffices to prove the statement 1.4 (a) (cf. Corollary 6.2) and thus (cf. (ii)) it suffices to show that by keeping fixed the triple $(G_{\dbZ_p}^{\ad},T'_{\dbZ_p},\mu')$, we can choose the family $\{(M,g\phi,G)|g\in G(W(k))\}$ of Shimura $F$-crystals over $k$ such that the set $R_G$ exists for it. Starting from the triple $(G_{\dbZ_p}^{\ad},T'_{\dbZ_p},\mu')$, most often we will only apply Theorem 10.1 and moreover:

\medskip
{\bf (iii.a)} make a good choice for $M_{\dbZ_p}$ and $G_{\dbZ_p}$ and for a reductive, closed subgroup scheme $G_{1\dbZ_p}$ of $\pmb{GL}_{M_{\dbZ_p}}$ which contains $G_{\dbZ_p}$; implicitly $M$, $G$, and $G_1$ will be $M_{\dbZ_p}\otimes_{\dbZ_p} W(k)$, $G_{W(k)}$, and $G_{1W(k)}$ (respectively);

\smallskip
{\bf (iii.b)} take the maximal torus $T_{\dbZ_p}$ of $G_{\dbZ_p}$ to be the inverse image of $T'_{\dbZ_p}$ in $G_{\dbZ_p}$;

\smallskip
{\bf (iii.c)} take $\mu:\dbG_m\to G$ to be the unique cocharacter that lifts $\mu':\dbG_m\to G^{\ad}$, that defines a cocharacter of $\pmb{GL}_M$ of weights $\{-1,0\}$, and that does not act through a non-trivial scalar multiplication on any direct summand of $M$ which is a $G$-module; implicitly $\phi$ will be $(1_{M_{\dbZ_p}}\otimes\sigma)\mu({1\over p})$. 

\medskip
It is always easy to check that the condition (i) of Theorem 10.1 and the parts of the condition (iv) of Theorem 10.1 which refer to $Z^0(G)$, $Z(G)$, and $Z(G_1)$ hold. Thus below we will not refer to these parts of the conditions of Theorem 10.1. For checking a stronger form of the condition (iii) of Theorem 10.1, let $\mu_{2k}:\dbG_m\to \pmb{GL}_{\bar M}$ be a cocharacter of $\pmb{GL}_{\bar M}$ of weights $\{-1,0\}$ which factors through $G_{1k}$ and for which we have $\text{Im}(d\mu_{2k})\subseteq\Lie(G_k)$.

Let $k_0\supseteq \dbF_p$ be the finite field extension such that $G^{\ad}_{\dbZ_p}$ is $\Res_{W(k_0)/\dbZ_p} \tilde G_{W(k_0)}$, with $\tilde G_{W(k_0)}$ an absolutely simple adjoint group scheme over $\Spec(W(k_0))$ (cf. (i) and Subsection 3.1). Let $k_1$ be the smallest field extension of $k_0$ such that $\tilde G_{W(k_1)}$ is split. As $\tilde G_{W(k_0)}$ is not of ${}^3D_4$ Lie type (cf. [30, Cor. 2, p. 182]), we have $[k_1:k_0]\le 2$.

Based on the type $\tau$ of $(\Lie(G^{\ad}),\phi)$ (see Subsection 4.4), we will distinguish six disjoint Cases.
 
\smallskip\noindent
{\bf Case 1: $\tau=A_n$ and $k_1=k_0$.} Let $M_0:=W(k_0)^{n+1}$. We identify $\tilde G_{W(k_0)}^{\sic}=\pmb{SL}_{M_0}$. Let $M_{\dbZ_p}$ be $M_0$ but viewed as a $\dbZ_p$-module. We identify $\Res_{W(k_0)/\dbZ_p} \tilde G^{\sic}_{W(k_0)}(\dbZ_p)$ (resp. $\Res_{W(k_0)/\dbZ_p} Z(\pmb{GL}_{M_0})(\dbZ_p)$) with the group of $W(k_0)$-linear automorphisms (resp. $W(k_0)$-linear scalar automorphisms) of $M_{\dbZ_p}$. Thus we can view $\Res_{W(k_0)/\dbZ_p} \tilde G^{\sic}_{W(k_0)}$ and $\Res_{W(k_0)/\dbZ_p} Z(\pmb{GL}_{M_0})$ as closed subgroup schemes of $\pmb{GL}_{M_{\dbZ_p}}$; let $G_{\dbZ_p}$ be the reductive, closed subgroup scheme of $\pmb{GL}_{M_{\dbZ_p}}$ generated by them. Let $G_{1\dbZ_p}:=\pmb{GL}_{M_{\dbZ_p}}$; thus $G_1=\pmb{GL}_M$. The monomorphism $G\hookrightarrow G_1$ is isomorphic to the standard monomorphism $\pmb{GL}_{n+1}^{[k_0:\dbF_p]}\hookrightarrow \pmb{GL}_{(n+1)[k_0:\dbF_p]}$. From this and [30, p. 85] we get that there exists a cocharacter $\mu:\dbG_m\to G$ as in (iii.c). 

The restriction of $\Tr$ to $\Lie(G)\times\Lie(G)$ is perfect, cf. Lemma 3.5 (i) and the description of $G\hookrightarrow G_1$. Thus we have a direct sum decomposition $\Lie(G_1)=\Lie(G)\oplus\Lie(G)^\perp$ of $G$-modules, where $\Lie(G)^\perp$ is the perpendicular of $\Lie(G)$ with respect to $\Tr$. As we have $\Tr(\phi(x),\phi(y))=\sigma(\Tr(x,y))$ for all $x$, $y\in\End_{W(k)}(M)$, the $B(k)$-vector space $\Lie(G)^\perp[{1\over p}]$ is normalized by $\phi$. Thus the condition (iv.a) of Theorem 10.1 holds. Also the condition (ii) of Theorem 10.1 holds, cf. Corollary 4.6. From Proposition 2.4 (b) applied with $S_0=G_k$ and $S_1=\im(\mu_{2k})$, we get that $\im(\mu_{2k})$ commutes with a maximal torus of $G_k$ and thus also with $Z(G_k)$. Thus $\mu_{2k}$ factors through the centralizer of $Z(G_k)$ in $G_{1k}$ which is $G_k$ itself; thus the condition (iii) of Theorem 10.1 holds. We conclude that the set $R_G$ exists, cf. Theorem 10.1.

\smallskip\noindent
{\bf Case 2: $p>2$ and  $\tau=C_n$.} We have $k_0=k_1$ and $\tilde G_{W(k_0)}^{\sic}$ is the $\pmb{Sp}$ group scheme of a symplectic space $(M_0,\psi_0)$ over $W(k_0)$ of rank $2n$. Let $G_{\dbZ_p}:=\Res_{W(k_0)/\dbZ_p} \pmb{GSp}(M_0,\psi_0)$. Let $M_{\dbZ_p}$ be as in Case 1. 

Let $G_{1\dbZ_p}:=\Res_{W(k_0)/\dbZ_p} \pmb{GL}_{M_0}$. Condition (ii) of Theorem 10.1 holds, cf. Case 1 applied to $(M,\phi,G_1)$. As $p>2$, the restriction of $\Tr$ to $\Lie(G^{\der})\times \Lie(G^{\der})=\Lie(G^{\ad})\times \Lie(G^{\ad})$ is perfect (cf. Lemma 3.5 (iv)); thus we have a direct sum decomposition of $G$-modules $\Lie(G_1)=\Lie(G^{\der})\oplus\Lie(G^{\der})^{\perp}$, where $\Lie(G^{\der})^{\perp}$ is the perpendicular of $\Lie(G^{\der})$ with respect to the restriction of $\Tr$ to $\Lie(G_1)\times\Lie(G_1)$. Let $\Lie(G^{\ad})^\perp:=\im(\Lie(G^{\der})^{\perp}\to \Lie(G_1^{\ad}))$. As $\Lie(Z(G_1))\subseteq\Lie(G_1^{\der})^{\perp}$, we have a direct sum decomposition of $G$-modules $\Lie(G_1^{\ad})=\Lie(G^{\ad})\oplus \Lie(G^{\ad})^\perp$. As in Case 1 we argue that $\phi$ normalizes $\Lie(G^{\ad})^\perp[{1\over p}]$. Thus condition (iv.b) of Theorem 10.1 holds. As the faithful representation $G\hookrightarrow \pmb{GL}_M$ is isomorphic to the standard faithful representation of rank $2n[k_0:\dbF_p]$ of a $\pmb{GSp}_{2n}^{[k_0:\dbF_p]}$ group scheme over $\Spec(W(k))$, from [30, p. 186] we get that there exists a cocharacter $\mu:\dbG_m\to G$ as in (iii.c). Due to the description of $G\hookrightarrow \pmb{GL}_M$, to check the condition (iii) of Theorem 10.1 we can assume $k_0=\dbF_p$. Let $\bar l_2$ be the image through $d\mu_{2k}$ of the standard generator of $\Lie(\dbG_m)$. As $\bar l_2\in\Lie(G_k)$, both $\bar l_2(\bar M)$ and $(1_{\bar M}+\bar l_2)(\bar M)$ are maximal isotropic subspaces of $\bar M=M_0\otimes_{\dbZ_p} k$ with respect to the reduction mod $p$ of $\psi_0$. It is easy to see that this implies that $\mu_{2k}$ factors through $G_k$. Thus the condition (iii) of Theorem 10.1 holds. We conclude that the set $R_G$ exists, cf. Theorem 10.1.

\smallskip\noindent
{\bf Case 3: $p>2$ and $\tau=D_n^{\dbH}$.} As $[k_1:k_0]\in\{1,2\}$, $\tilde G_{k_0}$ is the adjoint group of the $\pmb{SO}$ group $\tilde G^{1}_{k_0}$ of a quadratic form on $k_0^{2n}$ that is either $x_1x_2+\cdots+x_{2n-1}x_{2n}$ or $x_1x_2+\cdots+x_{2n-3}x_{2n-2}+x_{2n-1}^2-\bar ax_{2n}^2$ with $\bar a\in k_0$ not a square (cf. [11, Thm. 2.2.6 and Sect. 2.7]). Based on Subsection 3.1, the group scheme $\tilde G_{W(k_0)}$ is the adjoint of the $\pmb{SO}$ group scheme $\tilde G^{1}_{W(k_0)}$ of a quadratic form on $M_0:=W(k_0)^{2n}$ that is either $x_1x_2+\cdots+x_{2n-1}x_{2n}$ or $x_1x_2+\cdots+x_{2n-3}x_{2n-2}+x_{2n-1}^2-ax_{2n}^2$, where $a\in W(k_0)$ lifts $\bar a$. This case is in essence the same as Case 2: we only have to replace $\pmb{GSp}_{2n}$ group schemes by $\pmb{GSO}$ group schemes of rank $n+1$. Only one thing needs to be dealt with separately. See [2, plate IV] for standard weights and nodes associated to the $D_n$ Lie type.  

The existence of $\mu$ as in (iii.c) is implied by the $D_n$ case of [30, p. 186] as follows. Let $\theta$ be the set of morphisms $\Spec(W(k))\to\Spec(W(k_0))$ of $\Spec(\dbZ_p)$-schemes. We identify 
$$G_{1W(k)}=\prod_{e\in \theta} \pmb{GL}_{M_0}\times_{\Spec(W(k_0))} {}_e \Spec(W(k))\;\;\text{and}\;\;G^{\der}=\prod_{e\in \theta} \tilde G^{1e}_{W(k_0)},$$ 
where $\tilde G^{1e}_{W(k_0)}:=\tilde G^1_{W(k_0)}\times_{\Spec(W(k_0))} {}_e \Spec(W(k))\leqslant \pmb{GL}_{M_0}\times_{\Spec(W(k_0))} {}_e \Spec(W(k))$. We first consider the case $n=4$ and $k_1=k_0$. In this case $\mu$ exists for only one out of the three possible choices for $\tilde G^1_{W(k_0)}$. More precisely, the non-trivial images of $\sigma^i(\mu')$'s with $i\in\dbZ$ in a simple factor $\tilde G^{1e\ad}_{W(k)}$ of $G^{\ad}_{W(k)}$, are minuscule cocharacters that correspond (up to $N_T(W(k))$-conjugation and a numbering of the extremal nodes) to the last nodes $3$ and $4$ of the Dynkin diagram of $\tilde G^{1e\ad}_{W(k)}$ with respect to the image of $T'_{W(k)}$ in $\tilde G^{1e\ad}_{W(k)}$ (i.e., are $\varpi_3^{\vee}$ and $\varpi_4^{\vee}$ cocharacters). The fact that only (resp. that at most two) extremal nodes show up is implied by loc. cit. The fact that at least two extremal nodes show up is implied by the definition of the $D_4^{\dbH}$ type, cf. Subsubsection 4.2.1. Thus $\tilde G^{\der}_{W(k_0)}$ is such that the weight of the $\tilde G^{1e}_{W(k)}$-module $M_0\otimes_{W(k_0)} {}_eW(k)$ is $\varpi_1$, cf. loc. cit. The case $i=0$ also shows that $\mu$ exists, cf. loc. cit. 

If $n\ge 5$ or if $n=4$ and $k_1\neq k_0$, then there exists a unique isogeny cover $\tilde G^1_{W(k_0)}$ of $\tilde G_{W(k_0)}$ as in the previous paragraph. The existence of $\mu$ is argued as in the case $n=4$ and $k_1=k_0$ (for $n\ge 5$, the mentioned minuscule cocharacters correspond either to one or to both nodes $n-1$ and $n$ and thus are $\varpi_{n-1}^{\vee}$ and $\varpi_n^{\vee}$ cocharacters).

\smallskip\noindent 
{\bf Case 4: $p>2$ and $\tau$ is $B_n$ or $D_n^{\dbR}$.} We first consider the case $\tau=D_n^{\dbR}$. We fix embeddings $W(k_0)\hookrightarrow W(k_1)\hookrightarrow W(k)$. Below the nodes (resp. the weights) are of the Dynkin diagram of $\tilde G_{W(k)}$ with respect to the  pull back of $T'$ to a maximal torus of $\tilde G_{W(k)}$. If $k_1\neq k_0$ and $n\ge 5$, then the group $\Aut_{W(k_0)}(W(k_1))=\Gal(k_1/k_0)$ of order 2 permutes the two weights of the half spin representations of $\tilde G^{\sic}_{W(k_1)}$. If $n=4$ and $k_1\neq k_0$, then  we consider the two weights that are permuted by $\Gal(k_1/k_0)$ and that define half spin representations. If $n=4$ and $k_1=k_0$, then we consider the two weights $\varpi_{i_1}$ and $\varpi_{i_2}$ that define half spin representations of $\tilde G^{\sic}_{W(k_0)}$ such that the non-trivial cocharacters of $\tilde G^{\ad}_{W(k)}$ defined by $\sigma^i(\mu)$'s with $i\in\dbZ$, are minuscule cocharacters corresponding (up to $N_T(W(k))$-conjugation) to the node $\alpha_{i_3}$, where $i_3$ is such that $\{i_1,i_2,i_3\}=\{1,3,4\}$ (i.e., are $\varpi_{i_3}^{\vee}$ cocharacters). 

Regardless of what $n\ge 4$ and $[k_1:k_0]\in\{1,2\}$ are, there exists a well defined spin representation of $\tilde G^{\sic}_{W(k_0)}$ on a free $W(k_0)$-module $M_1$ of rank $2^n$ (associated to the considered pair of weights). We write $M_1\otimes_{W(k_0)} W(k_1)=M_1^{(1)}\oplus M_1^{(2)}$, where $M_1^{(1)}$ and $M_1^{(2)}$ are the two half spin representations of $\tilde G^{\sic}_{W(k_1)}$. Let $M_{\dbZ_p}$ be $M_1$ but viewed as a $\dbZ_p$-module. Let $\eps:=-(-1)^{[k_1:k_0]}$. Let $\pmb{GGSpin}^{\eps}_{2n}$ be the subgroup scheme of $\pmb{GL}_{M_1}$ generated by $\tilde G^{\sic}_{W(k_0)}$ and the torus whose pull back to $\Spec(W(k_1))$ is $Z(\pmb{GL}_{M_1^{(1)}})\times_{\Spec(W(k_1))} Z(\pmb{GL}_{M_1^{(2)}})$. Let $G_{\dbZ_p}:=\Res_{W(k_0)/\dbZ_p} \pmb{GGSpin}_{2n}^{\eps}$. Let $G_{1\dbZ_p}$ be the $\Res_{W(k_0)/\dbZ_p}$ of the closed subgroup scheme of $\pmb{GL}_{M_1}$ whose pull back to $\Spec(W(k_1))$ is $\pmb{GL}_{M_1^{(1)}}\times_{\Spec(W(k_1))} \pmb{GL}_{M_1^{(2)}}$. If $[k_1:k_0]=2$, then  $G_{1\dbZ_p}$ is the $\Res_{W(k_1)/\dbZ_p}$ of a $\pmb{GL}_{2^{n-1}}$ group scheme over $\Spec(W(k_1))$. If $k_1=k_0$, then $G_{1\dbZ_p}$ is the $\Res_{W(k_0)/\dbZ_p}$ of a $\pmb{GL}^2_{2^{n-1}}$ group scheme over $\Spec(W(k_0))$. Let $\mu:\dbG_m\to G$ be as in (iii.c), cf. [30, p. 186] and the definition of the $D_n^{\dbR}$ type. The condition (ii) of Theorem 10.1 holds, cf. Case 1. 

The restriction of $\Tr$ to either $\Lie(G)\times \Lie(G)$ or $\Lie(G_1)\times \Lie(G_1)$ is perfect; this follows from Lemma 3.5 (v) once we remark that $G_1$ is a product of $\pmb{GL}_{2^{n-1}}$ group schemes  and that the restriction $\Tr_1$ of $\Tr$ to $\Lie(G_1)\times\Lie(G_1)$, is the standard trace form. Thus we have a direct sum decomposition $\Lie(G_1)=\Lie(G)\oplus\Lie(G)^\perp$ of $G$-modules, where $\Lie(G)^\perp$ is the perpendicular of $\Lie(G)$ with respect to $\Tr_1$. As in Case 1 we argue that $\phi$ normalizes $\Lie(G)^\perp[{1\over p}]$. Thus the condition (iv.a) of Theorem 10.1 holds.

To check the condition (iii) of Theorem 10.1 we can assume $k_0=\dbF_p$. Let $G^{\text{n}}$ be the normalizer of $\Lie(G)$ in $G_1$. It is a closed subgroup scheme of $G_1$ that contains $G$. We check that $G^{\text{n}}$ is smooth and $\Lie(G^{\text{n}})=\Lie(G)$. It suffices to show that $\Lie(G^{\text{n}}_k)=\Lie(G_k)$. Let $\bar x\in\Lie(G)^\perp/p\Lie(G)^\perp$ be such that it normalizes $\Lie(G_k)$. We have $[\bar x,\Lie(G_k)]\subseteq \Lie(G_k)\cap (\Lie(G)^\perp/p\Lie(G)^\perp)=0$; thus $\bar x$ centralizes $\Lie(G_k)$. Let $T^0_k$ be a maximal torus of $G^{\der}_k$. The weights of the $T^0_k$-module $\bar M$ are distinct and their non-zero differences are (as $p>2$) not divisible by $p$ in the $\dbZ$-lattice of characters of $T^0_k$. This can be checked over $\dbC$ and thus it follows from [9, Part III, Exercise 18.4, Prop. 20.15]: the $\dbZ$-lattice of characters of $T^0_k$ is the $\dbZ$-lattice of $\oplus_{i=1}^n \dbQ l_i$ generated by $l_1,\ldots,l_{n-1}$, and ${1\over 2}(l_1+\cdots+l_n)$; the weights are ${1\over 2}(l_1+\cdots+l_n)-\sum_{i\in I} l_i$ with $I$ an arbitrary subset of $\{1,\ldots,n\}$. Due to the fact that the non-zero differences of these weights are not divisible by $p$, the centralizer of $T^0_k$ in $\pmb{GL}_{\bar M}$ is a maximal torus $T_k^1$ of $\pmb{GL}_{\bar M}$. Thus the centralizer of $T^0_k$ in $\End_k(\bar M)$ is $\Lie(T^1_k)$. Therefore $\bar x$ belongs to the Lie algebra of the intersections of all such centralizers $T^1_k$'s. This intersection of tori commutes with all tori of $G_k$ and thus also with $G_k$. For $i\in\{1,2\}$, the weights of the action of $T^0_k$ on $M_1^{(i)}\otimes_{W(k_1)} k$ are permuted transitively by the Weil group of $G^{\sic}_k$ with respect to $T^0_k$ (cf. loc. cit.) and thus the $G_k$-module $M_1^{(i)}\otimes_{W(k_1)} k$ is irreducible. Thus the centralizer of $G_k$ in $\pmb{GL}_{\bar M}$ is $Z(G_k)$. Therefore $\bar x\in\Lie(Z(G_k))\cap (\Lie(G)^\perp/p\Lie(G)^\perp)=0$; thus $\bar x=0$. We conclude that $G^{\text{n}}$ is smooth and that $\Lie(G^{\text{n}})=\Lie(G)$. Thus $G$ is the identity component of $G^{\text{n}}$. But $\mu_{2k}$ normalizes $\Lie(G_k)$ (cf. Proposition 2.4 (b)) and thus it factors through the identity component $G_k$ of $G^{\text{n}}_k$. Therefore the condition (iii) of Theorem 10.1 holds. Therefore the set $R_G$ exists, cf. Theorem 10.1. 

Suppose $\tau=B_n$. As $\tilde G^{\sic}_{W(k_0)}$ is split, we have $k_1=k_0$. Let $M_1$ be the $\tilde G^{\sic}_{W(k_0)}$-module that defines the spin representation of $\tilde G^{\sic}_{W(k_0)}$. Let $\pmb{GSpin}_{2n+1}$ be the subgroup scheme of $\pmb{GL}_{M_1}$ generated by $\tilde G^{\sic}_{W(k_0)}$ and $Z(\pmb{GL}_{M_1})$. Let $G_{\dbZ_p}:=\Res_{W(k_0)/\dbZ_p} \pmb{GGSpin}_{2n+1}$. Let $G_{1\dbZ_p}:=\Res_{W(k_0)/\dbZ_p} \pmb{GL}_{M_1}$. Using [9, Prop. 20.20] as a substitute for [9, Prop. 20.15], the rest of the argument is (i.e., the last two paragraphs are) as above. 

\smallskip\noindent
{\bf Case 5: $p=2$ and $\tau$ is $B_n$, $C_n$, $D_n^{\dbH}$, or $D_n^{\dbR}$.} In this Case we will show directly that there exist elements $w_1,\ldots,w_{[W_G:W_P]}$ of $W_G$ such that: (i) the orbits $o_{\bar g_{w_1}},\ldots,o_{\bar g_{w_{[W_G:W_P]}}}$ of $\dbT_{G_k,\sigma}$ are distinct and (ii) for all $i\in\{1,\ldots,[W_G:W_P]\}$ we have $\dbS(w_i)=d_i^{\text{cl}}-d_u$ (see Subsubsection 4.1.1 and the properties (iv) and (v) of Subsection 6.1 for notations). Based on the Basic Theorem A and Theorem 6.1 (b), this will imply that $\dbO=\{o_{\bar g_{w_i}}|1\le i\le [W_G:W_P]\}$. Thus the set $R_G$ exists, cf. Lemma 5.1 (a) and Lemma 4.5 (d). 

If $p>2$, the existence of the elements $w_1,\ldots,w_{[W_G:W_P]}$ follows from Cases 2 to 4 and Theorem 6.1 (a). To check for $p=2$ the existence of such elements $w_1,\ldots,w_{[W_G:W_P]}$, we first remark that for $w\in W_G$ the permutation $\pi_{w}$ of Subsubsection 4.1.1 is only the natural action of $w\sigma$ on the root system $\Phi$ of Subsubsection 3.3.1. Thus the number $\dbS(w)$ of Subsubsection 4.1.1 can be computed in terms of $w$ and of the automorphism $\sigma:\Phi\arrowsim\Phi$. The relation $\scrR$ on $W_G$ which keeps track on when two orbits $o_{\bar g_{w_{i_1}}}$ and $o_{\bar g_{w_{i_2}}}$ are distinct (see Basic Theorem B; here $w_{i_1},w_{i_2}\in W_G$), is determined by the action of $\sigma$ on $\Phi$ and by the cocharacter $\mu':\dbG_m\to T'_{W(k)}$. Here $\mu'$ is an element of the abelian group of cocharacters $X_*(T'_{W(k)})$ of $T'_{W(k)}$ and $\Phi$ is a subset of the abelian group of characters $X^*(T'_{W(k)})$ of $T'_{W(k)}$. We conclude that the existence of elements $w_1,\ldots,w_{[W_G:W_P]}$ having the desired properties, is encoded in the data provided by the following quintuple 
$$(X_*(T'_{W(k)}),X^*(T'_{W(k)}),\mu',\Phi,\sigma),$$ 
where $\sigma$ is an automorphism of the root system $\Phi$ and thus also of the Weyl group $W_G=W_{G^{\ad}}$ of $\Phi$. But this quintuple does not depend on the prime $p$ i.e., for each odd prime $l$ there exists a Shimura $F$-crystal $(M^{(l)},\phi^{(l)},G^{(l)})$ over an algebraic closure of $\dbF_l$ whose analogue quintuple is $(X_*(T'_{W(k)}),X^*(T'_{W(k)}),\mu',\Phi,\sigma)$. This is a direct consequence of the constructions performed in Cases 2 to 4 which were based only on the data provided by such quintuples. Here are two examples. If $\tau=B_n$, then the $\dbZ_l$-structure $G^{{(l)}\ad}_{\dbZ_l}$ of $G^{{(l)}\ad}$ is the Weyl restriction $\Res_{W(\dbF_{l^{[k_0:\dbF_p]}})/\dbZ_l}$ of a split, absolutely simple adjoint group of $B_n$ Lie type and $(M^{(l)},\phi^{(l)},G^{(l)})$ is obtained as in Case 4 using spin representations. If $\tau=D_n^{\dbH}$ and $[k_1:k_0]=2$, then the $\dbZ_l$-structure $G^{{(l)}\ad}_{\dbZ_l}$ of $G^{{(l)}\ad}$ is the Weyl restriction $\Res_{W(\dbF_{l^{[k_0:\dbF_p]}})/\dbZ_l}$ of a non-split, absolutely simple adjoint group of ${}^2D_n$ Lie type and $(M^{(l)},\phi^{(l)},G^{(l)})$ is obtained as in Case 3 using representations of rank $2n$ and quadratic forms.

Thus the existence of elements $w_1,\ldots,w_{[W_G:W_P]}$ having the properties (i) and (ii), follows from the first sentence of the previous paragraph. Therefore the set $R_G$ exists.  

\smallskip\noindent
{\bf Case 6: $\tau=A_{n-1}$ with $n\ge 3$ and $k_1\neq k_0$.} Let $M_1:=W(k_1)^n$. We identify $\tilde G_{W(k_1)}^{\sic}=\pmb{SL}_{M_1}$. Let $M_{\dbZ_p}$ be $M_1$ but viewed as a $\dbZ_p$-module. We identify naturally $G^{\der}_{\dbZ_p}=\Res_{W(k_0)/\dbZ_p} \tilde G^{\sic}_{W(k_0)}$ with a closed subgroup scheme of $\Res_{W(k_1)/\dbZ_p} \tilde G^{\sic}_{W(k_1)}$ and thus of $\pmb{GL}_{M_{\dbZ_p}}$. Up to a $\dbG_m(W(k_0))$-multiple, there exists a unique perfect, alternating form on $M_{\dbZ_p}$ that is $W(k_0)$-linear and that is fixed by $G^{\der}$. This implies that there exists a unique reductive subgroup scheme $G_{1\dbZ_p}$ of $\pmb{GL}_{M_{\dbZ_p}}$ that contains $G^{\der}_{\dbZ_p}$ and such that $G_{1W(k_0)}$ is a $\pmb{GSp}_{2n}^{[k_0:\dbF_p]}$ group scheme over $\Spec(W(k_0))$. Let $G_{\dbZ_p}$ be generated by $G_{\dbZ_p}^{\der}$ and by the center of the centralizer of $G_{\dbZ_p}^{\der}$ in $G_{1\dbZ_p}$. The monomorphism $G^{\der}\hookrightarrow G_1^{\der}$ is a product of $[k_0:\dbF_p]$ copies of the standard monomorphism $\pmb{SL}_{n}\hookrightarrow \pmb{Sp}_{2n}$. The ranks of $Z(G)$ and $Z(G_1)$ are $2[k_0:\dbF_p]$ and $[k_0:\dbF_p]$ (respectively). From the last two sentences we get that the ranks of $G$ and $G_1$ are $(n+1)[k_0:\dbF_p]$. Thus $T=T_1$. Let $\mu:\dbG_m\to G$ be as in Case 1. Let $\Lie(G_1)=\Lie(T)\oplus\bigoplus_{\alpha\in\Phi_1} \grg_{\alpha}$ be the root decomposition relative to $T$; we have $\Phi\subseteq \Phi_1$. The $T$-submodule
$$\Lie(G)^\perp:=\oplus_{\alpha\in\Phi_1\setminus\Phi} \grg_{\alpha}$$ 
of $\Lie(G_1)$ is uniquely determined by the identity $\Lie(G_1)=\Lie(G)\oplus\Lie(G)^\perp$. As $\phi$ normalizes $\Lie(T)$, $\Lie(G_{B(k)})$, and $\Lie(G_{1B(k)})$, it also normalizes $\Lie(G)^\perp[{1\over p}]$ i.e., the condition (iv.a) of Theorem 10.1 holds. Let $\alpha\in\Phi$ and $\alpha_1\in\Phi_1\setminus\Phi$ be such that $\alpha+\alpha_1\in\Phi_1$. If $\alpha+\alpha_1\in\Phi$, then $\grg_{\alpha_1}[{1\over p}]=[\grg_{-\alpha}[{1\over p}],\grg_{\alpha+\alpha_1}[{1\over p}]]\subseteq\Lie(G)[{1\over p}]$ (cf. Chevalley rule for the equality part) contradicts the relation $\alpha_1\in\Phi_1\setminus\Phi$. Thus $(\Phi+(\Phi_1\setminus\Phi))\cap \Phi_1\subseteq\Phi_1\setminus\Phi$. From this and the Chevalley rule we get that the direct sum decomposition $\Lie(G_1)=\Lie(G)\oplus\Lie(G)^\perp$ is left invariant by $\Lie(G)$. Thus this direct sum decomposition is normalized by any maximal torus of $G$ and therefore also by $G$ itself. Thus $\Lie(G)^\perp$ is a $G$-module. The centralizer of a maximal torus of $G_k$ in $G_{1k}$ is this torus itself and thus it is contained in $G_k$. From this and Proposition 2.4 (b) we get that $\mu_{2k}$ factors through $G_k$. Thus the condition (iii) of Theorem 10.1 holds. Also the condition (ii) of Theorem 10.1 holds, cf. Case 2 if $p>2$ and cf. Case 5 if $p=2$. The rest is as in Case 1 i.e., Theorem 10.1 implies that the set $R_G$ exists. This ends the proof of the statement 1.4 (a) and thus also of Basic Theorem C and of Corollaries 1.1 and 1.2.\endproof

\smallskip
\noindent
{\bf Corollary 11.1.} {\it  {\bf (a)} Conjecture 2.1 holds for $(\scrG,\scrP,\scrT,F)=(G_k,P_k,T_k,\sigma)$. 

\smallskip
{\bf (b)} We have $\dbO=\{o_{\bar g_w}|w\in W_G\}$ and the set $\dbO$ has precisely $[W_G:W_P]$ elements. 

\smallskip
{\bf (c)} There exists a unique equivalence class $[w]\in \scrR\backslash W_G$ called the pivotal class (equivalently there exists a unique orbit $o_{\bar g_w}\in\dbO$ called the pivotal orbit) such that $\dim(o_{\bar g_w})=\dim(P_k)$ (equivalently such that $\dim(o_{\bar g_w})\le\dim(P_k)$). Moreover, for $g\in G(W(k))$ we have $<\scrC_g>=<\scrC_{g_w}>$ if and only if $<\scrC_{\bar g}>=<\scrC_{\bar g_w}>$ (i.e., with the notations of Corollary 1.2, the subset $\text{Mod}_p^{-1}(\scrC_{\bar g_w})$ of $\dbI_{\infty}$ has only one element $\text{Lift}(\scrC_{\bar g_w})$).

\smallskip
{\bf (d)} Let $g\in G(W(k))$. Then $<\scrC_g>\in \im(\text{Lift})$ if and only if there exists a maximal torus $\tilde T$ of $G$ such that $g\phi$ normalizes $\Lie(\tilde T)$. Thus $\im(\text{Lift})$ does not depend on the choices of $\mu:\dbG_m\to G$ and $T$; therefore the section $\text{Lift}:\dbI\hookrightarrow \dbI_{\infty}$ of $\text{Mod}_p:\dbI_{\infty}\twoheadrightarrow\dbI$ is indeed canonical.}
\medskip

{\bf Proof:}
Part (a) follows from Basic Theorem C and Lemma 5.1 (b). Part (b) only translates Corollary 1.1 and Lemma 5.1 (a). Based on (a) and the fact that the composite map $b_{\text{can}}\circ a_{\text{can}}:\scrR\backslash W_G\to \dbO$ is an isomorphism, (c) only combines Remark 6.4 (c) with Theorem 8.3.

We check (d). As $g_w\phi$ normalizes $\Lie(T)$ for all $w\in W_G$ (cf. Subsubsection 4.2.1), we only need to check that if there exists a maximal torus $\tilde T$ of $G$ such that $g\phi$ normalizes $\Lie(\tilde T)$, then there exists $w\in W_G$ such that $<\scrC_g>=<\scrC_{g_w}>$. As $g\phi$ normalizes $\Lie(\tilde T)$, we have $\Lie(\tilde T)\subseteq p\Lie(G)+\Lie(P)=\Lie(G)\cap (g\phi)^{-1}(\Lie(G))$. Thus $\Lie(\tilde T_k)\subseteq\Lie(P_k)$. Let $A_k'$ and $A_k''$ be subgroups of $G_k$ as in the proof of Theorem 7.6. Let $\tilde T_k':=\tilde T_k\cap A_k'$ and $\tilde T_k'':=\tilde T_k\cap A_k''$; they are maximal tori of $A_k'$ and $A_k''$. We have $A_k''\leqslant L_k\leqslant P_k$; thus $\tilde T_k''\leqslant P_k$. The intersection $A_k'\cap P_k$ is a parabolic subgroup of $A_k'$ whose Lie algebra contains $\Lie(\tilde T_k')$. As $\tilde T_k'$ is generated by cocharacters of $\pmb{GL}_{\bar M}$ of weights $\{-1,0\}$, it centralizes a maximal torus $I_k'$ of $A_k'\cap P_k$ (cf. Proposition 2.4 (b)) and thus it is contained in $I_k'$. As $\tilde T_k$ is generated by $\tilde T_k'$ and $\tilde T_k''$, we get that $\tilde T_k\leqslant P_k$. Let $\bar h\in P(k)$ be such that $\bar h\tilde T_k\bar h^{-1}=T_k$. Let $h\in P(W(k))\leqslant \dbL$ be a lift of $\bar h$ such that $h\tilde Th^{-1}=T$, cf. Subsubsection 3.3.2. By replacing $g$ with $hg\phi(h^{-1})$, we can assume $\tilde T=T$. Thus $g\phi$ and $\phi$ both normalize $\Lie(T)$; thus $g\in N_T(W(k))$. Let $w\in W_G$ be such that $g\in T(W(k))g_w$. We have $<\scrC_g>=<\scrC_{g_w}>$, cf. Lemma 4.1. Thus (d) holds.
\endproof 

\medskip
\noindent
{\bf 11.2. Three functorial complements} 

{\it {\bf (a)}} In Case 6 of Subsection 11.1, as $G_{1\dbZ_p}$ we can also  take $\Res_{W(k_1)/\dbZ_p} \pmb{GL}_{M_1}$. Conditions (i) and (ii) of Theorem 10.1 hold. But the condition (iii) of Theorem 10.1 does not hold if and only if $p=2$, $n$ is even, and $L^{\ad}$ has simple factors of $A_{n\over 2}$ Lie types. The condition (iv) of Theorem 10.1 does not hold if $p=2$ and $n$ is even. 

\smallskip
{\bf (b)} Let $\tilde G$ be a reductive, closed subgroup scheme of $\pmb{GL}_M$ that contains $G$ and $(M,\phi,\tilde G)$ is a Shimura $F$-crystal. Let $\tilde G_{\dbZ_p}$ be the reductive, closed subgroup scheme of $\pmb{GL}_{M_{\dbZ_p}}$ whose pull back to $\Spec(W(k))$ is $\tilde G$, cf. Subsection 4.2. Let $\tilde T_{\dbZ_p}$ be a maximal torus of the centralizer of $T_{\dbZ_p}$ in $\tilde G_{\dbZ_p}$, cf. Proposition 3.3 (a) and (e). Let $\tilde T:=\tilde T_{W(k)}$. Let $N_{\tilde T}$, $W_{\tilde G}$, $\tilde\scrR$, $\tilde {\text{Lift}}:\tilde\dbI\hookrightarrow\tilde\dbI_{\infty}$, $\tilde a_{\text{can}}:\tilde\scrR\backslash W_{\tilde G}\tilde\to\tilde\dbI$, and $\tilde\dbO$ be the analogues of $N_T$, $W_G$, $\scrR$, $\text{Lift}:\dbI\hookrightarrow \dbI_{\infty}$, $a_{\text{can}}:\scrR\backslash W_G\tilde\to\dbI$, and $\dbO$ (respectively) but obtained working with $(M,\tilde g\phi,\tilde G)$'s and $\tilde T$ instead of $\scrC_g$'s and $T$; here $\tilde g\in \tilde G(W(k))$ and $g\in G(W(k))$. If $w\in W_G$, let $\tilde G^{g_w}_{\dbZ_p}$ be the reductive, closed subgroup scheme of $\pmb{GL}_{M^{g_w}_{\dbZ_p}}$ whose pull back to $\Spec(W(k))$ is $\tilde G$ (cf. Subsection 4.2) and let $\tilde T^{g_w}_{\dbZ_p}$ be a maximal torus of the centralizer of $T^{g_w}_{\dbZ_p}$ in $\tilde G^{g_w}_{\dbZ_p}$. It is easy to see that $g_w\phi$ normalizes the Lie algebra of the maximal torus $\tilde T^{g_w}_{W(k)}$ of $\tilde G$. From this and Corollary 11.1 (d) we get that for each $w\in W_G$ we have $<(M,g_w\phi,\tilde G)>\in\text{Lift}(\tilde\dbI)$.

The maps $b_{\text{can}}:\dbI\arrowsim\dbO$, $\text{Mod}_p:\dbI_{\infty}\twoheadrightarrow\dbI$, and $\text{Lift}:\dbI\hookrightarrow \dbI_{\infty}$ are functorial. In particular, we have an identity 
$$\tilde {\text{Lift}}\circ \text{Func}=\text{Func}_{\infty}\circ \text{Lift}:\dbI\to\tilde\dbI_{\infty},$$ 
where $\text{Func}:\dbI\to\tilde\dbI$ and $\text{Func}_{\infty}:\dbI_{\infty}\to\tilde\dbI_{\infty}$ are defined by $\text{Func}(<\scrC_{\bar g}>)=<(\bar M,\bar g\bar\phi,\bar\vartheta\bar g^{-1},\tilde G_k)>$ and $\text{Func}_{\infty}(<\scrC_g>)=<(M,g\phi,\tilde G)>$. Let $a_{\text{can}}:\scrR\backslash W_G\to\dbI$ be as in Corollary 1.1. 

If moreover $\tilde T$ is the centralizer of $T$ in $\tilde G$, then $W_G$ is naturally a subgroup of $W_{\tilde G}$. In this case we also have 
$$\text{Func}\circ a_{\text{can}}=\tilde a_{\text{can}}\circ\text{Nat}:\scrR\backslash W_G\to \tilde\dbI,$$
where $\text{Nat}:\scrR\backslash W_G\to \tilde\scrR\backslash W_{\tilde G}$ is the natural map induced by the inclusion $W_G\hookrightarrow W_{\tilde G}$. 

\smallskip
{\bf (c)} Referring to Corollary 1.1, one can use the ideas of either [22, Thms. 4.7 and 5.5] or [20, Sect. 2] to show that there exists a natural bijection $b_G:\scrR\backslash W_G\arrowsim W_P\backslash W_G$. We emphasize that the definition of $b_G$ is complicated and that $b_G$ is rarely functorial and thus it is not truly canonical. 

\bigskip
\noindent
{\boldsectionfont \S12. Applications to stratifications}

\bigskip
We use the notations of Subsection 4.1. See Subsection 5.1 for the action $\dbT_{G_k,\sigma}$ and for the topological space $\dbO$ of orbits of $\dbT_{G_k,\sigma}$ endowed with the order $\le $. Let $Y$ be a separated, reduced $k$-scheme. Let $F_Y:Y\to Y$ be its Frobenius endomorphism. Let $\bar\scrM_Y$ be the pull back to $Y$ of the vector bundle over $\Spec(k)$ defined by $\bar M$. If $\scrB$ is a vector bundle over $Y$, let $\scrB_*$ be its pull-back to a vector bundle over a $Y$-scheme $*$. 

Subsection 12.1 defines analogues of $\scrC_g$'s over $Y$ called $D$-bundles mod $p$ over $Y$; a standard example is presented in Example 12.2. Basic Theorem D and Corollary 12.3 assume $Y$ is of finite type over $k$ and list the main properties of the stratifications of $Y$ defined by $D$-bundles mod $p$ over $Y$. We end with two remarks.

\medskip
\noindent
{\bf 12.1. $D$-bundles mod $p$} 

By a $D$-bundle mod $p$ over $Y$ associated to $\{\scrC_g|g\in G(W(k))\}$ we mean a quadruple $\scrW:=(\scrK,\scrF,\scrV,\scrA)$, where

\medskip
{\bf (i)} $\scrK$ is a vector bundle over $Y$ of rank $r$,

\smallskip
{\bf (ii)} $\scrF:F_Y^*(\scrK)\to\scrK$ and $\scrV:\scrK\to F_Y^*(\scrK)$ are homomorphisms of vector bundles,

\smallskip
{\bf (iii)} $\scrA=\{(U_i,\rho_i)|i\in I\}$ is an atlas of maps such that the following four things hold:

\medskip\noindent
{\bf (iii.a)} each $U_i=\Spec(R_i)$ is an \'etale $Y$-scheme which is an affine scheme  and $\rho_i:\scrK_{U_i}\arrowsim \bar\scrM_{U_i}$ is an isomorphism with respect to which there exists a morphism $\bar g_i:U_i\to G_k$ such that the pull back of $\scrF$ to $U_i$ becomes an $R_i$-linear map $\bar M\otimes_k {}_{\sigma} R_i\to \bar M\otimes_k R_i$ that maps $x\otimes 1$ to $\bar g_i(\bar{\phi}(x)\otimes 1)$ and the pull back of $\scrV$ to $U_i$ becomes an $R_i$-linear map $\bar M\otimes_k R_i\to \bar M\otimes_k {}_{\sigma} R_i$ that maps $\bar g_i(\bar x\otimes 1)$ to $\bar \vartheta(x)\otimes 1$, where $\bar x\in\bar M$ and where we identify $\bar g_i$ with an element of $G_k(R_i)$;

\smallskip\noindent
{\bf (iii.b)} if $i$, $j\in I$ and $U_{ij}:=U_i\times_Y U_j=\Spec(R_{ij})$, then $\rho_{jU_{ij}}\circ \rho_{iU_{ij}}^{-1}$ is the automorphism of $\bar\scrM_{U_{ij}}$ defined by an element of $P_k(R_{ij})$;

\smallskip\noindent
{\bf (iii.c)} the set $\{\im(U_i\to Y)|i\in I\}$ is an open cover of $Y$;

\smallskip\noindent
{\bf (iii.d)} it is not strictly included in any other atlas satisfying the properties (iii.a) to (iii.c).

\medskip
Due to (i) and (iii.b), $\scrK$ has a natural structure of a $P_k$-bundle over $Y$. Let $y\in Y(k)$. Let $i\in I$ be such that there exists a point $y_i\in U_i(k)$ that maps to $y$, cf. (iii.c). Let $\bar g_{iy_i}:=\bar g_i\circ y_i\in G_k(k)$. Let $\scrC_{iy_i}:=(\bar M,\bar g_{iy_i}\bar{\phi},\bar \vartheta\bar g_{iy_i}^{-1},G_k)$. Due to (iii.b), the inner isomorphism class $<\scrC_{iy_i}>\in\dbI$ depends only on $y$ and not on the choices of $i\in I$, $y_i\in U_i(k)$, and $\bar g_i$; thus we can define $<\scrC_y>:=<\scrC_{iy_i}>$ as well as $o(y):=o_{\bar g_{iy_i}}\in\dbO$. Let $\dbS(y)$ be the dimension of the reduced stabilizer subgroup of $\bar g_{iy_i}$ under $\dbT_{G_k,\sigma}$. 

\smallskip
\noindent
{\bf Definition 12.1.} 
{\bf (a)} We say $\scrW$ is {\it uni+versal} if for each point $y\in Y(k)$ there exists a map $(U_i,\rho_i)\in\scrA$ with the properties that: (i) $y\in \im(U_i(k)\to Y(k))$ and (ii) we can choose $\bar g_i:U_i\to G_k$ such that its composite with the quotient $k$-morphism $G_k\to P_k\backslash G_k$, is \'etale. 

\smallskip
{\bf (b)} We say  $\scrW$ is {\it full} if for each point $y\in Y(k)$ there exists a pair $(U_i,\rho_i)\in\scrA$ such that $y\in \im(U_i(k)\to Y(k))$ and we can choose $\bar g_i$ to be \'etale. 

\medskip
If $\scrW$ is uni+versal (resp. full), then $Y$ is a smooth $k$-scheme which is equidimensional of dimension $\dim(N_k)=\dim(P_k\backslash G_k)$ (resp. of dimension $\dim(G_k)$).

\smallskip
\noindent
{\bf 12.1.1. Augmentations} 

Suppose $\scrW$ is uni+versal. We will augment $Y$ and $\scrW$ to get naturally a $D$-bundle mod $p$ that is full. Let $\tilde Y:=P_k\times_k Y$. Let $q_1:\tilde Y\to P_k$ and $q_2:\tilde Y\to Y$ be the two projections. Let $\tilde{\scrW}=(\tilde{\scrK},\tilde{\scrF},\tilde{\scrV},\tilde{\scrA})$ be the pull back of $\scrW$ via $q_2$, with $\tilde{\scrA}$ a maximal atlas which contains $q_2^*(\scrA)$. Let $i\in I$. Let $h_i:P_k\times_k U_i\to G_k$ be the composite of the first projection $P_k\times_k U_i\twoheadrightarrow P_k$ with the inclusion $P_k\hookrightarrow G_k$. Let $I_{h_i}:\bar\scrM_{P_k\times_k U_i}\arrowsim\bar\scrM_{P_k\times_k U_i}$ be the automorphism defined by inner conjugation through $h_i$. Let $\tilde\rho_i:=I_{h_i}\circ (q_2^*(\rho_i)):\tilde{\scrK}_{P_k\times_k U_i}\arrowsim\bar\scrM_{P_k\times_k U_i}.$

Let $\tilde{\bar g}_i:=h_i(\bar g_i\circ q_2)h_{i0}^{-1}:P_k\times_k U_i\to G_k$, where $h_{i0}$ is the composite of the restriction of $q_1\circ F_{\tilde Y}$ to $P_k\times_k U_i$ with the epimorphism $P_k\twoheadrightarrow P_k/U_k=L_k$ and with the monomorphism $L_k\hookrightarrow G_k$. If the composite of $\bar g_i:U_i\to G_k$ with the quotient $k$-morphism $G_k\to P_k\backslash G_k$ is \'etale, then the morphism $\tilde{\bar g}_ih_{i0}=h_i(\bar g_i\circ q_2):P_k\times_k U_i\to G_k$ is also \'etale; as the images of the tangent maps of $\tilde{\bar g}_i$ and $\tilde{\bar g}_ih_{i0}$ at an arbitrary $k$-valued point of $P_k\times_k U_i$ have the same dimensions, the morphism $\tilde{\bar g}_i:P_k\times_k U_i\to G_k$ is also \'etale. This implies that $\tilde{\scrW}$ is full. 

\smallskip
\noindent
{\bf Example 12.2.}
 Suppose $r=2d$. Let $d_0\in\dbN\cup\{0\}$ be the relative dimension of $N$ over $\Spec(W(k))$. Let $Z$ be a smooth scheme over $\Spec(W(k))$ of relative dimension $d_0$. Let $Y:=Z_k$. Let $\grA$ be an abelian  scheme over $Z$ of relative dimension $d$. Let $\scrM_0:=\underline{H}^1_{\text{dR}}(\grA/Z)$; it is a vector bundle over $Z$ of rank $r$. Let $\pmb{GL}_{\scrM_0}$ be the general linear group scheme over $Z$ defined by $\scrM_0$. Let $\scrK$ be the reduction mod $p$ of $\scrM_0$; it is a vector bundle over $Y$ of rank $r$. The evaluation of the Dieudonn\'e crystal of $\grA$ at the trivial thickening of $Y$, (when viewed without connection) consists of two pairs $\scrF:F_Y^*(\scrK)\to\scrK$ and $\scrV:\scrK\to F_Y^*(\scrK)$ of homomorphisms of vector bundles. Suppose that the vector bundle $\scrM_0$ over $Z$ has in the \'etale topology of $Z$ a structure of a $G$-bundle. 

Let $C_y$ be the completion of the local ring of a point $y\in Y(k)$ in $Z$. Let $Z_y=\Spec(C_y)$. We fix an identification $C_y=W(k)[[x_1,\ldots,x_{d_0}]]$ and we consider the Frobenius lift $\Phi_{C_y}$ of $C_y$ that is compatible with $\sigma$ and that takes $x_i$ to $x_i^p$ for all $i\in\{1,\ldots,d_0\}$. Let $M_{0y}:=H^1_{\text{dR}}(\grA\times_Z Z_y/C_y)$; it is a free $C_y$-module of rank $r$.  Let $\Phi_{y}:M_{0y}\to M_{0y}$ be the $\Phi_{C_y}$-linear endomorphism defined naturally by $\grA\times_{Z} Z_{y}$. Let $\scrM_{Z_y}$ be the $G$-bundle over $Z_y$ defined by $M\otimes_{W(k)} C_y$. We also assume that the following condition holds:

\medskip
{\bf (*)} for each $k$-valued point $y$ of $Y$, there exists an isomorphism $\scrM_0\times_{Z} Z_{y}\arrowsim \scrM_{Z_y}$ of $G$-bundles defined by an isomorphism $M_{0y}\arrowsim M\otimes_{W(k)} C_y$ under which $\Phi_{0y}$ gets identified with the $\Phi_{0y}$-linear endomorphism $n_{\text{univ}}(h_{y}\phi\otimes\Phi_{C_y})$ of $M\otimes_{W(k)} C_y$, where $h_{y}\in G(W(k))$ and where $n_{\text{univ}}:Z_y\to  N$ is formally \'etale. 

\medskip
Based on (*) and Artin's approximation theorem, we get that there exists an \'etale morphism $U_{y}=\Spec(R_{y})\to Y$, an isomorphism $\rho_{y}:\scrK_{U_y}\to \bar\scrM_{U_y}$ of $G_k$-bundles, and a morphism $\bar g_y:U_y\to G_k$ such that the following two things hold: (i) the pull back of $\scrF$ to $U_y$ becomes (via $\rho_y$) an $R_y$-linear map $\bar M\otimes_k {}_{\sigma} R_y\to \bar M\otimes_k R_y$ that maps $x\otimes 1$ to $\bar g_y(\bar{\phi}(x)\otimes 1)$ and the pull back of $\scrV$ to $U_y$ becomes (via $\rho_y$) an $R_y$-linear map $\bar M\otimes_k R_y\to \bar M\otimes_k {}_{\sigma} R_y$ that maps $\bar g_y(\bar x\otimes 1)$ to $\bar \vartheta(x)\otimes 1$, where $\bar x\in\bar M$ and where we identify $\bar g_y$ with an element of $G_k(R_y)$, and (ii) the composite of $\bar g_y:U_y\to G_k$ with the quotient morphism $G_k\to G_k/P_k$ is \'etale.

Let $\scrA=\{(U_i,\rho_i)|i\in I\}$ be a maximal atlas such that $Y(k)\subseteq I$ and the property (iii) of Subsection 12.1 holds for it. The quadruple $\scrW=(\scrK,\scrF,\scrV,\scrA)$ is a uni+versal $\scrD$-bundle mod $p$ over $Y$.

\medskip
\noindent
{\bf 12.2. Basic Theorem D} 

{\it Suppose $Y$ is of finite type over $k$. Let $\scrW:=(\scrK,\scrF,\scrV,\scrA)$ be a $D$-bundle mod $p$ over $Y$ associated to the family $\{\scrC_g|g\in G(W(k))\}$. Then the following six basic properties hold:

\medskip
{\bf (a)} There exists a stratification $\scrY$ of $Y$ in reduced, locally closed subschemes such that $y_1$, $y_2\in Y(k)$ belong to the same stratum of $\scrY$ if and only if $<\scrC_{y_1}>=<\scrC_{y_2}>$  (i.e., if and only if $o(y_1)=o(y_2)$, cf. Lemma 5.1 (a)). 

\smallskip
{\bf (b)} If $\scrW$ is full (resp. uni+versal), then each stratum $\grs$ of $\scrY$ is smooth and equidimensional of dimension $\dim(o(y))=\dim(G_k)-\dbS(y)$ (resp. $\dim(N_k)-\dbS(y)$) with $y\in\grs(k)$. Moreover, $\scrY$ has a unique open stratum.

\smallskip
{\bf (c)} Let $\grs_1$ and $\grs_2$ be two strata of $\scrY$. Let $y_j\in\grs_j(k)$, $j=\overline{1,2}$. Let $o_1'\in \dbO$ be such that $o_1'\le o(y_2)$. If $\grs_1$ specializes to $\grs_2$, then $o(y_1)\le o(y_2)$. If $\scrW$ is either full or uni+versal, then there exists a stratum $\grs_1'$ of $\scrY$ that specializes to $\grs_2$ and such that we have $o_1'=o(y_1')$ for some point $y_1'\in\grs_1'(k)$. In particular, if $\scrW$ is either full or uni+versal, then $\grs_1$ specializes to $\grs_2$ if and only if $o(y_1)\le o(y_2)$.

\smallskip
{\bf (d)}  The number of strata of $\scrY$ is at most $[W_G:W_P]$.

\smallskip
{\bf (e)} We assume $\scrW$ is uni+versal, $Y$ is proper (over $k$), and the strata of $\scrY$ are quasi-affine. Then $\scrY$ has a unique closed stratum $\grs_0$ of dimension $0$ and to which all strata of $\scrY$ specialize. 

\smallskip
{\bf (f)} Suppose the hypotheses of (e) hold. Then $\scrY$ has $[W_G:W_P]$ strata if and only if the action $\dbT_{G_k,\sigma}$ has a unique closed orbit.} 

{\bf Proof:}
Let $y\in Y(k)$. Let $I_0$ be a finite subset of $I$ such that $\{\im(U_i\to Y)|i\in I_0\}$ is an open cover of $Y$. Let $i\in I_0$. Let $\grs_{iy}$ be $\bar g_i^{-1}(o(y))$ endowed with its reduced structure. It is a reduced, locally closed subscheme of $U_i$. For $y_i\in\im(U_i(k)\to Y(k))$ we have $y_i\in\im(\grs_{iy}(k)\to Y(k))$ if and only if $<\scrC_y>=<\scrC_{y_i}>$, cf. Lemma 5.1 (a). Thus for $i, j\in I_0$, we have $\grs_{iy}\times_Y U_j=\grs_{jy}\times_Y U_i$. Thus we can uniquely define a stratification $\scrY$ of $Y$ via the following property: the stratum $\grs_y$ of $\scrY$ to which $y$ belongs is the reduced, locally subscheme of $Y$ whose $k$-valued points are $\cup_{i\in I_0} \im(\grs_{iy}(k)\to Y(k))$. Thus (a) holds.

We prove (b). We first assume $\scrW$ is full. The $k$-morphism $\grs_{iy}\to o(y)$ is \'etale. Thus, as $o(y)$ is connected and smooth, each $\grs_{iy}$ is smooth and equidimensional of dimension $\dim(o(y))=\dim(G_k)-\dbS(y)$ (cf. (12) for the equality part). This implies that the stratum $\grs_y$ is smooth and equidimensional of dimension $\dim(o(y))=\dim(G_k)-\dbS(y)$. As $\dbT_{G_k,\sigma}$ has a dense, open orbit (see Example 5.6), $\scrY$ has a unique open, dense stratum. 

Let now $\scrW$ be uni+versal. Let $\tilde Y=P_k\times_k Y$ and $\tilde\scrW$ be as in Subsubsection 12.1.1. Let $\tilde{\scrY}$ be the stratification of $\tilde Y$ defined by $\tilde\scrW$, cf. (a). It is the pull back of $\scrY$ via the projection $q_2:\tilde Y\twoheadrightarrow Y$. Thus each stratum $\grs$ of $\scrY$ is the quotient of a unique stratum $\tilde{\grs}=P_k\times_k \grs$ of $\tilde\scrY$ under the left multiplication action of $P_k$ on $\tilde{\grs}$. As $\tilde{\grs}$ is smooth and equidimensional, $\grs$ is also so and its dimension is
$\dim(\tilde{\grs})-\dim(P_k)=\dim(G_k)-\dbS(y)-\dim(P_k)=\dim(N_k)-\dbS(y).$
As $\tilde{\scrY}$ has a unique dense, open stratum, the same applies to $\scrY$. Thus (b) holds.

The first part of (c) is obvious. It suffices to prove the second part of (c) for the full case, cf. the augmentation process of the previous paragraph. Thus we can assume $\scrW$ is full. Let $(U_i,\rho_i)\in\scrA$ be such that $y_2\in \im(U_i(k)\to Y(k))$ and $\bar g_i:U_i\to G_k$ is \'etale. As $o_1'\le o(y_2)$, the inverse image $\bar g_i^{*}(o_1')$ is a locally closed subscheme of $U_i$ whose image in $Y$ is contained in a stratum $\grs_1'$ of $\scrY$ that specializes to $y_2$. If $y_1'\in\grs_1'(k)$, then we have $o(y_1')=o_1'$. Thus (c) holds. Part (d) follows from Corollary 1.1. 

We prove (e) and (f). As each stratum of $\scrY$ is equidimensional, $\scrY$ has a closed stratum $\grs_0$. As $Y$ is proper, $\grs_0$ is also proper. As $\grs_0$ is proper, quasi-affine, it has dimension $0$. If $y\in\grs_0(k)$, then $\dbS(y)=\dim(N_k)$ (cf. (b)). Thus $\dim(o(y))=\dim(G_k)-\dim(N_k)=\dim(P_k)$, cf. (12) and Corollary 11.1 (b). Only the pivotal orbit in $\dbO$ has dimension $\dim(P_k)$, cf. Corollary 11.1 (c). Thus $\scrY$ has a unique closed stratum $\grs_0$ defined by: for $y\in Y(k)$, we have $y\in\grs_0$ if and only if $o_y$ is the pivotal orbit in $\dbO$. Due to the uniqueness of $\grs_0$, each stratum of $\scrY$ specializes to $\grs_0$. Thus (e) holds. As $\dbO$ has $[W_G:W_P]$ orbits (cf. Corollary 11.1 (b)), (f) follows from (c) and (e).
\endproof 

\smallskip
\noindent
{\bf Corollary 12.3.} {\it  Let $Y(k)^{\text{top}}$ be the topological space underlying $Y(k)$. 

\medskip
{\bf (a)} The map $o_Y:Y(k)^{\text{top}}\to\dbO$ that associates to $y\in Y(k)^{\text{top}}$ the orbit $o(y)\in\dbO$, is continuous.

\smallskip
{\bf (b)} If $\scrW$ is uni+versal, then the map $o_Y$ is also open.}
\medskip

{\bf Proof:}
Part (a) only translates the first part of the property 12.2 (c). To check (b), let $U_Y$ be an open subscheme of $Y$. Let $\scrW_{U_Y}:=(\scrK_{U_Y},\scrF_{U_Y},\scrV_{U_Y},\scrA_{U_Y})$ be the restriction of $\scrW=(\scrK,\scrF,\scrV,\scrA)$ to ${U_Y}$; here $\scrA_{U_Y}$ is the subset of $\scrA$ formed by those maps $(U_i,\rho_i)$ that satisfy $\im(U_i\to Y)\subseteq U_Y$ (cf. condition (iii.d) of Subsection 12.1). As $\scrW$ is uni+versal, $\scrW_{U_Y}$ is also uni+versal. By applying the second part of the property 12.2 (c) to $\scrW_{U_Y}$, we get that $o_Y(U_Y(k)^{\text{top}})=\im(o_{U_Y})$ is an open subset of $\dbO$. Thus (b) holds.
\endproof 

\smallskip
\noindent
{\bf Remark 12.4.} {\bf (a)} Let $\scrI$ be an integral canonical model in unramified mixed characteristic $(0,p)$ of a Shimura variety of Hodge type, cf. [35, Defs. 3.2.3 6) and 3.2.6]. Let $\grA_{\scrI}$ be an abelian scheme over $\scrI$ obtained naturally via an embedding of our Shimura variety of Hodge type into a Siegel modular variety. Let $Z_0$ be a smooth quotient of $\scrI$ of finite type such that $\grA_{\scrI}$ is the pull back of an abelian scheme $\grA_0$ over $Z_0$. We assume that $\scrI$ and $Z_0$ are such that there exists a reductive, closed subgroup scheme $\grG_0$ of $\pmb{GL}_{\underline{H}^1_{\text{dR}}(\grA_0/Z_0)}$ which is naturally the (integral) crystalline realization of the reductive group over $\dbQ$ that defines our Shimura variety of Hodge type. 

Let $Z$ be an open subscheme of $Z_{0W(k)}$ such that $Y:=Z_k$ is connected. Let $\grA:=\grA_0\times_{Z_0} Z$. The fact that there is a Shimura $F$-crystal $(M,\phi,G)$ over $k$ such that the vector bundle $\underline{H}^1_{\text{dR}}(\grA/Z)$ over $Z$ has a natural structure of a $G$-bundle and the condition (*) of Example 12.2 holds, is a direct consequence of the deformation theory of [35, Subsect. 5.4] and of our connectedness and reductiveness assumptions. Thus in applications, $Y$ is an open closed subscheme of $Z_{0k}$ and $\scrW$ is a uni+versal $\scrD$-bundle mod $p$ over $Y$ as in Example 12.2.

The stratification $\scrY$ of $Y$ defined by $\scrW$ is the pull back of a stratification $\scrY_0$ of a connected component of the special fibre of $Z_0$. By combining [35, Subsubsect. 3.2.12 and Thm. 5.1] with [26, Thm. (1.2)] one gets that all strata of $\scrY_0$ are quasi-affine. The tools of [42, Sect. 5] are general enough to imply (based on Corollary 1.1 and the property (d) of Subsection 12.2) that $\scrY_0$ has exactly $[W_G:W_P]$ strata (even if $Z$ is not proper over $W(k)$). 

\smallskip
{\bf (b)} The proof of Basic Theorem D is independent of the existence of an integrable and nilpotent connection $\nabla$ on $\scrK$ with respect to which $\scrF$ and $\scrV$ are horizontal. But if $\scrW$ is not defined by some Barsotti--Tate group over $Y$ endowed with crystalline tensors, the strata of $\scrY$ are not necessarily quasi-affine even if $\scrW$ is uni+versal.

\medskip
\noindent
{\bf Acknowledgement.}
 We would like to thank U of Arizona and Binghamton University for good working conditions, J.-P. Serre and O. Gabber for the proof of Lemma 2.5, and the referee for many valuable comments and suggestions.

\bigskip
\references{37}
{\nspace{

\Ref[1]
A.~Borel,
\sl Linear algebraic groups. Second edition,
\rm Grad. Texts in Math., Vol. {\bf 126}, Springer-Verlag, New York, 1991.

\Ref[2]
N.~Bourbaki,
\sl Lie groups and Lie algebras, Chapters 4--6, 
\rm Elements of Mathematics, Springer-Verlag, Berlin, 2002.

\Ref[3]
S.~Bosch, W.~L\"utkebohmert, and M.~Raynaud,
\sl N\'eron models,
\rm Ergebnisse der Mathematik und ihrer Grenzgebiete (3), Vol. {\bf 21}, Springer-Verlag, Berlin, 1990.

\Ref[4]
F.~Bruhat and J.~Tits, 
\sl Groupes r\'eductifs sur un corps local. II. Sch\'emas en groupes. Existence d'une donn\'ee radicielle valu\'ee,
\rm Inst. Hautes \'Etudes Sci. Publ. Math., Vol. {\bf 60},  5--184, 1984.

\Ref[5]
P.~Deligne,
\sl Vari\'et\'es de Shimura: interpr\'etation modulaire, et
techniques de construction de mod\`eles canoniques,
\rm Automorphic forms, representations and $L$-functions (Oregon State Univ., Corvallis, OR, 1977), Part 2,   247--289, Proc. Sympos. Pure Math., Vol. {\bf 33}, Amer. Math. Soc., Providence, RI, 1979.

\Ref[6]
M.~Demazure, A.~Grothendieck, et al. 
\sl Sch\'emas en groupes, Vol. I--III, 
\rm Lecture Notes in Math., Vol. {\bf 151--153}, Springer-Verlag, Berlin-New York, 1970.

\Ref[7]
P. Deligne and G. Lusztig, 
\sl Representations of reductive groups over finite fields,
\rm Ann. of Math. (2) {\bf 103}  (1976), no. 1, 103--161.

\Ref[8]
J.-M.~Fontaine,
\sl Groupes $p$-divisibles sur les corps locaux, 
\rm J. Ast\'erisque {\bf 47/48}, Soc. Math. de France, Paris, 1977.

\Ref[9]
W.~Fulton and J.~Harris,
\sl Representation theory. A first course,
\rm Grad. Texts in Math., Vol. {\bf 129}, Readings in Mathematics, Springer-Verlag, New York, 1991.

\Ref[10]
D.~Goss,
\sl Basic structures of function field arithmetic,
\rm Ergebnisse der Mathematik und ihrer Grenzgebiete (3), Vol. {\bf 35}, Springer-Verlag, Berlin, 1996. 

\Ref[11]
D.~Gorenstein, R.~Lyons, and R.~Solomon,
\sl The classification of the finite simple groups, Number 3,
\rm Math. Surv. and Monog., Vol. {\bf 40}, no. 3, Amer. Math. Soc., Providence, RI, 1994. 

\Ref[12]
X. ~He,
\sl G-stable pieces and partial flag varieties,
\rm to appear in Comp. Math.

\Ref[13]
J.~E.~Humphreys, 
\sl Introduction to Lie algebras and representation theory,
\rm Grad. Texts in Math., Vol. {\bf 9}, Springer-Verlag, New York-Berlin, 1972.

\Ref[14]
J.~E.~Humphreys, 
\sl Conjugacy classes in semisimple algebraic groups, 
\rm Math. Surv. and Monog., Vol. {\bf 43}, Amer. Math. Soc., Providence, RI, 1995.

\Ref[15]
L.~Illusie, 
\sl D\'eformations des groupes de Barsotti--Tate (d'apr\'es A. Grothendieck), 
\rm Seminar on arithmetic bundles: the Mordell conjecture (Paris, 1983/84),  151--198, J. Ast\'erisque {\bf 127}, Soc. Math. de France, Paris, 1985.

\Ref[16] 
R.~E.~Kottwitz, 
\sl Points on some Shimura varieties over finite fields, 
\rm J. of Amer. Math. Soc. {\bf 5} (1992), no. 2,  373--444.

\Ref[17] 
H.~Kraft, 
\sl Kommutative algebraische p-Gruppen (mit Anwendungen auf p-divisible Gruppen und abelsche Variet\"aten), 
\rm manuscript 86 pages, Univ. Bonn, 1975.

\Ref[18]
M.-A.~Knus, A.~Merkurjev, M.~Rost, and J.-P.~Tignol, 
\sl The book of involutions, 
\rm Coll. Public., Vol. {\bf 44}, Amer. Math. Soc., Providence, RI, 1998.

\Ref[19] 
R.~Langlands and M.~Rapoport,
\sl Shimuravariet\"aten und Gerben, 
\rm J. reine angew. Math. {\bf 378} (1987),  113--220. 

\Ref[20]
G.~Lusztig, 
\sl Parabolic character sheaves. I,
\rm Mosc. Math. J. {\bf 4}  (2004),  no. 1,  153--179, 311.

\Ref[21] 
H.~Matsumura, 
\sl Commutative algebra. Second edition, 
\rm Mathematics Lecture Note Series, Vol. {\bf 56}. Benjamin/Cummings Publishing Co., Inc., Reading, Mass., 1980.

\Ref[22]
B.~Moonen,
\sl Group schemes with additional structures and Weyl group cosets,
\rm Moduli of abelian varieties (Texel Island, 1999),  255--298, Progr. of Math., Vol. {\bf 195}, Birkh\"auser, Basel, 2001.

\Ref[23]
B.~Moonen,
\sl A dimension formula for Ekedahl--Oort strata,
\rm Ann. Inst. Fourier (Grenoble) {\bf 54} (2004), no. 3,  666--698.  

\Ref[24]
M.-H. Nicole and A.~Vasiu,
\sl Minimal truncations of supersingular $p$-divisible groups,
\rm Indiana Univ. Math. J. {\bf 56} (2007), no. 6, 2887--2898.

\Ref[25] 
F.~Oort, 
\sl Newton polygons and formal groups: conjectures by Manin and Grothendieck, 
\rm Ann. of Math. (2) {\bf 152} (2000), no. 1,  183--206.

\Ref[26] 
F.~Oort, 
\sl A stratification of a moduli space of abelian varieties, 
\rm Moduli of abelian varieties (Texel Island, 1999),  345--416, Progr. of Math., Vol. {\bf 195}, Birkh\"auser, Basel, 2001.

\Ref[27] 
F.~Oort, 
\sl Minimal p-divisible groups, 
\rm Ann. of Math. (2) {\bf 161} (2005), no. 2,  1021--1036.

\Ref[28]
R.~Pink,
\sl $l$-adic algebraic monodromy groups, cocharacters, and the Mumford--Tate conjecture,
\rm J. reine angew. Math. {\bf 495} (1998),  187--237.

\Ref[29] 
M.~Raynaud, 
\sl Anneaux Locaux Hens\'eliens, 
\rm Lecture Notes in Math., Vol. {\bf 169}, Springer-Verlag, Berlin-New York, 1970.

\Ref[30] 
J.-P.~Serre, 
\sl Groupes alg\'ebriques associ\'es aux modules de Hodge--Tate, 
\rm Journ\'ees de G\'eom. Alg. de Rennes (Rennes, 1978),  155--188, J. Ast\'erisque {\bf 65}, Soc. Math. de France, Paris, 1979. 

\Ref[31]
T. A. ~Springer,
\sl An extension of Bruhat's lemma,
\rm J. Algebra  {\bf 313}  (2007),  no. 1, 417--427.

\Ref[32]
R.~Steinberg,
\sl Endomorphisms of linear algebraic groups,
\rm Memoirs of the Amer. Math. Soc., Vol.  {\bf 80}, Amer. Math. Soc., Providence, RI, 1968.

\Ref[33] 
J.~Tits, 
\sl Classification of algebraic semisimple groups, 
\rm Algebraic Groups and Discontinuous Subgroups (Boulder, CO, 1965),   33--62, Proc. Sympos. Pure Math., Vol. {\bf 9}, Amer. Math. Soc., Providence, RI, 1966.

\Ref[34] 
J.~Tits,
\sl Reductive groups over local fields, 
\rm Automorphic forms, representations and $L$-functions (Oregon State Univ., Corvallis, OR, 1977), Part 1,   29--69, Proc. Sympos. Pure Math., Vol. {\bf 33}, Amer. Math. Soc., Providence, RI, 1979.

\Ref[35]
A.~Vasiu, 
\sl Integral canonical models of Shimura varieties of preabelian
 type, 
\rm Asian J. Math. {\bf 3} (1999), no. 2,  401--518.

\Ref[36]
A.~Vasiu,
\sl On two theorems for flat, affine group schemes over a discrete valuation ring,
\rm Centr. Eur. J. Math. {\bf 3} (2005),  no. 1,  14--25.

\Ref[37]
A.~Vasiu,
\sl Crystalline Boundedness Principle,
\rm Ann. Sci. \'Ecole Norm. Sup. {\bf 39} (2006), no. 2,  245--300.

\Ref[38]
A.~Vasiu,
\sl Level $m$ stratifications, 
\rm J. Alg. Geom. {\bf 17} (2008), no. 4, 599--641.

\Ref[39]
A.~Vasiu,
\sl Integral canonical models of unitary Shimura varieties,
\rm Asian J. Math. {\bf 12} (2008), no. 2, 151--176.

\Ref[40]
A.~Vasiu,
\sl Geometry of Shimura varieties of Hodge type over finite fields,
\rm Higher-dimensional geometry over finite fields,  197--243, NATO Sci. Peace Secur. Ser. D Inf. Commun. Secur., {\bf 16}, IOS, Amsterdam, 2008.

\Ref[41]
A.~Vasiu,
\sl Reconstructing p-divisible groups from their truncations of small level,
\rm Comment. Math. Helv. {\bf 85} (2010), no. 1, 165--202.

\Ref[42]
A.~Vasiu,
\sl Manin problems for Shimura varieties of Hodge type,
\rm manuscript, 

\rm http://arxiv.org/abs/math/0209410.

\Ref[43] 
T.~Wedhorn, 
\sl The dimension of Oort strata of Shimura varieties of PEL-type,
\rm Moduli of abelian varieties (Texel Island, 1999),  441--471, Progr. of Math., Vol. {\bf 195}, Birkh\"auser, Basel, 2001.

\Ref[44] 
\sl Th\'eorie des topos et cohomologie \'etale des sch\'emas. Tome 2 et 3, 
\rm S\'eminaire de G\'eométrie Alg\'ebrique du Bois-Marie 1963--1964 (SGA 4). Dirig\'e par M. Artin, A. Grothendieck, et J. L. Verdier. Avec la collaboration de N. Bourbaki, P. Deligne et B. Saint-Donat. Lecture Notes in Math., Vols. {\bf 270} and {\bf 305}, Springer-Verlag, Berlin-New York, 1972 and 1973.

}}
\noindent
\hbox{Adrian Vasiu,}
\hbox{Department of Mathematical Sciences, Binghamton University,}

\hbox{Binghamton, New York 13902-6000, U.S.A.}
\hbox{e-mail: adrian\@math.binghamton.edu}

\enddocument